\documentclass[a4paper, 10pt, USenglish, reqno]{amsart}

\PassOptionsToPackage{noend}{algpseudocode}
\PassOptionsToPackage{bookmarksdepth=4}{hyperref}
\PassOptionsToPackage{nameinlink,capitalize}{cleveref}
\PassOptionsToPackage{dvipsnames,table}{xcolor}

\usepackage{float}
\usepackage{mathrsfs}
\usepackage{tcolorbox}
\usepackage[
	eqreset=section,
	thmreset=section,
	noload={todonotes,caption},
	noqed,
]{myPreamble}
\renewcommandx{\includetikz}[3][1={},3={}]{%
	\ifstrempty{#3}{\includegraphics[#1]{Pics/#2.pdf}}{\includegraphics[#1]{Pics/#3.pdf}}%
}

\crefname{figure}{Figure}{Figures}
\Crefname{figure}{Figure}{Figures}
\crefname{ALG@line}{step}{steps}

\makeatletter
	\renewcommand\appendixproof@toc{subsection}
	\def\operator@font{\rm}
\makeatother
\addqedto[\(\triangle\)]{rem}

\crefname{algorithm}{Algorithm}{Algorithms}

	\makeatletter
		\renewcommand\alglinenumber[1]{{\footnotesize\textbf{\thealgorithm}.\oldstylenums{\arabic{ALG@line}}:}}
		\renewcommand\theALG@line{\thealgorithm.\oldstylenums{\arabic{ALG@line}}}
		\providecommand\theHALG@line{\thealgorithm.\oldstylenums{\arabic{ALG@line}}}
	\makeatother
	
	\newenvironment{boxedalgorithm}[2][]{%
		\refstepcounter{algorithm}%
		\begin{tcolorbox}[
			title={\textbf{Algorithm \thealgorithm.} #2},
			left=0pt, right=0pt, top=0pt, bottom=0pt,
			colframe = MidnightBlue,
			colback = MidnightBlue!5,
			colbacktitle = MidnightBlue!15,
			coltitle = black,
			#1,
		]
	}{%
		\end{tcolorbox}
	}

	\NewEnviron{casebox}[1]{%
		\setlength\fboxrule{1.25pt}%
		\fbox{%
			\setlength\parindent{0pt}%
			\begin{minipage}[t]{0.47\linewidth}
				\hspace*{0pt}\hfill{\sc #1}\hfill\hspace*{0pt}%
				\par\noindent
				\BODY
			\end{minipage}%
		}%
	}

	\newlist{directions}{enumerate}{1}
	\AtBeginEnvironment{directions}{%
		\let\olditem\item
		\renewcommand\item[1][]{%
			\def\itemtitle{#1}%
			\olditem
			\ifstrempty{#1}{}{%
				\phantomsection
				\addcontentsline{toc}{paragraph}{\itemtitle}%
				\textbf{#1.}~~\ignorespaces%
			}%
		}
	}
	\setlist[directions,1]{%
		label={\bf\alph*.},
		ref=\thesubsubsection\alph*,
		align=parleft,
		leftmargin=*,
	}
	\crefformat{directionsi}{#2\S#1#3}
	\Crefformat{directionsi}{#2\S#1#3}
	\crefname{directionsi}{\S}{\S}
	\Crefname{directionsi}{\S}{\S}

	\makeatletter
		\newcommand{\@DRS}[1][\gamma]{\operatorname{DRS_{#1}}}
		\newcommand{\@@DRS}[1][\nicefrac1\gamma]{\operatorname{DRS^*_{#1}}}
		\newcommand{\DRS}{\@ifstar\@@DRS\@DRS}

		\newcommand{\@DRE}[1][\gamma]{\varphi_{#1}^{\text{\sc dr}}}
		\newcommand{\@@DRE}{\psi_{\gamma_*}^{\text{\sc dr}}}
		\newcommand{\DRE}{\@ifstar\@@DRE\@DRE}
		
		\newcommand{\sDRS}{% smooth
			\mathchoice{\vphantom{\psi}}{\vphantom{\psi}}{}{}%
			\@ifstar{\smash\psi_1}{\smash\varphi_1}%
		}
		
		\newcommand{\nsDRS}{% nonsmooth
			\mathchoice{\vphantom{\psi}}{\vphantom{\psi}}{}{}%
			\@ifstar{\smash\psi_2}{\smash\varphi_2}%
		}

		\newcommand{\C}{\@ifstar\@@C\@C}
		\newcommandx{\@C}[2][1={\empty},2={\empty}]{%
			C
			\ifstrempty{#1}{}{
				\left(
					\ifx#1\empty \gamma L_{\sDRS}\else #1\fi,
					\ifx#2\empty \lambda\else #2\fi
				\right)
			}
		}
		\newcommandx{\@@C}[2][1={\empty},2={\empty}]{%
			C
			\ifstrempty{#1}{}{
				\left(
					\ifx#1\empty \nicefrac{1}{\gamma\mu_{\sDRS}}\else #1\fi,
					\ifx#2\empty \lambda\else #2\fi
				\right)
			}
		}

		\newcommand{\D}{\@ifstar\@@D\@D}
		\newcommandx{\@D}[2][1={\empty},2={\empty}]{%
			C
			\ifstrempty{#1}{}{
				\left(
					\ifx#1\empty \nicefrac{L_{\epicomp Af}}{\beta}\else #1\fi,
					\ifx#2\empty \lambda\else #2\fi
				\right)
			}
		}
		\newcommandx{\@@D}[2][1={\empty},2={\empty}]{%
			C
			\ifstrempty{#1}{}{
				\left(
					\ifx#1\empty \nicefrac{\beta\|A\|^2}{\mu_f}\else #1\fi,
					\ifx#2\empty \lambda\else #2\fi
				\right)
			}
		}
	\makeatother
	\newcommand{\FBE}[1][\gamma]{\varphi_{#1}^{\text{\sc fb}}}

	\DeclareMathOperator{\refl}{\mathcal R}											% Rf(x) = f(-x)
	\DeclareMathOperator{\infconv}{\square}											% Rf(x) = f(-x)
	\newcommandx{\Res}[2][1={\gamma},2={dr}]{R_{#1}^{\text{\sc #2}}}%
	
	\renewcommand{\L}{\mathcal L}

	\newcommand{\ADMM}[1][\beta]{\operatorname{ADMM_{#1}}}

	\newcommand{\sADMM}{f}% smooth
	\newcommand{\nsADMM}{g}% nonsmooth
	\newcommand{\LL}{\mathscr L}% Lagrangian

	\newcommand{\epicomp}[2]{(#1#2)}

	\renewcommand{\FBsmooth}{\sDRS}
	\renewcommand{\FBnonsmooth}{\nsDRS}

	\usetikzlibrary{intersections}

	\definecolor{DRS}{rgb}{0,0.6056,0.9787}% blue
	\definecolor{Nesterov}{rgb}{0.7644,0.4441,0.8243}% purple
	\definecolor{LBFGS}{rgb}{0.8889,0.4356,0.2781}% red
	\colorlet{Anderson}{Nesterov!50!white}
	\definecolor{Broyden}{rgb}{0.2422,0.6433,0.3044}% green

	\pgfplotsset{%
		compat = newest,
	}
	\pgfplotsset{%
		mylegend/.style = {
			legend style = {
				solid,
				rounded corners,
				color        = {rgb,1:red,0.0;green,0.0;blue,0.0},
				line width   = 1,
				draw opacity = 1.0,
				fill         = {rgb,1:red,1.0;green,1.0;blue,1.0},
				fill opacity = 0.6,
				text opacity = 1.0,
				font         = {{\fontsize{8 pt}{10.4 pt}\selectfont}},
				text         = {rgb,1:red,0.0;green,0.0;blue,0.0},
				cells        = {anchor = center},
				at           = {(1.02, 1)},
			},
		},
		every axis/.append style = {mylegend},% <--- automatically adds it to every plot (overwrites previous legend settings)
	}
	\pgfplotsset{%
		consensus/.style = {% shared by all consensus plots
			height = 43mm,
			width  = 75mm,
			xlabel shift = {-3pt},
			point meta max = nan,
			point meta min = nan,
			axis background/.style = {
				fill    = {rgb,1:red,1.0;green,1.0;blue,1.0},
				opacity = 1.0,
			},
			scaled x ticks = false,
			x tick style = {
				color   = {rgb,1:red,0.0;green,0.0;blue,0.0},
				opacity = 1.0,
			},
			x tick label style = {
				color   = {rgb,1:red,0.0;green,0.0;blue,0.0},
				opacity = 1.0,
				rotate  = 0,
			},
			xlabel style = {
				at           = {(ticklabel cs:0.5)},
				anchor       = near ticklabel,
				font         = {{\fontsize{11 pt}{14.3 pt}\selectfont}},
				color        = {rgb,1:red,0.0;green,0.0;blue,0.0},
				draw opacity = 1.0,
				rotate       = 0,
			},
			xmajorgrids = true,
			xtick align = inside,
			xticklabel style = {
				draw opacity = 1.0,
				rotate       = 0,
				font         = {{\fontsize{8 pt}{10.4 pt}\selectfont}},
				color        = {rgb,1:red,0.0;green,0.0;blue,0.0},
			},
			x grid style = {
				solid,
				color        = {rgb,1:red,0.0;green,0.0;blue,0.0},
				draw opacity = 0.1,
				line width   = 0.5,
			},
			axis x line* = left,
			x axis line style = {
				solid,
				color        = {rgb,1:red,0.0;green,0.0;blue,0.0},
				draw opacity = 1.0,
				line width   = 1,
			},
			scaled y ticks = false,
			y tick style = {
				color = {rgb,1:red,0.0;green,0.0;blue,0.0},
				opacity = 1.0,
			},
			y tick label style = {
				color   = {rgb,1:red,0.0;green,0.0;blue,0.0},
				opacity = 1.0,
				rotate  = 0,
			},
			ylabel style = {
				at           = {(ticklabel cs:0.5)},
				anchor       = near ticklabel,
				font         = {{\fontsize{11 pt}{14.3 pt}\selectfont}},
				color        = {rgb,1:red,0.0;green,0.0;blue,0.0},
				draw opacity = 1.0,
				rotate       = 0,
			},
			ymajorgrids = true,
			ytick align = inside,
			yticklabel style = {
				font         = {{\fontsize{8 pt}{10.4 pt}\selectfont}},
				color        = {rgb,1:red,0.0;green,0.0;blue,0.0},
				draw opacity = 1.0,
				rotate       = 0,
			},
			y grid style = {
				solid,
				color        = {rgb,1:red,0.0;green,0.0;blue,0.0},
				draw opacity = 0.1,
				line width   = 0.5,
			},
			axis y line*      = left,
			y axis line style = {
				solid,
				color        = {rgb,1:red,0.0;green,0.0;blue,0.0},
				draw opacity = 1.0,
				line width   = 1,
			},
			colorbar = false
		},
		consensus_res/.style = {% shared by all consensus RESIDUAL plots
			consensus,
			log basis y = 10,
			legend style/.append style = {anchor = north west},
			ylabel = {$\beta\|x-z\|$},
			xshift =  1.0mm,
			yshift = -1.0mm,
		},
		consensus_cost/.style = {% shared by all consensus COST plots
			consensus,
			ylabel = {$f(z)$},
			legend style/.append style = {anchor = north east},
			xshift      = 74.7mm,
			yshift      = -1.0mm,
		},
		20newsgroup/.style = {% shared by 20newsgroup plots
			xmin =  -1.97,
			xmax = 102.97,
			xtick       = {0.0,25.0,50.0,75.0,100.0},
			xticklabels = {$0$,$25$,$50$,$75$,$100$},
		},
		20newsgroup_res/.style = {% 20newsgroup residual plots (manually add ymode = log)
			20newsgroup,
			consensus_res,
			ymin = 0.000577,
			ymax = 2728.141,
			ytick       = {0.01,1.0,100.0},
			yticklabels = {$10^{-2}$,$10^0$,$10^2$},
		},
		20newsgroup_cost/.style = {% 20newsgroup cost plots
			20newsgroup,
			consensus_cost,
			ymin = -104.702,
			ymax =  -98.204,
			ytick       = {-104.0,-103.0,-102.0,-101.0,-100.0,-99.0},
			yticklabels = {$-104$,$-103$,$-102$,$-101$,$-100$,$-99$},
			ylabel shift = {-3pt},
		},
		nips/.style = {% shared by nips plots
			xmin = -13.97,
			xmax = 514.97,
			xtick       = {0.0,100.0,200.0,300.0,400.0,500.0},
			xticklabels = {$0$,$100$,$200$,$300$,$400$,$500$},
		},
		nips_res/.style = {% nips residual plots (manually add ymode = log)
			nips,
			consensus_res,
			xmin = -13.97,
			xmax = 514.97,
			xtick       = {0.0,100.0,200.0,300.0,400.0,500.0},
			xticklabels = {$0$,$100$,$200$,$300$,$400$,$500$},
		},
		nips_cost/.style = {% nips cost plots
			nips,
			consensus_cost,
			ymin = -9.393,
			ymax = -9.378,
			ytick       = {-9.39,-9.385,-9.38},
			yticklabels = {$-9.39$,{},$-9.38$},
			ylabel shift = {-10pt},
		},
	}
	\pgfplotsset{%
		colormap = {contourColormap}{[1pt]%
			rgb(0pt)=(0.2422,0.1504,0.6603);%
			rgb(1pt)=(0.25039,0.164995,0.707614);%
			rgb(2pt)=(0.257771,0.181781,0.751138);%
			rgb(3pt)=(0.264729,0.197757,0.795214);%
			rgb(4pt)=(0.270648,0.214676,0.836371);%
			rgb(5pt)=(0.275114,0.234238,0.870986);%
			rgb(6pt)=(0.2783,0.255871,0.899071);%
			rgb(7pt)=(0.280333,0.278233,0.9221);%
			rgb(8pt)=(0.281338,0.300595,0.941376);%
			rgb(9pt)=(0.281014,0.322757,0.957886);%
			rgb(10pt)=(0.279467,0.344671,0.971676);%
			rgb(11pt)=(0.275971,0.366681,0.982905);%
			rgb(12pt)=(0.269914,0.3892,0.9906);%
			rgb(13pt)=(0.260243,0.412329,0.995157);%
			rgb(14pt)=(0.244033,0.435833,0.998833);%
			rgb(15pt)=(0.220643,0.460257,0.997286);%
			rgb(16pt)=(0.196333,0.484719,0.989152);%
			rgb(17pt)=(0.183405,0.507371,0.979795);%
			rgb(18pt)=(0.178643,0.528857,0.968157);%
			rgb(19pt)=(0.176438,0.549905,0.952019);%
			rgb(20pt)=(0.168743,0.570262,0.935871);%
			rgb(21pt)=(0.154,0.5902,0.9218);%
			rgb(22pt)=(0.146029,0.609119,0.907857);%
			rgb(23pt)=(0.138024,0.627629,0.89729);%
			rgb(24pt)=(0.124814,0.645929,0.888343);%
			rgb(25pt)=(0.111252,0.6635,0.876314);%
			rgb(26pt)=(0.0952095,0.679829,0.859781);%
			rgb(27pt)=(0.0688714,0.694771,0.839357);%
			rgb(28pt)=(0.0296667,0.708167,0.816333);%
			rgb(29pt)=(0.00357143,0.720267,0.7917);%
			rgb(30pt)=(0.00665714,0.731214,0.766014);%
			rgb(31pt)=(0.0433286,0.741095,0.73941);%
			rgb(32pt)=(0.0963952,0.75,0.712038);%
			rgb(33pt)=(0.140771,0.7584,0.684157);%
			rgb(34pt)=(0.1717,0.766962,0.655443);%
			rgb(35pt)=(0.193767,0.775767,0.6251);%
			rgb(36pt)=(0.216086,0.7843,0.5923);%
			rgb(37pt)=(0.246957,0.791795,0.556743);%
			rgb(38pt)=(0.290614,0.79729,0.518829);%
			rgb(39pt)=(0.340643,0.8008,0.478857);%
			rgb(40pt)=(0.3909,0.802871,0.435448);%
			rgb(41pt)=(0.445629,0.802419,0.390919);%
			rgb(42pt)=(0.5044,0.7993,0.348);%
			rgb(43pt)=(0.561562,0.794233,0.304481);%
			rgb(44pt)=(0.617395,0.787619,0.261238);%
			rgb(45pt)=(0.671986,0.779271,0.2227);%
			rgb(46pt)=(0.7242,0.769843,0.191029);%
			rgb(47pt)=(0.773833,0.759805,0.16461);%
			rgb(48pt)=(0.820314,0.749814,0.153529);%
			rgb(49pt)=(0.863433,0.7406,0.159633);%
			rgb(50pt)=(0.903543,0.733029,0.177414);%
			rgb(51pt)=(0.939257,0.728786,0.209957);%
			rgb(52pt)=(0.972757,0.729771,0.239443);%
			rgb(53pt)=(0.995648,0.743371,0.237148);%
			rgb(54pt)=(0.996986,0.765857,0.219943);%
			rgb(55pt)=(0.995205,0.789252,0.202762);%
			rgb(56pt)=(0.9892,0.813567,0.188533);%
			rgb(57pt)=(0.978629,0.838629,0.176557);%
			rgb(58pt)=(0.967648,0.8639,0.16429);%
			rgb(59pt)=(0.96101,0.889019,0.153676);%
			rgb(60pt)=(0.959671,0.913457,0.142257);%
			rgb(61pt)=(0.962795,0.937338,0.12651);%
			rgb(62pt)=(0.969114,0.960629,0.106362);%
			rgb(63pt)=(0.9769,0.9839,0.0805)%
		}%
	}
\showgrayfalse

\def\TheTitle{Douglas-Rachford splitting and ADMM for nonconvex optimization}
\def\TheSubTitle{Accelerated and Newton-type linesearch algorithms}
\def\TheShortTitle{DRS and ADMM for nonconvex optimization: Accelerated and Newton-type algorithms}

\def\TheFunding{%
	A. Themelis is supported by the Japan Society for the Promotion of Science (JSPS) KAKENHI grant JP21K17710.
\\
	L. Stella's work was done prior to joining Amazon.
\\
	P. Patrinos is supported by the
	Research Foundation Flanders (FWO) research projects G0A0920N, G086518N, and G086318N;
	Research Council KU Leuven C1 project No. C14/18/068;
	Fonds de la Recherche Scientifique -- FNRS and the Fonds Wetenschappelijk Onderzoek -- Vlaanderen under EOS project 30468160 (SeLMA);
	European Union's Horizon 2020 research and innovation programme under the Marie Sk\l odowska-Curie grant agreement No. 953348.%
}
\def\TheKeywords{%
	Nonsmooth nonconvex optimization,
	Douglas-Rachford splitting,
	ADMM,
	quasi-Newton methods%
}
\def\TheSubjclass{%
	90C06, %(large scale),
	90C25, %(cvx prog),
	90C26, %(noncvx prog),
	49J52, %(CoV: nonsm. analysis),
	49J53.%(CoV: set-valued \& var.analysis)%
}

\edef\TheFullTitle{\TheTitle: \TheSubTitle}
\newcommand{\Title}[1]{\expandafter\expandafter\expandafter\title\expandafter[\TheShortTitle]{#1}}
\expandafter\expandafter\expandafter\Title\expandafter{\TheFullTitle}
\author{%
	Andreas Themelis,
	Lorenzo Stella
	and
	Panagiotis Patrinos%
}
\thanks{%
	\TheFunding\\
	{\tt
		\href{mailto:andreas.themelis@kuleuven.be}{andreas.themelis@kuleuven.be},\hfill
		\href{mailto:panos.patrinos@esat.kuleuven.be}{panos.patrinos@esat.kuleuven.be},\hfill
		\href{mailto:lorenzostella@gmail.com}{lorenzostella@gmail.com}%
	}%
}
\subjclass{\TheSubjclass}

\begin{document}

	\begin{abstract}

		Although the performance of popular optimization algorithms such as Douglas-Rachford splitting (DRS) and the ADMM is satisfactory in small and well-scaled problems, ill conditioning and problem size pose a severe obstacle to their reliable employment.
		Expanding on recent convergence results for DRS and ADMM applied to nonconvex problems, we propose two linesearch algorithms to enhance and robustify these methods by means of quasi-Newton directions.
		The proposed algorithms are suited for nonconvex problems, require the same black-box oracle of DRS and ADMM, and maintain their (subsequential) convergence properties.
		Numerical evidence shows that the employment of L-BFGS in the proposed framework greatly improves convergence of DRS and ADMM, making them robust to ill conditioning.
		Under regularity and nondegeneracy assumptions at the limit point, superlinear convergence is shown when quasi-Newton Broyden directions are adopted.
%%<< [APPENDED to input], level=-1, file=TeX/Text/Abstract.tex
	\end{abstract}
	\maketitle

	% ~~~~~~~~~~~~~~~~~~~~~~~~~~~~~~~~~~~~~~~~~~~~~~~~~~~~~~~~~~~~~~~~~~~~~~~~~~~ %

	\section{Introduction}

		Due to their simplicity and versatility, the \emph{Douglas-Rachford splitting} (DRS) and the \emph{alternating direction method of multipliers} (ADMM) have gained much popularity in the last decades.
		Although originally designed for convex problems, their generalizations and extensions to nonconvex problems have recently attracted much attention, see e.g.
		\cite{hesse2013nonconvex,bauschke2014local,bauschke2014method,hesse2014alternating,li2016douglas,li2017peaceman,themelis2020douglas} for DRS and
		\cite{li2015global,hong2016convergence,guo2017convergence,goncalves2019convergence,wang2019global,themelis2020douglas} for ADMM.
		The former algorithm addresses the following composite minimization problems
		\begin{equation}\label{eq:P}
			\minimize_{s\in\R^p} \varphi(s)\equiv\sDRS(s)+\nsDRS(s),
		\end{equation}
		for some \(\func{\sDRS,\nsDRS}{\R^p}{\Rinf}\) (\(\Rinf\coloneqq\R\cup\set\infty\) denotes the extended-real line).
		Starting from some \(s\in\R^p\), one iteration of DRS applied to \eqref{eq:P} with \emph{stepsize} \(\gamma>0\) and \emph{relaxation} \(\lambda>0\) amounts to
		\[
			\tag{DRS}\label{DRS}
			\renewcommand\arraystretch{1.3}
			\begin{cases}[lcl]
				u
			& {}\in{} &
				\prox_{\gamma\sDRS}(s)
			\\
				v
			& {}\in{} &
				\prox_{\gamma\nsDRS}(2u-s)
			\\
				s^+
			& {}={} &
				s+\lambda(v-u).
			\end{cases}
		\]
		Here, \(\prox_h\) denotes the \DEF{proximal mapping} of function \(h\); cf. \cref{sec:Notation}.
		
		The ADMM addresses optimization problems that can be formulated as
		\begin{equation}\label{eq:CP}
			\minimize_{(x,z)\in\R^m\times\R^n}{
				f(x)+g(z)
			}
		\quad
			\stt Ax+Bz=b,
		\end{equation}
		for some \(\func{f}{\R^m}{\Rinf}\), \(\func{g}{\R^m}{\Rinf}\), \(A\in\R^{p\times m}\), \(B\in\R^{p\times n}\) and \(b\in\R^p\).
		Equivalently, problem \eqref{eq:CP} amounts to the minimization of \(\func{\Phi}{\R^m\times\R^n}{\Rinf}\) given by
		\[
			\Phi(x,z)
		{}\coloneqq{}
			\begin{cases}[l@{~~}l]
				f(x)+g(z) & \text{if \(Ax+Bz=b\),}
			\\
				\infty & \text{otherwise.}
			\end{cases}
		\]
		Starting from a triplet \((x,y,z)\in\R^m\times\R^p\times\R^n\), one iteration of ADMM applied to \eqref{eq:CP} with \emph{penalty} \(\beta>0\) and \emph{relaxation} \(\lambda>0\) amounts to
		\[\tag{ADMM}\label{ADMM}
			\renewcommand\arraystretch{1.3}
			\begin{cases}[lcl]
				y^{\nicefrac+2}
			&{}={}&
				y - \beta(1-\lambda)(Ax+Bz-b)
			\\
				x^+
			&{}\in{}&
				\argmin\LL_\beta({}\cdot{},y^{\nicefrac+2},z)
			\\
				y^+
			&{}={}&
				y^{\nicefrac+2} + \beta(Ax^++Bz-b)
			\\
				z^+
			&{}\in{}&
				\argmin\LL_\beta(x^+,{}\cdot{},y^+),
			\end{cases}
		\]
		where \(\LL_\beta\) is the \(\beta\)-augmented Lagrangian of \eqref{eq:CP}, namely
		\begin{align}\label{eq:L}
			\LL_\beta(x,z,y)
		{}\coloneqq{} &
			f(x)+g(z)+\innprod{y}{Ax+Bz-b}+\tfrac\beta2\|Ax+Bz-b\|^2.
		\end{align}
		
		Although apparently more general, \ref{ADMM} is known to be equivalent to \ref{DRS} applied to the Fenchel dual of \eqref{eq:P} when the problem is convex; the identity of the two algorithms has been recently shown to hold in general through a simple change of variable \cite{yan2016self,themelis2020douglas} (see \cref{thm:ADMM_DRS} for the details).
		
		% \begin{rem}[Classical (un-relaxed) versions]
			In both algorithms, the relaxation \(\lambda\) serves as an averaging factor, and \(\lambda=1\) corresponds to the standard form in which the methods are most known.
			In the particular case of \ref{ADMM}, when \(\lambda=1\) one has \(y^{\nicefrac+2}=y\); the intermediate variable \(y^{\nicefrac+2}\) can thus be removed from the formulation in this case, and one recovers the familiar 3-step update of the algorithm.
			The half-update \(y^{\nicefrac+2}\) for the general case is introduced for simplification purposes; other sources, including the pioneering work \cite{eckstein1992douglas} (see Theorem 8 therein), avoid this by replacing \(y^{\nicefrac+2}\) with \(y\) and all occurrences of \(Ax^+\) in the \(z\)- and \(y\)-updates with \(\lambda Ax^+-(1-\lambda)(Bz-b)\).
		% \end{rem}
		Independently of \(\lambda\), the core of the respective iterations can be summarized in the following familiar oracles:
		\begin{align}
		\label{eq:DRS}
			\DRS(\tilde s)
		{}={} &
			\set{(u,v)\in\R^p\times\R^p}[~
				\begin{array}{l}
					u\in\prox_{\gamma\sDRS}(\tilde s)\\
					v\in\prox_{\gamma\nsDRS}(2u-\tilde s)
				\end{array}
			]
		\shortintertext{and}
		\label{eq:ADMM}
			\ADMM(\tilde y,\tilde z)
		{}={} &
			\set{(x,y,z)\in\R^m\times\R^p\times\R^n}[~
				\begin{array}{l}
					x\in\argmin\LL_\beta({}\cdot{},\tilde y,\tilde z)\\
					y=\tilde y+\beta(Ax+B\tilde z-b)\\
					z\in\argmin\LL_\beta(x,y,{}\cdot{})
				\end{array}
			].
		\end{align}
		
		One of the main advantages of \ref{DRS} and \ref{ADMM} lies in their ``splitting'' nature, in the sense that they exploit the additive structure of the respective problems \eqref{eq:P} and \eqref{eq:CP} by performing operations involving only either one component.
		For this reason, these methods typically involve simple operations and are thus amenable to address large-scale problems.
		On the other hand, although every iteration is relatively cheap, whenever the problem is not well scaled convergence up to a satisfactory tolerance may require prohibitively many iterations.

		\subsection{Contributions}

			Aware of these pros an cons and in contrast to an \emph{open-loop} preconditioining approach, in this paper we propose linesearch variants that allow us to integrate \ref{DRS} and \ref{ADMM} with fast update directions, stemming for instance from quasi-Newton schemes.
			To the best of our knowledge, the resulting algorithms are the first that (1) are compatible with \emph{fully nonconvex} problems, (2) maintain the \emph{same complexity} of the original \ref{DRS} and \ref{ADMM} iterations, and (3) preserve their \emph{global (subsequential) convergence guarantees}.
			Moreover, under regularity and nondegeneracy assumptions at the limit point, they converge superlinearly when a modified Broyden's scheme is used to compute the update directions.
			Extensive numerical simulations show that limited-memory quasi-Newton methods such as L-BFGS are also very effective in practice.

		\subsection{Paper organization}

			A preliminary discussion on the methodology is given in the following subsection, where the proposed linesearch \cref{alg:DRS,alg:ADMM} are also presented.
			The section will then conclude with a list of notational conventions adopted throughout the paper.
			\Cref{sec:Blocks} recaps the key properties of the \ref{DRS} and \ref{ADMM} oracles, which constitute the building blocks of the proposed algorithms.
			Once these preliminaries are dealt with, \cref{sec:Algorithms} offers a thorough overview of the proposed algorithms, together with a list of possible choices for the update directions that can considerably improve convergence speed.
			Parameter-free variants of the algorithms for out-of-the-box implementations are also provided in the concluding subsection.
			\Cref{sec:Convergence} contains the convergence results of the two algorithms, with some auxiliary material deferred to \cref{sec:Appendix}.
			In \cref{sec:Simulations} we provide numerical evidence in support of the efficacy of the proposed algorithms with simulations on nonconvex sparse linear regression and sparse principal component analysis problems, and on a strongly convex model predictive control problem.
			\Cref{sec:Conclusions} concludes the paper.

		\subsection{Methodology overview and proposed linesearch algorithms}\label{sec:AlgIntro}

			\begin{figure}[b]
				\def\vhrulesep{3pt}%
				\AtBeginEnvironment{algorithmic}{%
					\linespread{1.424}\selectfont
				}%
				\normalsize
				\begin{boxedalgorithm}{Linesearch DRS}
					\phantomsection
					\label[algorithm]{alg:DRS}

			\begin{algorithmic}[1]
			\item[]\hspace*{-\algorithmicindent}%
				\begin{minipage}{\linewidth+\algorithmicindent}
					{\sc Require}~
			% 	\\
			% 	\noindent
					\(s^0\in\R^p\);~
					tol. \(\varepsilon>0\);~
					relax. \(\lambda\in(0,2)\);~
					max \# backtracks \(i_{\rm max}\leq\infty\)%
				\\
				\noindent
					\begin{casebox}{under \cref{ass:DRS} {\normalfont(set \(\pi\coloneqq1\))}}
						stepsize \(\gamma\) as in \eqref{eq:DRE:gamma}%
					\\
						constant \smash{\(0<c<\C\)} as in \eqref{eq:DRE:C}
					\end{casebox}
					\hfill
					\begin{casebox}{under \smash{\cref{ass:DRS*}} {\normalfont(set \(\pi\coloneqq-1\))}}
						stepsize \(\gamma>\nicefrac{1}{\mu_{\sDRS}}\)%
					\\
						constant \smash{\(0<c<\C*\)} as in \eqref{eq:DRE:C}
					\end{casebox}
				\end{minipage}
			% ~~~~~~~~~~~~~~~~~~~~~~~~~~~~~~~~~~~~~~~~~~~~~~~~~~~~~~~~~~~~~~~~~~~~~~~~~~~~~~~~~~~~~~~~~~~~~~~~~
			% \vspace*{\vhrulesep}
			% \hrule
			% \vspace*{\vhrulesep}
			% ~~~~~~~~~~~~~~~~~~~~~~~~~~~~~~~~~~~~~~~~~~~~~~~~~~~~~~~~~~~~~~~~~~~~~~~~~~~~~~~~~~~~~~~~~~~~~~~~~
			% \item[]\hspace*{-\algorithmicindent}%
			% {\color{orange}%
			% 	{\sc Provide}~
			% 	\(x_*\) satisfying
			% 	\(
			% 		\dist(0,\hat\partial\varphi(x_*))
			% 		{}\leq{}
			% 		\varepsilon
			% 	\)
			% 	{\color{red}\bf under Assumption I}%\cref{ass:DRS}}
			% }%
			% ~~~~~~~~~~~~~~~~~~~~~~~~~~~~~~~~~~~~~~~~~~~~~~~~~~~~~~~~~~~~~~~~~~~~~~~~~~~~~~~~~~~~~~~~~~~~~~~~~
			\vspace*{\vhrulesep}
			\hrule
			\vspace*{\vhrulesep}
			% ~~~~~~~~~~~~~~~~~~~~~~~~~~~~~~~~~~~~~~~~~~~~~~~~~~~~~~~~~~~~~~~~~~~~~~~~~~~~~~~~~~~~~~~~~~~~~~~~~
			% \addtocounter{ALG@line}{-1}%
			\item[]\hspace*{-\algorithmicindent}%
				\begin{minipage}{\linewidth+\algorithmicindent}
					{\sc Initialize}~
			% 	\\
					\(k=0\);~
					compute
					\(
						(u^0,v^0)
					{}\in{}
						\DRS(s^0)
					\)~
					and evaluate \(\DRE(s^0)\) as in \eqref{eq:DRE}%
				\end{minipage}
			% ~~~~~~~~~~~~~~~~~~~~~~~~~~~~~~~~~~~~~~~~~~~~~~~~~~~~~~~~~~~~~~~~~~~~~~~~~~~~~~~~~~~~~~~~~~~~~~~~~
			\hrule
			\vspace*{\vhrulesep}
			\vspace*{\vhrulesep}
			% ~~~~~~~~~~~~~~~~~~~~~~~~~~~~~~~~~~~~~~~~~~~~~~~~~~~~~~~~~~~~~~~~~~~~~~~~~~~~~~~~~~~~~~~~~~~~~~~~~
			\State\label{step:DRS:1}%
				\(r^k=u^k-v^k\);~~
				{\sc if}~~\(\frac1\gamma\|r^k\|\leq\varepsilon\)~~{\sc then~ return} \((s^k,u^k,v^k)\);%
			% 	\(x_*\coloneqq v^k\)
				~~{\sc end if}%
			%%%%
			\State\label{step:DRS:nominal}%
				Set \(\bar s^{k+1}=s^k-\lambda r^k\)
			\Comment{Nominal DRS step}
			\State
				Select an update direction \(d^k\in\R^p\) and set \(\tau_k=1\) and \(i_k=0\)
			\State\label{step:DRS:1/2}%
				Define the candidate update as
				\(
					s^{k+1}
				{}={}
					(1-\tau_k)\bar s^{k+1} + \tau_k(s^k+d^k)
				\)
			\State\label{step:DRS:uv}%
				Compute
				\(
					(u^{k+1},v^{k+1})
				{}\in{}
					\DRS(s^{k+1})
				\)
				and use it to evaluate
				\(
					\DRE(s^{k+1})
				\)
				as in \eqref{eq:DRE}
			%%%%
			\If{~~\(\pi\DRE(s^{k+1})\leq\pi\DRE(s^k)-\tfrac c\gamma\|r^k\|^2\)~~}\label{step:DRS:LS}%
				\Comment{Update accepted}
				\Statex\hspace*{\algorithmicindent}
					\(k\gets k+1\)~ and go to \cref{step:DRS:1}
			\ElsIf{~~\(i_k=i_{\rm max}\)~~}
				\Comment{Max \#backtrackings: do nominal DRS step}
				\Statex\hspace*{\algorithmicindent}
					Set \(s^{k+1}=\bar s^{k+1}\) and
					\(
						(u^{k+1},v^{k+1})
					{}\in{}
						\DRS(s^{k+1})
					\)
				\Statex\hspace*{\algorithmicindent}
					\(k\gets k+1\)~ and go to \cref{step:DRS:1}%
			\Else
				\Comment{Backtrack and retry}
				\Statex\hspace*{\algorithmicindent}\label{step:DRS:tau/2}
					\(\tau_k\gets\nicefrac{\tau_k}{2}\),~
					\(i_k\gets i_k+1\),~
					and restart from \cref{step:DRS:1/2}
			\EndIf
			\end{algorithmic}
				\end{boxedalgorithm}
			\end{figure}	
			\begin{figure}[t]
				\def\vhrulesep{3pt}%
				\AtBeginEnvironment{algorithmic}{%
					\linespread{1.424}\selectfont
				}%
				\normalsize
				\begin{boxedalgorithm}{Linesearch ADMM}
					\phantomsection
					\label[algorithm]{alg:ADMM}

			\begin{algorithmic}[1]
			\item[]\hspace*{-\algorithmicindent}%
				\begin{minipage}{\linewidth+\algorithmicindent}
					{\sc Require}~
			% 	\\
			% 	\noindent
					\((x^{-1},y^{-1},z^{-1})\in\R^n\times\R^p\times\R^m\);~
					\(\varepsilon>0\);~
					\(\lambda\in(0,2)\);~
					\(i_{\rm max}\leq\infty\)%
				\\
				\noindent
					\begin{casebox}{under \cref{ass:ADMM} {\normalfont(set \(\pi\coloneqq1\))}}
						penalty \(\beta\) as in \eqref{eq:ADMM:beta}
					\\
						constant \smash{\(0<c<\D\)} as in \eqref{eq:DRE:C}
					\end{casebox}
					\hfill
					\begin{casebox}{under \smash{\cref{ass:ADMM*}} {\normalfont(set \(\pi\coloneqq-1\))}}
						penalty \(\beta<\nicefrac{\mu_f}{\|A\|^2}\)%
					\\
						constant \smash{\(\mathtight[0.55]0<c<\D*\)} as in \eqref{eq:DRE:C}%
					\end{casebox}
				\end{minipage}
			% ~~~~~~~~~~~~~~~~~~~~~~~~~~~~~~~~~~~~~~~~~~~~~~~~~~~~~~~~~~~~~~~~~~~~~~~~~~~~~~~~~~~~~~~~~~~~~~~~~
			\vspace*{\vhrulesep}
			\hrule
			\vspace*{\vhrulesep}
			% ~~~~~~~~~~~~~~~~~~~~~~~~~~~~~~~~~~~~~~~~~~~~~~~~~~~~~~~~~~~~~~~~~~~~~~~~~~~~~~~~~~~~~~~~~~~~~~~~~
			% \item[]\hspace*{-\algorithmicindent}%
			% 	{\sc Provide}
			% 	\(x_*,y_*,z_*\):
			% 	\(\|Ax_*+Bz_*-b\|\leq\frac\varepsilon\beta\),
			% 	\(-\trans Ay_*\in\hat\partial\sADMM(x_*)\),
			% 	\(\dist(-\trans By_*,\hat\partial\nsADMM(z_*))\leq\|B\|\varepsilon\)%
			% ~~~~~~~~~~~~~~~~~~~~~~~~~~~~~~~~~~~~~~~~~~~~~~~~~~~~~~~~~~~~~~~~~~~~~~~~~~~~~~~~~~~~~~~~~~~~~~~~~
			% \vspace*{\vhrulesep}
			% \hrule
			% ~~~~~~~~~~~~~~~~~~~~~~~~~~~~~~~~~~~~~~~~~~~~~~~~~~~~~~~~~~~~~~~~~~~~~~~~~~~~~~~~~~~~~~~~~~~~~~~~~
			% \addtocounter{ALG@line}{-1}%
			\item[]\hspace*{-\algorithmicindent}%
				{\sc Initialize}~
				\begin{tabular}[t]{@{}l@{}}
					\(k=0\);~
					\(r^{-1}=Ax^{-1}+Bz^{-1}-b\);~
					\(
						y^{-\nicefrac12}
					{}={}
						y^{-1}-\beta(1-\lambda)r^{-1}
					\);%
				\\
					\(
						(x^0,y^0,z^0)
					{}\in{}
						\ADMM(y^{-\nicefrac12},z^{-1})
					\)~
					and evaluate \(\LL_\beta(x^0,y^0,z^0)\)
				\end{tabular}
			% ~~~~~~~~~~~~~~~~~~~~~~~~~~~~~~~~~~~~~~~~~~~~~~~~~~~~~~~~~~~~~~~~~~~~~~~~~~~~~~~~~~~~~~~~~~~~~~~~~
			\hrule
			\vspace*{\vhrulesep}
			\vspace*{\vhrulesep}
			% ~~~~~~~~~~~~~~~~~~~~~~~~~~~~~~~~~~~~~~~~~~~~~~~~~~~~~~~~~~~~~~~~~~~~~~~~~~~~~~~~~~~~~~~~~~~~~~~~~
			\State\label{step:ADMM:1}%
				\(r^k=Ax^k+Bz^k-b\);~~
				{\sc if}~~\(\beta\|r^k\|\leq\varepsilon\)~~{\sc then~ return} \((x^k,y^k,z^k)\);%
			% 	\((x_*,y_*,z_*)=(x^k,y^k,z^k)\)
				~~{\sc end if}%
			%%%%
			\State\label{step:ADMM:nominal}%
				Set \(\bar y^{k+\nicefrac12}=y^k-\beta(1-\lambda)r^k\)
			\State
				Select an update direction \(d^k\in\R^p\) and set \(\tau_k=1\) and \(i_k=0\)
			\State\label{step:ADMM:1/2}%
				Define the candidate update
				\(
					y^{k+\nicefrac12}
				{}={}
					(1-\tau_k)\bar y^{k+\nicefrac12}
					{}+{}
					\tau_k(y^k-\beta(r^k+d^k))
				\)
			\State
				Compute
				\(
					(x^{k+1}\!,y^{k+1}\!,z^{k+1})
				{}\in{}
					\ADMM(y^{k+\nicefrac12}\!,z^k)
				\)
				and evaluate \(\LL_\beta(x^{k+1}\!,y^{k+1}\!,z^{k+1})\)%
			%%%%
			\If{~~\(\pi\LL_\beta(x^{k+1},z^{k+1},y^{k+1})\leq\pi\LL_\beta(x^k,z^k,y^k)-\beta c\|r^k\|^2\)~~}\label{step:ADMM:LS}%
				\Comment{Update accepted}
				\Statex\hspace*{\algorithmicindent}
					\(k\gets k+1\) ~and go to \cref{step:ADMM:1}
			\ElsIf{~~\(i_k=i_{\rm max}\)~~}
				\Comment{Max \#backtrackings: do nominal ADMM step}
				\Statex\hspace*{\algorithmicindent}
					\(
						(x^{k+1},y^{k+1},z^{k+1})
					{}\in{}
						\ADMM(\bar y^{k+\nicefrac12},z^k)
					\)
				\Statex\hspace*{\algorithmicindent}
					\(k\gets k+1\)~ and go to \cref{step:ADMM:1}%
			\Else
				\Comment{Backtrack and retry}
				\Statex\hspace*{\algorithmicindent}\label{step:ADMM:tau/2}
					\(\tau_k\gets\nicefrac{\tau_k}{2}\),~
					\(i_k\gets i_k+1\),~
					and restart from \cref{step:ADMM:1/2}
			\EndIf
			\end{algorithmic}
				\end{boxedalgorithm}
			\end{figure}
			
			Although a complete and rigorous discussion will be given in the dedicated \cref{sec:Algorithms} after the needed preliminaries have been dealt with, the overall methodology is quite simple and can be informally summarized in few sentences.
			The proposed algorithms leverage on favorable properties of the Douglas-Rachford envelope (DRE) \cite{patrinos2014douglas,themelis2020douglas}, a continuous, real-valued, exact penalty function for problem \eqref{eq:P}.
			Given a stepsize \(\gamma>0\), the DRE associated to problem \eqref{eq:P} at a point \(s\in\R^p\) is given by
			\[
				\DRE(s)
			{}\coloneqq{}
				\sDRS(u)
				{}+{}
				\nsDRS(v)
				{}+{}
				\tfrac1\gamma
				\innprod{s-u}{v-u}
				{}+{}
				\tfrac{1}{2\gamma}\|v-u\|^2,
			\]
			where \((u,v)\in\DRS(s)\) is the result of a \ref{DRS}-update with stepsize \(\gamma\) at \(s\).
			Under the assumptions on the problem dealt in this paper (see \cref{sec:Blocks} for the details), the right-hand side in the above equation is the same for any \((u,v)\in\DRS(s)\), making \(\DRE\) a well-defined function of variable \(s\).
			In fact, the DRE turns out to be continuous and such that
			\begin{equation}\label{eq:Lyapunov}
				\DRE(\bar s^+)
			{}\leq{}
				\DRE(s)-\tfrac{C}{\gamma}\|v-u\|^2
			\end{equation}
			holds for any \(u,v,\bar s^+\) resulting from a \ref{DRS} update at \(s\), provided that the stepsize \(\gamma\) and the relaxation \(\lambda\) are selected as instructed in \cref{alg:DRS}'s initialization.
			By replacing \(C\) with a smaller value \(c<C\), thanks to the aforementioned continuity of \(\DRE\) not only will this inequality hold for \(\bar s^+\), but also for all the points in a neighborhood.
			This enables the flexibility to choose an arbitrary update direction \(d\), stemming for instance from a quasi-Newton method of choice, and backtrack the sought update \(s+d\) towards \(\bar s^+\) until it enters said neighborhood.
			
			Skimming through the steps of \cref{alg:DRS}, the decrease constant \(c<C\) is required at initialization time, and computed based on the problem assumptions.
			The backtracking instead occurs at \cref{step:DRS:1/2}, where by halvening \(\tau_k\) enough times the average between the standard \ref{DRS}-update \(\bar s^{k+1}\) and the sought custom update \(s^k+d^k\) will be close enough to \(\bar s^{k+1}\), causing the inequality in \cref{step:DRS:LS} to be satisfied.
			Surprisingly, the analysis of \cref{alg:DRS} will be enough to cover the \ref{ADMM} counterpart of \cref{alg:ADMM}, for the steps in both methods will be shown to coincide.

		\subsection{Connections with previous work}\label{sec:Previous}

			This work is based on the theoretical analysis of the original \ref{DRS} and \ref{ADMM} algorithms in the nonconvex setting developed in the paper \cite{themelis2020douglas} by some of the authors.
			In fact, part of the content originally appeared in a preprint (non peer-reviewed) version, and was then removed because of page limitation and lack of focus.
			The current version benefits from a thorough polishing and many additions, which include finite-termination and superlinear convergence analysis, strongly convex case through self-duality arguments, simplified assumptions, and numerical simulations.
			An overview of the results in (the latest, published, version of) \cite{themelis2020douglas} needed in this paper is given in \cref{sec:Blocks}, which is dedicated to the properties of \ref{DRS} and \ref{ADMM} oracles and the relations existing among the two.
			
			The proposed linesearch follows the same rationale of the predecessors {\rm PANOC} \cite{stella2017simple} and {\rm NAMA} \cite{stella2019newton} algorithms, which instead of \ref{DRS} and \ref{ADMM} oracles were based on proximal gradient and alternating minimization steps, respectively.
			Their rationale hinges on the same mechanisms explained in \cref{sec:AlgIntro}, namely the identification of a continuous ``envelope'' serving as Lyapunov function for the nominal iterations, \ie, ensuring an inequality in the likes of \eqref{eq:Lyapunov}, and a backtracking to bias the desired update close enough to the nonminal algorithmic step in such a way to enforce a decrease condition on the envelope, as done here in \cref{step:DRS:LS}.
			An initial attempt to encompass these methods under the same lens is given in the doctoral dissertation \cite{themelis2018proximal} with the \emph{Continuous-Lyapunov Descent} framework ({\rm CLyD}), which however is not based on peer-reviewed material and offers limited theoretical results.
			
			Related, but substantially different, is the {\rm minFBE} algorithm of \cite{stella2017forward}, based on forward-backward iterations.
			Under a convexity assumption on the nonsmooth term, the nonsmooth minimization problem can be cast as the minimization of the continuously differentiable ``forward-backward envelope'' function (FBE), where classical smooth minimization techniques (with standard linesearch strategies) can directly be applied.
			Apart from the additional convexity assumption, not required in other methods, this early work suffers from a complication of the algorithmic oracle required for the evaluation of the gradient of the FBE, which entails Hessian evaluations of the smooth function (and, as a byproduct, an additional twice differentiability assumption).

		\subsection{Notation and known facts}\label{sec:Notation}

			We denote as \(\Rinf\coloneqq\R\cup\set{\infty}\) the extended-real line.
			With \(\id\) we indicate the identity function \(x\mapsto x\) defined on a suitable space, and with \(\I\) the identity matrix of suitable size.
			The distance of a point \(x\in\R^n\) to a nonempty set \(E\subseteq\R^n\) is given by \(\dist(x,E)=\inf_{z\in E} \|z-x\|\).
			The relative interior of \(E\), denoted \(\relint E\), is the interior of \(E\) relative to the smallest affine space containing \(E\).
			
			For a sequence \(\seq{x^k}\) we write \(\seq{x^k}\subset E\) to indicate that \(x^k\in E\) for all \(k\in\N\).
			We say that \(\seq{x^k}\subset\R^n\) converges at \(R\)-linear rate (to a point \(x_\star\)) if there exists \(c>0\) and \(\rho\in(0,1)\) such that \(\|x^k-x_\star\|\leq c\rho^k\) holds for every \(k\), and at superlinear rate (to \(x_\star\)) if either \(s^k=s_\star\) for some \(k\) or \(\frac{\|s^{k+1}-s_\star\|}{\|s^k-s_\star\|}\to0\) as \(k\to\infty\).
			
			A function \(\func{h}{\R^n}{\Rinf}\) is \DEF{proper} if \(h\not\equiv\infty\), in which case its \DEF{domain} is defined as the set \(\dom h\coloneqq\set{x\in\R^n}[h(x)<\infty]\), and is \DEF{lower semicontinuous (lsc)} if for any \(\bar x\in\R^n\) it holds that \(h(\bar x)\leq\liminf_{x\to\bar x}h(x)\).
			A point \(x_\star\in\dom h\) is a \DEF{local minimum} for \(h\) if \(h(x)\geq h(x_\star)\) holds for all \(x\) in a neighborhood of \(x_\star\).
			If the inequality can be strengthened to \(h(x)\geq h(x_\star)+\tfrac\mu2\|x-x_\star\|^2\) for some \(\mu>0\), then \(x_\star\) is a \DEF{strong local minimum}.
			We say that \(h\) is \DEF{level bounded} if the sublevel sets
			\(
				\set{x\in\R^n}[
					h(x)\leq\alpha
				]
			\)
			are bounded for all \(\alpha\in\R\), a condition which is equivalent to \(\liminf_{\|x\|\to\infty}f(x)=\infty\).
			
			We denote by \(\hat\partial h\) the \DEF{regular subdifferential} of \(h\), where
			\begin{equation}\label{eq:hatpartial}
				v\in\hat\partial h(\bar x)
			\quad\Leftrightarrow\quad
				\liminf_{\bar x\neq x\to\bar x}{
					\frac{h(x)-h(\bar x)-\innprod{v}{x-\bar x}}{\|x-\bar x\|}
				}
			{}\geq{}
				0.
			\end{equation}
			The (limiting) \DEF{subdifferential} of \(h\) is \(\ffunc{\partial h}{\R^n}{\R^n}\), where
			\(
				v\in\partial h(x)
			\)
			iff \(x\in\dom h\) and there exists a sequence \(\seq{x^k,v^k}\) with \(v^k\in\hat\partial h(x^k)\) such that
			\(
				(x^k,h(x^k),v^k)
			{}\to{}
				(x,h(x),v)
			\)
			as \(k\to\infty\).
			A necessary condition for local minimality of \(x\) for \(h\) is \(0\in\hat\partial h(x)\), see \cite[Thm. 10.1]{rockafellar2011variational}.
			
			The \DEF{proximal mapping} of \(\func{h}{\R^n}{\Rinf}\) with parameter \(\gamma>0\) is the set-valued mapping defined as
			\begin{align}
				\prox_{\gamma h}(x)
			{}\coloneqq{} &
				\argmin_{w\in\R^n}\set{
					h(w) + \tfrac{1}{2\gamma}\|w-x\|^2
				}.
			\shortintertext{%
				The value function of the corresponding minimization problem, namely the \DEF{Moreau envelope} with stepsize \(\gamma\), is denoted as
			}
				h^\gamma(x)
			{}\coloneqq{} &
				\inf_{w\in\R^n}\set{
					h(w) + \tfrac{1}{2\gamma}\|w-x\|^2
				}.
			\end{align}
			The necessary optimality condition for the problem defining \(\prox_{\gamma h}\) together with the calculus rule of \cite[Ex. 8.8]{rockafellar2011variational} implies
			\begin{equation}\label{eq:prox:subdiff}
				\tfrac1\gamma(x-\bar x)
			{}\in{}
				\hat\partial h(\bar x)
			\quad
				\forall\bar x\in\prox_{\gamma h}(x).
			\end{equation}
			
			With \(\conj h\) we indicate the convex conjugate of function \(h\), pointwise defined as \(\conj h(y)=\sup\set{\innprod{y}{{}\cdot{}}-h}\).
			We say that \(h\) is \DEF{strongly convex} if there exists a constant \(\mu_h>0\) such that \(h-\frac{\mu_h}{2}\|{}\cdot{}\|^2\) is convex, in which case we may say that \(h\) is \DEF{\(\mu_h\)-strongly convex} to make the constant \(\mu_h\) explicit.

	\section{The building blocks: DRS and ADMM}\label{sec:Blocks}

		\ref{DRS} and \ref{ADMM} are nowadays considered textbook algorithms of the realm of convex optimization, and their properties are well documented in the literature.
		For instance, it is common knowledge that both algorithms are mutually equivalent when applied to the respective dual formulations, and that convergence is guaranteed for arbitrary stepsize and penalty parameters under minimal assumptions.
		These algorithms are in fact well understood through an elegant and powerful link with monotone operator theory, a connection in which convexity plays an indispensable role and which can explain convergence through a Fejér-type monotonicity of the generated sequences.
		This property entails the existence of a constant \(c>0\) such that
		\begin{equation}\label{eq:Fejer}
			\|s^+-s_\star\|^2
		{}\leq{}
			\|s-s_\star\|^2
			{}-{}
			c\|s-s^+\|^2
		\end{equation}
		holds for every solution \(s_\star\); by telescoping the inequality and with no information about the whereabouts of any solution \(s_\star\) required, from the  lower boundedness of the (squared) norm it is immediate to deduce that the \emph{residual} \(\|s^+-s\|\) vanishes.
		
		Under some smoothness assumption, a new descent condition in the form of \eqref{eq:Fejer} was shown to hold even for nonconvex problems, with the squared distance \(\|{}\cdot{}-s_\star\|^2\) being replaced by another lower bounded function; upon adopting the same telescoping arguments, this led to new convergence results in the absence of convexity.
		In this paper we show that the new descent condition also leads to linesearch extensions of \ref{DRS} and \ref{ADMM} that preserve the same convergence properties and oracle complexity.
		In this section we present all the preliminary material that is needed for their development.
		We begin with \ref{DRS}, first by offering a brief recap of the key inequalities of the nonconvex analysis developed in \cite{themelis2020douglas} and then by showing, through duality arguments, that the needed smoothness requirement can be replaced by strong convexity.
		Although convergence of \ref{DRS} is well known in the latter case, it allows us to generalize the standing assumptions of the proposed linesearch algorithms.
		Finally, by means of a primal equivalence first noticed in \cite{yan2016self} that identifies \ref{DRS} and \ref{ADMM}, we will obtain a similar analysis for the latter algorithm.

		\subsection{Douglas-Rachford splitting}
			
			\subsubsection{The nonconvex case}

				For convex problems, both \ref{DRS} and \ref{ADMM} are well known to converge for arbitrary stepsize and penalty parameters under minimal assumptions.
				In the nonconvex setting, the works \cite{li2016douglas,li2017peaceman} and later \cite{themelis2020douglas} extended the analysis to the nonconvex case when the functions satisfy the following requirements.
				
				\begin{ass}[Requirements for \ref{DRS}: the smooth \(\sDRS\) case]\label{ass:DRS}%
					In problem \eqref{eq:P}:
					\begin{enumeratass}
					\item\label{ass:sDRS}
						\(\func{\sDRS}{\R^p}{\R}\) has \(L_{\sDRS}\)-Lipschitz continuous gradient.
					\item\label{ass:nsDRS}
						\(\func{\nsDRS}{\R^p}{\Rinf}\) is lsc.
					\item\label{ass:phiDRS}
						Problem \eqref{eq:P} admits a solution: \(\argmin\varphi\neq\emptyset\).
					\end{enumeratass}
				\end{ass}
				
				In the setting of \cref{ass:DRS}, \cite{li2016douglas,li2017peaceman} pioneered the idea of employing an augmented Lagrangian function as Lyapunov potential for \ref{DRS} iterations, namely
				\[
					\L_c(u,v,y)
				{}\coloneqq{}
					\sDRS(u)
					{}+{}
					\nsDRS(v)
					{}+{}
					\innprod{y}{v-u}
					{}+{}
					\tfrac c2\|v-u\|^2
				\]
				for some \(c\in\R\) (not necessarily positive).
				It was shown that \(\L_c\bigl(u,v,\gamma^{-1}(s-u)\bigr)\) decreases along the iterates generated by \ref{DRS} for sufficiently small \(\gamma\), and subsequential convergence of the algorithm to stationary points was thus inferred.
				The results have been improved in \cite{themelis2020douglas}, where considering \(c=\nicefrac1\gamma\) is shown to lead to a tight convergence analysis that cannot be further improved without restricting the problem assumptions.
				Moreover, the augmented Lagrangian is regarded as a function of the sole variable \(s\), and thus coincides with the \emph{Douglas-Rachford envelope} (DRE) of \cite{patrinos2014douglas}, namely
				\begin{equation}\label{eq:DRE}
					\DRE(s)
				{}\coloneqq{}
					\sDRS(u)
					{}+{}
					\nsDRS(v)
					{}+{}
					\tfrac1\gamma
					\innprod{s-u}{v-u}
					{}+{}
					\tfrac{1}{2\gamma}\|v-u\|^2,
				\end{equation}
				where \(u\) and \(v\) are the result of a \ref{DRS}-update with stepsize \(\gamma\) starting at \(s\) (\(u\) and \(v\) are independent of the relaxation \(\lambda\)).
				As detailed in \cite[Prop. 2.3 and Rem. 3.1]{themelis2020douglas}, when \cref{ass:DRS} is satisfied and \(\gamma<\nicefrac{1}{L_{\sDRS}}\) one has that \(\prox_{\gamma\sDRS}\) is Lipschitz continuous and \(\prox_{\gamma\nsDRS}\) is a well-defined set-valued mapping, in the sense that \(\prox_{\gamma\nsDRS}(x)\) is nonempty for every \(x\in\R^p\).
				In fact, in order for the DRE to be well defined it suffices that \(\prox_{\gamma\sDRS}\) is single valued, as the expression \eqref{eq:DRE} can easily be seen to equal
				\begin{equation}\label{eq:DREMoreau}
					\DRE(s)
				{}={}
					\sDRS^\gamma(s)
					{}-{}
					\tfrac1\gamma\|s-u\|^2
					{}+{}
					\nsDRS^\gamma(2u-s),
				\quad
					\text{with }u=\prox_{\gamma\sDRS}(s).
				\end{equation}
				Nevertheless, under \cref{ass:DRS} the DRE enjoys a close kinship with the cost function \(\varphi\), as summarized next.
				
				\begin{fact}[{\cite[Prop. 3.2 and Thm. 3.4]{themelis2020douglas}}]\label{thm:DREequiv}%
					Suppose that \cref{ass:DRS} holds.
					Then, for all \(\gamma<\nicefrac{1}{L_{\sDRS}}\) the DRE \(\DRE\) is real valued, locally Lipschitz, and satisfies the following:%
					\begin{enumerate}
					\item\label{thm:DREequiv:inf}%
						\(
							\inf\varphi=\inf\DRE
						\)
						~and~
						\(
							\argmin\varphi=\prox_{\gamma\sDRS}(\argmin\DRE)
						\).
					\item\label{thm:DREequiv:LB}%
						\(\DRE\) is level bounded iff \(\varphi\) is level bounded.
					\end{enumerate}
				\end{fact}
				
				In the same spirit as the preceding works \cite{li2016douglas,li2017peaceman}, the convergence analysis of nonconvex \ref{DRS} in \cite{themelis2020douglas} revolves around the following result, which assesses that that the DRE decreases along the iterations by a quantity which is proportional to the fixed-point residual \(\|s^k-s^{k+1}\|^2\).
				For simplicity of exposition, we use the simplified bounds on the stepsize \(\gamma\) as in \cite[Rem. 4.2]{themelis2020douglas}, which only discern whether \(\sDRS\) is convex or not.
				Tight ranges are given in \cite[Thm. 4.1]{themelis2020douglas}, and require the knowledge of the \emph{hypoconvexity} modulus of \(\sDRS\).
				
				\begin{fact}[Sufficient decrease on the DRE {\cite[Thm. 4.1]{themelis2020douglas}}]\label{thm:DRE:SD}%
					Suppose that \cref{ass:DRS} is satisfied, and consider one \ref{DRS} update
					\(
						s
					{}\mapsto{}
						(u,v,s^+)
					\)
					for some stepsize
					\(
						\gamma
					{}<{}
						\nicefrac{1}{L_{\sDRS}}
					\)
					and relaxation \(\lambda\in(0,2)\).
					Then,
					\begin{equation}\label{eq:DRE:SD}
						\DRE(s^+)
					{}\leq{}
						\DRE(s)
						{}-{}
						\tfrac1\gamma
						\C\|r\|^2,
					\end{equation}
					where \(r\coloneqq u-v\) and
					\begin{equation}\label{eq:DRE:C}
						\C[\alpha][\lambda]
					{}\coloneqq{}
						\tfrac{\lambda}{(1+\alpha)^2}
						\left(
							\tfrac{2-\lambda}{2}
							{}-{}
							\alpha
							\cdot
							\begin{ifcases}
								\max\set{\alpha-\nicefrac\lambda2,0}
							&
								\sDRS\text{ is convex}
							\\
								1
							\otherwise
							\end{ifcases}
						\right).
					\end{equation}
					In particular, the constant \(\C[]\) is strictly positive provided that
					\begin{equation}\label{eq:DRE:gamma}
						\gamma
					{}<{}
						\begin{ifcases}
							\frac{1}{L_{\sDRS}}
						&
							\sDRS\text{ is convex}
						\\[5pt]
							\tfrac{2-\lambda}{2L_{\sDRS}}
						\otherwise.
						\end{ifcases}
					\end{equation}
				\end{fact}
				
				Combined with the lower boundedness of the DRE (\cref{thm:DREequiv}), the vanishing of \(s^k-s^{k+1}=\lambda(u^k-v^k)\) for the nominal \ref{DRS} algorithm was readily deduced.
				This is enough to guarantee \emph{subsequential convergence} of \ref{DRS} to stationary points, in the sense that whenever a sequence \(\seq{s^k}\) satisfies \(u^k-v^k\to0\) with \((u^k,v^k)\in\DRS(s^k)\), any accumulation point \(u^\star\) of \(\seq{v^k}\) satisfies the stationarity condition \(0\in\partial\varphi(u^\star)\).
				As detailed in the dedicated \cref{sec:Convergence}, the analysis of the here proposed \ref{DRS}-based \cref{alg:DRS} follows the same line of proof, revolving around the \emph{decrease condition} in \cref{thm:DRE:SD} and the \emph{continuity} of the DRE in \cref{thm:DREequiv}.
				Before detailing the arguments, in the remaining subsections we set the ground for extending the theory beyond \cref{ass:DRS}.
				First, with duality arguments we replace the smoothness condition with a strong convexity requirement; then, by means of the same change of variable adopted in \cite{themelis2020douglas} we express \ref{ADMM} operations \eqref{eq:ADMM} in terms of the \ref{DRS} oracle \eqref{eq:DRS}, thus obtaining the \ref{ADMM}-based \cref{alg:ADMM} as a simple byproduct.

			\subsubsection{The strongly convex case}\label{sec:SelfEquivalence}

				As observed in \cite{yan2016self}, for convex problems \ref{DRS} (or, equivalently, its sibling \ref{ADMM}) is equivalent to itself applied to the dual formulation.
				To see this, let \(\sDRS\) and \(\nsDRS\) be (proper, lsc) convex functions, and consider the following dual formulation of \eqref{eq:P}:\footnote{%
					Specifically, \eqref{eq:D} is the dual of
					\(
						\minimize_{x,z\in\R^n}\sDRS(x)+\nsDRS(z)
					\)
					\(
						\stt x-z=0
					\).
				}
				\begin{equation}\label{eq:D}
					\minimize_{y\in\R^n}\psi(y)\coloneqq\sDRS*(y)+\nsDRS*(y),
				\end{equation}
				with \(\nsDRS*=\conj\nsDRS\) and \(\sDRS*\coloneqq\conj{(\refl\sDRS)}=\refl\conj\sDRS\) (see \cite[Prop. 13.23(iv)]{bauschke2017convex} for the last equality), where \(\refl\) denotes the ``mirroring'' \(\refl\sDRS\coloneqq\sDRS(-{}\cdot{})\).
				The Moreau identity \cite[Thm. 14.3(ii)]{bauschke2017convex} yields
				\begin{equation}\label{eq:prox*}
					\prox_{\delta\sDRS*}(y)
				{}={}
					y+\delta\prox_{\nicefrac\sDRS\delta}(-\nicefrac y\delta)
				\quad\text{and}\quad
					\prox_{\delta\nsDRS*}(y)
				{}={}
					y-\delta\prox_{\nicefrac\nsDRS\delta}(\nicefrac y\delta).
				\end{equation}
				Therefore, \ref{DRS} applied to \eqref{eq:D} produces the following triplet:
				{%
					\[\tag{\text{\rm DRS\textsuperscript*}}\label{eq:DRS*}
						s_*
					{}\mapsto{}
						\begin{cases}[l >{{}}c<{{}} l]
							u_* &=& \prox_{\gamma_*\sDRS*}(s_*)\\
							v_* &=& \prox_{\gamma_*\nsDRS*}(2u_*-s_*)\\
							s_*^+ &=& s_*+\lambda_*(v_*-u_*),
						\end{cases}
					\quad\Leftrightarrow\quad
						\begin{cases}[l >{{}}c<{{}} l]
							u_* &=& s_*+\gamma_*\prox_{\nicefrac{\sDRS}{\gamma_*}}(-\nicefrac{s_*}{\gamma_*})\\
							v_* &=& 2u_*-s_*-\gamma_*\prox_{\nicefrac{\nsDRS}{\gamma_*}}(\nicefrac{(2u_*-s_*)}{\gamma_*})\\
							s_*^+ &=& s_*+\lambda_*(v_*-u_*),
						\end{cases}
					\]
				}%
				leading to the following result.
				
				\begin{thm}[Self-duality of DRS]\label{thm:selfdual}%
					Suppose that \(\sDRS\) and \(\nsDRS\) are (proper, lsc) convex functions, and consider a \ref{eq:DRS*}-update \(s_*\mapsto(u_*,v_*,s_*^+)\).
					Let \(s\coloneqq-\frac{s_*}{\gamma_*}\) and consider a \ref{DRS}-update \(s\mapsto(u,v,s^+)\) with stepsize \(\gamma=\nicefrac{1}{\gamma_*}\) and relaxation \(\lambda=\lambda_*\).
					Then, the variables are related as follows:
					\begin{equation}\label{eq:DRS-DRS*}
						\begin{cases}[r @{{}={}} l]
							\lambda & \lambda_*\\
							\gamma & \frac{1}{\gamma_*}\\
							s & -\frac{s_*}{\gamma_*}\\
							u & \frac{u_*-s_*}{\gamma_*}\\
							v & \frac{2u_*-s_*-v_*}{\gamma_*},
						\end{cases}
					\qquad\text{or, equivalently,}\qquad
						\begin{cases}[r @{{}={}} l]
							\lambda_* & \lambda\\
							\gamma_* & \frac1\gamma\\
							s_* & -\frac s\gamma\\
							u_* & \frac{u-s}{\gamma}\\
							v_* & \frac{2u-s-v}{\gamma}.
						\end{cases}
					\end{equation}
					Moreover,
					\(
						\grayout{
							\begin{cases}[r @{{}\in{}} l]
								-u_* & \partial\sDRS(u) \\
								v_*  & \partial\nsDRS(v)
							\end{cases}
						\quad\text{and}\quad
						}
						\DRE*(s_*) = -\DRE(s),
					\)
					where \(\DRE*\) is the DRE with stepsize \(\gamma_*\) associated to the dual problem \eqref{eq:D}.
					\begin{proof}
						The identities in \eqref{eq:DRS-DRS*} follow by a direct application of those in \eqref{eq:DRS*}.
						Next,
						\[
							\sDRS*^{\gamma_*}(s_*)
						{}={}
							(\conj{(\refl\sDRS)})^{\nicefrac1\gamma}(-\nicefrac s\gamma)
						{}={}
							(\conj{\sDRS})^{\nicefrac1\gamma}(\nicefrac s\gamma)
						{}={}
							\tfrac{1}{2\gamma}\|s\|^2
							{}-{}
							\sDRS^\gamma(s),
						\]
						and similarly
						\(
							\nsDRS*^{\gamma_*}(2u_*-s_*)
						{}={}
							\tfrac{1}{2\gamma}\|2u-s\|^2
							{}-{}
							\nsDRS^\gamma(2u-s)
						\).
						Therefore, from \eqref{eq:DREMoreau} we have
						\begin{align*}
							\DRE*(s_*)
						{}={} &
							\fillwidthof[c]{
								\tfrac{1}{2\gamma}\|s\|^2
								{}-{}
								\sDRS^\gamma(s)
							}{
								\sDRS*^{\gamma_*}(s_*)
							}
							{}-{}
							\tfrac{1}{\gamma_*}\|u_*-s_*\|^2
							\fillwidthof[c]{
								{}+{}
								\tfrac{1}{2\gamma}\|2u-s\|^2
								{}-{}
								\nsDRS^\gamma(2u-s)
							}{
								{}+{}
								\nsDRS*^{\gamma_*}(2u_*-s_*)
							}
						\\
						{}={} &
							\overbrace{
								\tfrac{1}{2\gamma}\|s\|^2
								{}-{}
								\sDRS^\gamma(s)
							}
							\overbrace{
								\fillwidthof[c]{
									{}-{}
									\tfrac{1}{\gamma_*}\|u_*-s_*\|^2
								}{
									{}-{}
									\tfrac{1}{\gamma}\|u\|^2
								}
							}
							\overbrace{
								{}+{}
								\tfrac{1}{2\gamma}\|2u-s\|^2
								{}-{}
								\nsDRS^\gamma(2u-s)
							}
						\\
						{}={} &
							{}-{}
							\sDRS^\gamma(s)
							{}+{}
							\tfrac1\gamma\|u-s\|^2
							{}-{}
							\nsDRS^\gamma(2u-s)
						\\
						{}\overrel*{\eqref{eq:DREMoreau}}{} &
							{}-{}
							\DRE(s)
						\end{align*}
						as claimed.
					\end{proof}
				\end{thm}
				
				The conjugate \(\conj h\) of a proper convex lsc function \(h\) is also proper convex and lsc; when \(h\) is \(\mu_h\)-strongly convex, then \(\conj h\) has also \(\nicefrac{1}{\mu_h}\)-Lipschitz-continuous gradient.
				This means that whenever \(\nsDRS\) is proper lsc and convex and \(\sDRS\) is additionally strongly convex, the theory presented in the previous subsection applies to the dual formulation \eqref{eq:D}, provided dual solutions exist.
				As shown in \cite[Cor. 27.6.(i)a]{bauschke2017convex}, this is guaranteed whenever a constraint qualification on the domains of the two functions holds,\footnote{%
					Owing to strong convexity, the additional requirement of nonemptyness of \(\argmin\sDRS+\nsDRS\) in the cited reference would be trivially satisfied.
				}
				which leads to the following convex, but fully nonsmooth, dual version of \cref{ass:DRS}.
				
				\bgroup
					\let\theoldass\theass
					\let\theoldHass\theHass
					\addtocounter{ass}{-1}%
					\renewcommand{\theass}{\theoldass*}%
					\renewcommand{\theHass}{\theoldHass*}%
					\begin{ass}[Requirements for \ref{DRS}: the strongly convex case]\label{ass:DRS*}%
						In problem \eqref{eq:P} the following hold:
						\begin{enumeratass}
						\item\label{ass:sDRS*}%
							\(\func{\sDRS}{\R^p}{\Rinf}\) is proper, lsc, and \(\mu_{\sDRS}\)-strongly convex.
						\item\label{ass:nsDRS*}%
							\(\func{\nsDRS}{\R^p}{\Rinf}\) is proper, lsc, and convex.
						\item\label{ass:DRS*:relint}%
							\(0\in\relint(\dom\sDRS-\dom\nsDRS)\).
						\end{enumeratass}
					\end{ass}
				\egroup
				
				When the functions are convex and without necessarily either one being strongly so, plain \ref{DRS} iterations are known to converge for any stepsize \(\gamma\) and relaxation \(\lambda\in(0,2)\) provided that \(\partial\sDRS+\partial\nsDRS\) attains a zero \cite[Cor. 28.3]{bauschke2017convex}.
				\Cref{ass:DRS*:relint} is a mild and easily verifiable domain qualification ensuring this latter property, and is therefore quite standard.
				What is instead crucial in our setting is \emph{strong} convexity, as it guarantees well definedness of \cref{alg:DRS} even in the absence of smoothness.
				The duality pairing strong convexity and smoothness was also exploited in \cite{patrinos2014douglas}, where the DRE was first introduced and \ref{DRS} was shown to be a scaled gradient descent algorithm on it.
				This observation led to the development of a Nesterov-type acceleration whenever \(\nsDRS\) is convex and \(\sDRS\) is strongly convex and quadratic.
				The analysis was then extended to \ref{ADMM} in the follow-up work \cite{pejcic2016accelerated} by means of duality arguments.
				
				A consequence of \cref{thm:selfdual} is the following dual version of \cref{thm:DRE:SD}.
				
				\begin{cor}\label{thm:DRE:SD*}%
					Suppose that \cref{ass:DRS*} holds and consider one \eqref{DRS} update \(s\mapsto(u,v,s^+)\) applied to the primal formulation \eqref{eq:P} with stepsize \(\gamma>\nicefrac{1}{\mu_{\sDRS}}\) and relaxation \(\lambda\in(0,2)\).
					Then,
					\[
						\DRE(s^+)
					{}\geq{}
						\DRE(s)
						{}+{}
						\tfrac1\gamma\C*\|r\|^2
					\]
					where \(r\coloneqq u-v\) and \(\C*[]\) is a strictly positive constant defined as in \eqref{eq:DRE:C}.
				\end{cor}
				
				\grayout{%
				% {\color{blue}%
					\begin{align*}
						\tfrac{1}{2\gamma}\|x_\star-v\|^2
						{}+{}
						\tfrac{\gamma\mu_{\sDRS}-1}{2\gamma}\|x_\star-u\|^2
					{}\leq{} &
						\inf\varphi
						{}-{}
						\DRE(s)
					{}={}
						\DRE*(s_*)
						{}-{}
						\inf\psi
					\\
					{}\leq{} &
						\tfrac{1+\gamma_*L_{\sDRS*}}{2\gamma_*}\|u_*-u_*^\star\|^2
					\end{align*}
					
					\begin{align*}
						\DRE*(s_*)
					{}={} &
						\sDRS*(u_*)
						{}+{}
						\nsDRS*(v_*)
						{}+{}
						\innprod{\nabla\sDRS*(u_*)}{v_*-u_*}
						{}+{}
						\tfrac{1}{2\gamma_*}\|v_*-u_*\|^2
					\\
					{}\leq{} &
						\sDRS*(v_*)
						{}+{}
						\nsDRS*(u_*^\star)
						{}+{}
						\innprod{2u_*-s_*-v_*}{v_*-u_*^\star}
						{}+{}
						\tfrac{1-\gamma_*L_{\sDRS*}}{2\gamma_*}\|v_*-u_*\|^2
					\\
					{}\leq{} &
						\psi_\star
						{}+{}
						\tfrac1\gamma_*\innprod{2u_*-s_*-v_*}{v_*-u_*^\star}
						{}+{}
						\innprod{\nabla\sDRS*(v_*)}{v_*-u_*^\star}
						{}+{}
						\tfrac{1-\gamma_*L_{\sDRS*}}{2\gamma_*}\|v_*-u_*\|^2
					\\
					{}={} &
						\psi_\star
						{}+{}
						\tfrac1\gamma_*\innprod{u_*-v_*}{v_*-u_*^\star}
						{}+{}
						\innprod{\nabla\sDRS*(v_*)-\nabla\sDRS*(u_*)}{v_*-u_*^\star}
						{}+{}
						\tfrac{1-\gamma_*L_{\sDRS*}}{2\gamma_*}\|v_*-u_*\|^2
					\\
					{}\leq{} &
						\psi_\star
						{}+{}
						\tfrac1\gamma_*\innprod{u_*-v_*}{v_*-u_*^\star}
						{}+{}
						\tfrac{\varepsilon}{2\gamma_*}\|v_*-u_*^\star\|^2
						{}+{}
						\tfrac{\gamma_*^2L_{\sDRS*}^2+\varepsilon(1-\gamma_*L_{\sDRS*})}{2\gamma_*\varepsilon}\|v_*-u_*\|^2
					\\
					\numberthis\label{eq:ineq}
					{}={} &
						\psi_\star
						{}+{}
						\tfrac1\gamma_*\innprod{u_*-v_*}{u_*-u_*^\star}
						{}+{}
						\tfrac{\varepsilon}{2\gamma_*}\|v_*-u_*^\star\|^2
						{}+{}
						\tfrac{\gamma_*^2L_{\sDRS*}^2-\varepsilon(1+\gamma_*L_{\sDRS*})}{2\gamma_*\varepsilon}\|v_*-u_*\|^2
					\\
					{}\leq{} &
						\psi_\star
						{}+{}
						\tfrac{1}{2\delta\gamma_*}\|u_*-u_*^\star\|^2
						{}+{}
						\tfrac{\varepsilon}{2\gamma_*}\|v_*-u_*^\star\|^2
						{}+{}
						\tfrac{\gamma_*^2L_{\sDRS*}^2-\varepsilon(1+\gamma_*L_{\sDRS*}-\delta)}{2\gamma_*\varepsilon}\|v_*-u_*\|^2.
					\end{align*}
				
					\begin{align*}
						\DRE*(s_*)
						{}-{}
						\psi_\star
					{}\leq{} &
						\tfrac1\gamma_*\innprod{u_*-v_*}{v_*-u_*^\star}
						{}+{}
						\innprod{\nabla\sDRS*(v_*)-\nabla\sDRS*(u_*)}{v_*-u_*^\star}
						{}+{}
						\tfrac{1-\gamma_*L_{\sDRS*}}{2\gamma_*}\|v_*-u_*\|^2
					\\
					{}={} &
						\innprod{v+\nabla\sDRS*(v_*)}{v_*-u_*^\star}
						{}+{}
						\tfrac{1-\gamma_*L_{\sDRS*}}{2\gamma_*}\|v_*-u_*\|^2
					\end{align*}

					\begin{rem}[Existence of dual solutions]%
						Suppose that \cref{ass:DRS*} holds, and let \(x^\star\) be the unique (primal) solution of \eqref{eq:P}.
						Then, \(s^\star\) is a fixed point of \(\DRS\) iff \(\prox_{\gamma\sDRS}(s^\star)=x^\star\), which is equivalent to having \(s^\star=x^\star+\gamma\xi^\star\) with \(\xi^\star\in\partial\sDRS(x^\star)\).
						The \(\DRS\) optimal triplet is thus \((x^\star+\gamma\xi^\star,x^\star,x^\star)\), corresponding to the \(\DRS*[\nicefrac1\gamma]\) triplet \((-\xi^\star-\nicefrac{x^\star}{\gamma},-\xi^\star,-\xi^\star)\) as it follows from \cref{thm:selfdual}.
						In light of \cite[Thm. 5.1.3 and Ex. 5.3.1]{auslender2002asymptotic}, dual solutions exist iff \(\partial q(0)\neq\emptyset\) for \(q(y)=\min_x\Phi(x,y)\coloneqq f(x)+g(x+y)\).
						Since \(x^\star\) is the unique minimizer of \(\Phi({}\cdot{},0)\), \cite[Thm. 10.13]{rockafellar2011variational} ensures that whenever \((0,\xi^\star)\in\partial\Phi(x^\star,0)\) one has \(\xi^\star\in\partial q(0)\).
						Notice that \((0,\xi^\star)\in\partial\Phi(x^\star,0)\) whenever \(\xi^\star\in\partial\nsDRS(x^\star)\), hence we conclude that dual solutions exist whenever either \(\partial\sDRS(x^\star)\neq\emptyset\) or, equivalently, \(\partial\nsDRS(x^\star)\neq\emptyset\).
					\end{rem}
				
					Suppose that \(\partial\sDRS(x^\star)\neq\emptyset\) (or, equivalently, \(\partial\nsDRS(x^\star)\neq\emptyset\)), and let \(\xi^\star\in\partial\sDRS(x^\star)\) be fixed.
					For \(\gamma_*=\nicefrac1\gamma\), consider the \(\DRS*[\gamma_*]\) optimal triplet \((-\xi^\star-\nicefrac{x^\star}{\gamma},-\xi^\star,-\xi^\star)\).
					Then, \eqref{eq:ineq} in primal terms reads
					\[
						\inf\varphi
						{}-{}
						\DRE(s)
					{}\leq{}
						\innprod{v-u}{v_*-u_*}
						{}+{}
						\innprod{v-u}{u_*-u_*^\star}
						{}+{}
						\tfrac{\varepsilon}{2\gamma_*}\|v_*-u_*^\star\|^2
						{}+{}
						\tfrac{\gamma_*^2L_{\sDRS*}^2+\varepsilon(1-\gamma_*L_{\sDRS*})}{2\gamma\varepsilon}\|v-u\|^2.
					\]
				
					There are two options here for the term \(\innprod{u_*-v_*}{v_*-u_*^\star}\): either (1) expanding it into three squares, or (2) using Young inequality.
					\begin{enumerate}[label=(\arabic*)]
					\item
						\begin{align*}
							\DRE*(s_*)
						{}\leq{} &
							\psi_\star
							{}+{}
							\tfrac{1-\gamma_*L_{\sDRS*}}{2\gamma_*}\|v_*-u_*\|^2
					%
					% 		{}+{}
					% 		\tfrac{1}{\gamma_*}\innprod{u_*-v_*}{v_*-u_*^\star}
							{}+{}
							\tfrac{1}{2\gamma_*}\|u_*-u_*^\star\|^2
							{}-{}
							\tfrac{1}{2\gamma_*}\|u_*-v_*\|^2
							{}-{}
							\tfrac{1}{2\gamma_*}\|v_*-u_*^\star\|^2
						\\
						&
							{}+{}
							\tfrac{\gamma_*L_{\sDRS*}^2}{2\varepsilon}\|v_*-u_*\|^2
							{}+{}
							\tfrac{\varepsilon}{2\gamma_*}\|v_*-u_*^\star\|^2
						\\
						{}={} &
							\psi_\star
							{}+{}
							\tfrac{(\gamma_*L_{\sDRS*}-\varepsilon)L_{\sDRS*}}{2\varepsilon}\|v_*-u_*\|^2
							{}+{}
							\tfrac{1}{2\gamma_*}\|u_*-u_*^\star\|^2
							{}+{}
							\tfrac{\varepsilon-1}{2\gamma_*}\|v_*-u_*^\star\|^2,
						\end{align*}
						which combined with \eqref{eq:Qlb} results in
						\[
							\tfrac{1}{2\gamma}\|x_\star-v\|^2
							{}+{}
							\tfrac{\gamma\mu_{\sDRS}-1}{2\gamma}\|x_\star-u\|^2
						{}\leq{}
							\tfrac{(\gamma_*L_{\sDRS*}-\varepsilon)L_{\sDRS*}}{2\varepsilon}\|v_*-u_*\|^2
							{}+{}
							\tfrac{1}{2\gamma_*}\|u_*-u_*^\star\|^2
							{}+{}
							\tfrac{\varepsilon-1}{2\gamma_*}\|v_*-u_*^\star\|^2,
						\]
					\end{enumerate}
				}%

		\subsection{Alternating direction method of multipliers}

			Although apparently more general, it is well known that for convex problems \ref{ADMM} coincides with \ref{DRS} applied to the dual formulation, and vice versa.
			More generally, problem \eqref{eq:CP} can be reduced to a \ref{DRS}-compatible form as
			\begin{equation}\label{eq:CP2P}
				\minimize_{s\in\R^p}{
					\epicomp Af(s)
					{}+{}
					\epicomp Bg(b-s)
				},
			\end{equation}
			where for \(\func h{\R^n}{\Rinf}\) and \(C\in\R^{p\times n}\) we indicate with \(\epicomp Ch\) the \emph{epicomposition}
			\[
				\func{\epicomp Ch}{\R^p}{\Rinf},
			\quad\text{defined as}\quad
				\epicomp Ch(s)
			{}\coloneqq{}
				\inf\set{h(x)}[Cx=s].
			\]
			It was shown in \cite[Thm. 1]{yan2016self} and later generalized in \cite[Thm. 5.5]{themelis2020douglas} that one iteration of \ref{ADMM} applied to \eqref{eq:CP} is equivalent to one step of \ref{DRS} applied to this new reformulation with stepsize \(\gamma=\nicefrac1\beta\), as stated next.
			
			\begin{fact}[Primal equivalence of \ref{DRS} and \ref{ADMM} {\cite[Thm. 5.5]{themelis2020douglas}}]\label{thm:ADMM_DRS}%
				Starting from a triplet \((x,y,z)\in\R^m\times\R^p\times\R^n\), consider an \ref{ADMM}-update applied to problem \eqref{eq:CP} with relaxation \(\lambda\) and penalty \(\beta>0\).
				Let
				\[
					\begin{cases}[l >{{}}c<{{}} l]
						s
					&\coloneqq&
						Ax-\nicefrac y\beta
					\\
						u
					&\coloneqq&
						Ax
					\\
						v
					&\coloneqq&
						b-Bz
					\end{cases}
				\quad\text{and, similarly,}\quad
					\begin{cases}[l >{{}}c<{{}} l]
						s^+
					&\coloneqq&
						Ax^+-\nicefrac{y^+}\beta
					\\
						u^+
					&\coloneqq&
						Ax^+
					\\
						v^+
					&\coloneqq&
						b-Bz^+.
					\end{cases}
				\]
				Then, the variables are related as follows:
				\[
					\begin{cases}[l >{{}}c<{{}} l]
						s^+
					&=&
						s+\lambda(v-u)
					\\
						u^+
					&\in&
						\prox_{\gamma\sDRS}(s^+)
					\\
						v^+
					&\in&
						\prox_{\gamma\nsDRS}(2u^+-s^+),
					\end{cases}
				\quad\text{where}\quad
					\begin{cases}[l >{{}}c<{{}} l]
						\sDRS
					&\coloneqq&
						\epicomp Af
					\\
						\nsDRS
					&\coloneqq&
						\epicomp Bg(b-{}\cdot{})
					\\
					\gamma
					&\coloneqq&
						\nicefrac1\beta.
					\end{cases}
				\]
				Moreover,
				\begin{enumerate}
				\item\label{thm:ADMM:-Ay}%
					\(
						-\trans Ay^+\in\hat\partial f(x^+)
					\), and
				\item\label{thm:ADMM:-By}%
					\(
						\dist(-\trans By^+,\hat\partial g(z^+))
					{}\leq{}
						\beta\|B\|
						\|Ax^++Bz^+-b\|
					\).
				\end{enumerate}
			\end{fact}

			This enabled the possibility to infer the convergence of nonconvex \ref{ADMM} from the simpler analysis of that of \ref{DRS}, when the reformulation \eqref{eq:CP2P} complies with the needed \ref{DRS} requirements.
			
			\begin{ass}[Requirements for \ref{ADMM}: the smooth \(\epicomp Af\) case]\label{ass:ADMM}%
				In problem \eqref{eq:CP}, \(A\in\R^{p\times m}\), \(B\in\R^{p\times n}\), \(b\in\R^p\), \(\func{f}{\R^m}{\R}\) and \(\func{g}{\R^n}{\Rinf}\) are such that:
				\begin{enumeratass}
				\item\label{ass:f}%
					\(A\) is surjective (full row rank) and \(\func{\epicomp Af}{\R^p}{\R}\) has \(L_{\epicomp Af}\)-Lipschitz gradient.
				\item\label{ass:g}%
					\(\func{\epicomp Bg}{\R^p}{\Rinf}\) is lsc.
				\item\label{ass:Phi}%
					a solution exists: \(\argmin\Phi\neq\emptyset\), where \(\Phi(x,z)=f(x)+g(z)+\indicator_{Ax+Bz=b}(x,z)\).
				\end{enumeratass}
			\end{ass}
			
			Similarly, in parallel to \cref{ass:DRS*} we may also consider the strongly convex case for \ref{ADMM}.
			As shown in \cite[Prop. 5.4]{themelis2020douglas}, for any matrix \(A\) the function \(\sDRS\coloneqq\epicomp Af\) is strongly convex whenever so is \(f\), in which case the strong convexity modulus is \(\nicefrac{\mu_f}{\|A\|^2}\).\footnote{%
				In the limiting case \(A=0\), one has that \(\epicomp Af=f(0)+\indicator_{\set0}\) is lsc and \(\infty\)-strongly convex for any \(f\), and properness amounts to the condition \(0\in\dom f\).
			}
			Moreover, with \(\nsDRS=\epicomp Bg(b-{}\cdot{})\) one has that \(\dom\sDRS-\dom\nsDRS=A\dom\sADMM+B\dom\nsADMM-b\), owing to the fact that \(\dom\sDRS=A\dom f\) and \(\dom\nsDRS=b-B\dom g\), see \cite[Prop. 12.36(i)]{bauschke2017convex}.
			Altogether, the following \ref{ADMM} counterpart of \cref{ass:DRS*} is obtained.
			
			\bgroup
				\let\theoldass\theass
				\let\theoldHass\theHass
				\addtocounter{ass}{-1}%
				\renewcommand{\theass}{\theoldass*}%
				\renewcommand{\theHass}{\theoldHass*}%
				\phantomsection
				\begin{ass}[Requirements for \ref{ADMM}: the strongly convex case]\label{ass:ADMM*}%
					In problem \eqref{eq:CP}, the following hold:
					\begin{enumeratass}
					\item\label{ass:f*}%
						\(\func{f}{\R^m}{\Rinf}\) is lsc, proper, and \(\mu_f\)-strongly convex.
					\item\label{ass:g*}%
						\(\func{\epicomp Bg}{\R^n}{\Rinf}\) is lsc, proper, and convex (e.g. when \(g\) is proper, convex and level bounded).
					\item\label{ass:ADMM*:relint}%
						\(b\in\relint(A\dom f+B\dom g)\).
					\end{enumeratass}
				\end{ass}
			\egroup

			\begin{fact}\label{thm:LL2DRE}%
				Starting from a triplet \((x,y,z)\in\R^m\times\R^p\times\R^n\), consider an \ref{ADMM}-update applied to problem \eqref{eq:CP} with relaxation \(\lambda\) and penalty \(\beta>0\).
				If
				\begin{enumerator}
				\item\label{thm:DRE2LL}
					either \cref{ass:ADMM} holds and \(\beta>L_{\epicomp Af}\),
				\item\label{thm:DRE2LL*}
					or \cref{ass:ADMM*} holds and \(\beta<\nicefrac{\mu_f}{\|A\|^2}\),
				\end{enumerator}
				then \(\DRE(s^+)=\LL_\beta(x^+,z^+,y^+)\) for \(s^+=Ax^+-y^+\nicefrac{}\beta\) and \(\gamma=\nicefrac1\beta\), with \(\DRE\) being the DRE associated to the equivalent problem formulation \eqref{eq:CP2P}.
			\end{fact}
			
			The following is a straightforward consequence of \cref{thm:DRE:SD,thm:LL2DRE}.
			\begin{cor}
				Starting from a triplet \((x^-,y^-,z^-)\in\R^m\times\R^p\times\R^n\), consider the \ref{ADMM} updates
				\(
					(x^-,y^-,z^-)
				{}\to{}
					(x,y,z)
				{}\to{}
					(x^+,z^+,y^+)
				\)
				with penalty \(\beta\) and relaxation \(\lambda\).
				Let \(r=Ax+Bz-b\) and \(\D[]\) be as in \eqref{eq:DRE:C}.
				\begin{enumerate}
				\item 
					If \cref{ass:ADMM} holds and \(\beta>L_{\epicomp Af}\), then
					\[
						\LL_\beta(x^+,z^+,y^+)
					{}\leq{}
						\LL_\beta(x,z,y)
						{}-{}
						\beta\D\|r\|^2,
					\]
					with the constant \(\D[]\) being strictly positive provided that
					\begin{equation}\label{eq:ADMM:beta}
						\beta
					{}>{}
						\begin{ifcases}
							L_{\epicomp Af}
						&
							f\text{ is convex}
						\\[5pt]
							\tfrac{2L_{\epicomp Af}}{2-\lambda}
						\otherwise.
						\end{ifcases}
					\end{equation}
				\item
					If \cref{ass:ADMM*} holds and \(\beta<\nicefrac{\mu_f}{\|A\|^2}\), then
					\[
						\LL_\beta(x^+,z^+,y^+)
					{}\geq{}
						\LL_\beta(x,z,y)
						{}+{}
						\beta\D*\|r\|^2,
					\]
					and \(\D*[]\) is strictly positive.
				\end{enumerate}
			\end{cor}

	\section{The linesearch algorithms}\label{sec:Algorithms}

		As shown in \cite{li2016douglas,li2017peaceman,themelis2020douglas}, both \ref{DRS} and \ref{ADMM} converge under mild assumptions that do not entail convexity of either functions.
		As a consequence, a wide range of nonsmooth and nonconvex problems can be addressed by iterating relatively simple operations.
		On the other hand, it is well known that even for convex problems the convergence can be prohibitively slow unless the problem is well scaled, which is rarely the case in practice.
		In contrast, fast local methods such as Newton-type exist that by exploiting higher-order information can suitably reshape the problem into a more convenient geometry.
		The major hindrance against their employment is their \emph{local} nature, in the sense that convergence is guaranteed only if the starting point is already close enough to a solution, on top of some regularity criteria around such solution.
		For this reason, fast local methods are typically paired with a \emph{linesearch} that globalizes convergence by ensuring a decrease condition on the cost or on a surrogate merit function.
		
		The purpose of \cref{alg:DRS,alg:ADMM}, presented in \cref{sec:AlgIntro}, is exactly to complement the global (subsequential) convergence and operational simplicity of \ref{DRS} and \ref{ADMM} with the fast local convergence of Newton-type schemes by means of a tailored linesearch.
		It should be noted that the proposed method differs considerably from classical linesearch strategies, which can only cope with \emph{directions of descent} and thus heavily hinge on differentiability requirements.
		Different from the convex setting of \cite{patrinos2014douglas}, where the envelope function is differentiable, in the setting dealt here it can only be guaranteed to be (locally Lipschitz) continuous (cf. \cref{thm:DREequiv}), and thus not suitable for a standard backtracking based, e.g., on the Armijo condition.
		For this reason, \cref{alg:DRS,alg:ADMM} here proposed adopt the novel linesearch protocol of the umbrella \emph{Continuous-Lyapunov Descent method} ({\sf CLyD}) \cite[\S4]{themelis2018proximal} already benchmarked with the PANOC \cite{stella2017simple} and NAMA \cite{stella2019newton} solvers, respectively based on the forward-backward splitting and the alternating minimization algorithm.
		
		The chosen linesearch allows us to discard differentiability requirements as it merely exploits continuity of the DRE and the sufficient decrease property.
		Before elaborating on the details in the next subsection, we first prove that \cref{alg:DRS,alg:ADMM} are linked through the same equivalence relating the underlying \ref{DRS} and \ref{ADMM} oracles which, for exposition clarity, will be referred to as the \emph{nominal steps}.
		
		\begin{prop}[Equivalence of \cref{alg:DRS} and \cref{alg:ADMM}]\label{thm:DRS=ADMM}%
			Suppose that \cref{ass:ADMM} [resp. \cref{ass:ADMM*}] holds and consider the iterates generated by \cref{alg:ADMM} with penalty \(\beta\) and relaxation \(\lambda\), starting from \((x^{-1},y^{-1},z^{-1})\).
			For each \(k\in\N\), let
			\begin{equation}\label{eq:ADMM2DRS}
				\begin{cases}[l >{{}}c<{{}} l]
					s^k
				&\coloneqq&
					Ax^k-y^k\nicefrac{}\beta
				{}={}
					b-Bz^{k-1}-y^{k-\nicefrac12}\nicefrac{}\beta
				\\
					u^k
				&\coloneqq&
					Ax^k
				\\
					v^k
				&\coloneqq&
					b-Bz^k.
				\end{cases}
			\end{equation}
			Then, \(\sDRS\coloneqq\epicomp Af\), \(\nsDRS\coloneqq\epicomp Bg(b-{}\cdot{})\), and \(\varphi\coloneqq\sDRS+\nsDRS\) satisfy \cref{ass:DRS} [resp. \cref{ass:DRS*}], and \(\seq{s^k,u^k,v^k}\) is a sequence generated by \cref{alg:DRS} with stepsize \(\gamma=\nicefrac1\beta\), relaxation \(\lambda\), starting from \(s^0\coloneqq b-Bz^{-1}-y^{-\nicefrac12}\nicefrac{}\beta\) and with the same sufficient decrease constant \(c\) and choice of directions \(\seq{d^k}\).
			Moreover, \(\DRE(s^k)=\LL_\beta(x^k,z^k,y^k)\) and \(Ax^k+Bz^k-b=u^k-v^k\) hold for every \(k\), hence both sequences \(\seq{r^k}\) and \(\seq{\tau_k}\) coincide in the two algorithms at every iteration.
			\begin{proof}
				Clearly, \cref{ass:DRS} [resp. \cref{ass:DRS*}] is satisfied if \cref{ass:ADMM} [resp. \cref{ass:ADMM*}] holds.
				It follows from \cref{thm:ADMM_DRS,thm:LL2DRE} that for any \(\tilde y,\tilde z\) such that
				\(
					s=-\tilde y\nicefrac{}\beta-B\tilde z+b
				\)
				one has
				\begin{subequations}\label{eq:operators}
					\begin{align}
						\DRS(s)
					{}={} &
						\set{(Ax,b-Bv)}[(x,y,z)\in\ADMM(\tilde y,\tilde z)],
					\\
						\DRE(s)
					{}={} &
						\LL_\beta(x,z,y)
					\quad
						\forall (x,y,z)\in\ADMM(\tilde y,\tilde z).
					\end{align}
				\end{subequations}
				In particular, the entire claim will follow once we prove that \(s^k\) and \((y^{k-\nicefrac12},z^{k-1})\) are related as in \eqref{eq:ADMM2DRS} for every \(k\).
				We proceed by induction.
				The case \(k=0\) follows from the definition of the initital iterate \(s^0\) as in the statement.
				Suppose now that the identity holds up to iteration \(k\), and let \(r^k=Ax^k+Bz^k-b\) be as in \cref{alg:ADMM}; then, observe that
				\newcommand{\hint}[1]{\text{\footnotesize(#1)}}%
				\begin{align*}
					\bar s^{k+1}
				{}={} &
					s^k+\lambda(v^k-u^k)
				&
					\hint{\cref{step:DRS:nominal}}
				\\
				{}={} &
					Ax^k-y^k\nicefrac{}\beta-\lambda(\overbracket*{Ax^k+Bz^k-b}^{r^k})
				&
					\hint{induction}
				\\
				{}={} &
					b-Bz^k-y^k\nicefrac{}\beta+(1-\lambda)r^k
				\\
				\numberthis\label{eq:bars}
				{}={} &
					b-Bz^k-\bar y^{k+\nicefrac12}\nicefrac{}\beta.
				&
					\hint{\cref{step:ADMM:nominal}}
				\end{align*}
				For \(\tau\in\R\), let
				\(
					\tilde s_\tau^{k+1}
				{}\coloneqq{}
					(1-\tau)\bar s^{k+1}+\tau(s^k+d^k)
				\)
				and
				\(
					\tilde y_\tau^{k+\nicefrac12}
				{}\coloneqq{}
					(1-\tau)\bar y^{k+\nicefrac12}
					{}+{}
					\tau(y^k-\beta(r^k+d^k))
				\)
				be the candidate updates with stepsize \(\tau\) at \cref{step:DRS:1/2} and \cref{step:ADMM:1/2}, respectively.
				We have
				\begin{align*}
					\tilde s_\tau^{k+1}
				{}={} &
					(1-\tau)\bar s^{k+1}+\tau(s^k+d^k)
				\\
				{}\overrel*{\eqref{eq:bars}}{} &
					(1-\tau)(b-\bar y^{k+\nicefrac12}\nicefrac{}\beta-Bz^k)
					{}+{}
					\tau
					\bigl(
						\overbracket*{
							Ax^k-y^k\nicefrac{}\beta
		% 					b-y^{k-\nicefrac12}\nicefrac{}\beta-Bz^{k-1}
						}^{
							s^k
						}
						{}+{}
						d^k
					\bigr)
				&
					\hint{induction}
				\\
				{}={} &
					b-Bz^k
					{}-{}
					\Bigl[
						(1-\tau)\bar y^{k+\nicefrac12}\nicefrac{}\beta
						{}+{}
						\tau(y^k\nicefrac{}\beta-Ax^k-Bz^k+b-d^k)
					\Bigr]
				\\
				{}={} &
					b-Bz^k
					{}-{}
					\Bigl[
						(1-\tau)\bar y^{k+\nicefrac12}\nicefrac{}\beta
						{}+{}
						\tau(y^k\nicefrac{}\beta-r^k-d^k)
					\Bigr]
				\\
				{}={} &
					b-Bz^k
					{}-{}
					\tilde y_\tau^{k+\nicefrac12}\nicefrac{}\beta.
				\end{align*}
				It then follows from \eqref{eq:operators} that
				\(
					\DRE(\tilde s_\tau^{k+1})
				{}={}
					\LL_\beta(\tilde x_\tau^{k+1},\tilde z_\tau^{k+1},\tilde y_\tau^{k+1})
				\)
				for any \((\tilde x_\tau^{k+1},\tilde z_\tau^{k+1},\tilde y_\tau^{k+1})\in\ADMM(\tilde y_\tau^{k+\nicefrac12},z^k)\).
				Combined with the fact that \(r^k=u^k-v^k\) holding by induction (that is, \(r^k\) is the same in both algorithms), we conclude that stepsize \(\tau=\tau_k\) is accepted at the \(k\)-th iteration of \cref{alg:ADMM} iff so happens in \cref{alg:DRS}.
				In particular, one has
				\[
					s^{k+1}
				{}={}
					\tilde s_{\tau_k}^{k+1}
				{}={}
					b-Bz^k
					{}-{}
					\tilde y_{\tau_k}^{k+\nicefrac12}
				{}={}
					b-Bz^k
					{}-{}
					y^{k+\nicefrac12},
				\]
				which completes the induction argument.
			\end{proof}
		\end{prop}

		\subsection{A continuity-based linesearch}\label{sec:LS}

			We now discuss the core of \cref{alg:DRS,alg:ADMM}, namely the linesearch strategy starting at \cref{step:DRS:LS,step:ADMM:LS}, respectively.
			Remarkably, not only is it flexible to \emph{any} update direction \(d\) (and not only of descent type), but it will be shown in \cref{thm:tau1} that it also accepts unit stepsize whenever \(d\) is ``good'', in a sense that will be made precise in \cref{sec:Superlinear}.
			In other words, the proposed linesearch is robust against the \emph{Maratos effect} \cite{maratos1978exact}, see also \cite[\S6.2]{izmailov2014newton}, a pathology typical of linesearch methods in nonsmooth optimization such as sequential quadratic programming that inhibits the achievement of fast convergence rates.
			
			Owing to the equivalence proven in \cref{thm:DRS=ADMM}, we limit the preliminary discussion to the \ref{DRS}-based \cref{alg:DRS}, as the rationale of the \ref{ADMM} counterpart uses the same arguments.
			Similarly, we may limit the discussion to the case in which \cref{ass:DRS} holds, as the complementary case of \cref{ass:DRS*} only differs from a change of sign in the DRE (cf. \cref{thm:selfdual}), covered by the initialization \(\pi=-1\), and a choice of stepsize \(\gamma\) and decrease constant \(c\) consistent with \cref{thm:DRE:SD*}, whence the different initialization prescribed in \cref{alg:DRS}.
			
			\begin{figure}
				\footnotesize
				\begin{minipage}[t]{0.325\linewidth}
					\def\fignum{1}%
					\includetikz[width=\linewidth]{LS}[LS1]
					\centering
					{\bf\textit{(a)}}~
					Nominal \ref{DRS} iterates
				\end{minipage}
				\hfill
				\begin{minipage}[t]{0.325\linewidth}
					\def\fignum{2}%
					\includetikz[width=\linewidth]{LS}[LS2]
					\centering
					{\bf\textit{(b)}}~
					Failed linesearch trials
				\end{minipage}
				\hfill
				\begin{minipage}[t]{0.325\linewidth}
					\def\fignum{3}%
					\includetikz[width=\linewidth]{LS}[LS3]
					\centering
					{\bf\textit{(c)}}~
					Stepsize accepted
				\end{minipage}
				\caption[Main steps of \cref{alg:DRS}]{%
					\emph{Main steps of \cref{alg:DRS}.}
					One call to the \ref{DRS} oracle at \(s\) yields the pair \((u,v)\) and the nominal \ref{DRS} update \(\bar s^+=s+\lambda(v-u)\).
					On the DRE, this implies a decrease by at least \(\tfrac1\gamma\C[]\|v-u\|^2\).
					Since \(\DRE\) is continuous and \(c\lneqq\C[]\), all points close enough to \(\bar s^+\) belong to the sublevel set \(\bigl[\DRE\leq\DRE(s)-c\|u-v\|^2\bigr]\) (cyan-shaded region).
					Therefore, \textbf{for any direction \(d\)}, all points close to \(\bar s^+\) in the line segment
					\(
						[\bar s^+,s+d]
					{}={}
						\set{(1-\tau)\bar s^++\tau(s+d)}[{\tau\in[0,1]}]
					\)
					belong to this set, hence the linesearch is accepted for small enough \(\tau\).
				}%
				\label{fig:LS}%
			\end{figure}

			Suppose that the current iterate is \(s\in\R^p\), and let \(d\in\R^p\) be an arbitrary candidate \emph{update direction} at \(s\).
			What \(d\) is, and how it is retrieved is irrelevant at the moment and will be discussed in detail in the dedicated \cref{sec:Directions}; suffice it to say that the choice of an update direction \(d\) represents our degree of freedom for extending \ref{DRS} while maintaining its (subsequential) convergence properties, and that ``ideally'' we would like to replace the \emph{nominal} \ref{DRS} update with the chosen \(s^+=s+d\), for we have reason to believe this choice will lead us closer to a solution.
			
			Let the nominal update be \(\bar s^+=s+\lambda(v-u)\) as in \cref{step:DRS:nominal}.
			% and suppose that \(v\neq u\) for otherwise a solution would be found and there would be no reason to investigate any further.
			Due to the sufficient decrease property on \(\DRE\) (cf. \cref{thm:DRE:SD}), it holds that
			\[
				\DRE(\bar s^+)
			{}\leq{}
				\DRE(s)
				{}-{}
				\tfrac1\gamma
				\C\|r\|^2
			\]
			where \(\C[]\) is as in \eqref{eq:DRE:C} and \(r=u-v\).
			However, nothing can be guaranteed as to whether \(\DRE(s+d)\) is also (sufficiently) smaller than \(\DRE(s)\) or not, nor can we hope to enforce the condition with a classical backtracking \(s+\tau d\) for small \(\tau>0\), as no notion of descent is known to \(\DRE\) (which is continuous but not necessarily differentiable, on top of the fact that the direction \(d\) is even arbitrary).
			Nevertheless, for any \(c\lneqq\C[]\), not only does \(\bar s^+\) satisfy the sufficent decrease with constant \(c\), but due to continuity of \(\DRE\) so do all the points around: loosely speaking,
			\begin{equation}\label{eq:SDa}
				\DRE(s')
			{}\leq{}
				\DRE(s)
				{}-{}
				\tfrac c\gamma
				\|r\|^2
			\quad
				\text{for all \(s'\) \emph{close} to \(\bar s^+\).}
			\end{equation}
			The idea is then to ``push'' the candidate update \(s+d\) towards the ``safe'' update \(\bar s^+\) until the \emph{relaxed} decrease condition \eqref{eq:SDa} holds.
			One way to do so is through a linesearch along the segment connecting the ``ideal'' update \(s+d\) and the ``safe'' nominal update \(\bar s^+\), as done in \cref{step:DRS:1/2}.
			The procedure is synopsized in \cref{fig:LS} for the toy problem of finding a point in the intersection of a line \(\ell\) and a circumference \(C\), cast in the form \eqref{eq:P} as
			\[
				\minimize_{s\in\R^2}{}
					\varphi(s)
				{}\coloneqq{}
					\underbrace{
						\tfrac12
						\dist^2(s,\ell)
					}_{\sDRS(s)}
					{}+{}
					\underbrace{
						\indicator_C(s)
						\vphantom{\tfrac12}
					}_{\nsDRS(s)}.
			\]
			The contour levels correspond to those of \(\DRE\); notice that, since \(\prox_{\gamma\sDRS}\equiv\id\) on \(\ell\) and that \(\argmin\varphi\subset\ell\), for this problem it holds that \(\argmin\varphi=\argmin\DRE\), cf. \cref{thm:DREequiv:inf}.

		\subsection{Iteration complexity}\label{sec:Complexity}

			In both the proposed algorithms, a step of the nominal method is required for the evaluation of \(\DRE\) or \(\LL_\beta\), where with ``nominal method'' we indicate \ref{DRS} for \cref{alg:DRS} and \ref{ADMM} for \cref{alg:ADMM}.
			Therefore, the number of nominal steps performed at iteration \(k\) corresponds to the number of backtrackings \(i_k\).
			In other words, the \(k\)-th iteration of the proposed algorithm is, in general, as expensive as \(i_k\) many nominal iterations.
			In order to bound complexity, a maximum number of backtrackings \(i_{\rm max}\) can be imposed, say \(i_{\rm max}=5\), so that whenever the linesearch condition fails \(i_{\rm max}\) many times, one can discard the direction \(d^k\) and proceed with a nominal update.
			Nevertheless, as explained in the following remark, whenever function \(\sDRS\) is (generalized) quadratic, not necessarily convex, the linearity of the \(u\)-update can conveniently be exploited to save computations in the linesearch.
			This is particularly appealing when, additionally, evaluating \(\prox_{\gamma\nsDRS}\) is cheap (projections onto simple sets, thresholding, \ldots), in which case each iteration is, at most, roughly twice as expensive as one of \ref{DRS}.
			Similar arguments also apply to the \ref{ADMM}-\cref{alg:ADMM}, and examples of this kind will be given in the Simulations \cref{sec:Simulations}.
			
			\begin{rem}[Exploiting linearity of \(\prox_{\gamma\sDRS}\)]\label{rem:linearprox}%
				If function \(\sDRS\) is quadratic (or \emph{generalized} quadratic, \ie quadratic plus the indicator of an affine subspace), then its proximal mapping is affine and the \(u\)-update at \cref{step:DRS:uv} can be expanded to
				\[
					u^{k+1}
				{}={}
					\prox_{\gamma\sDRS}
					\bigl[
						(1-\tau_k)\bar s^{k+1}+\tau_k(s^k+d^k)
					\bigr]
				{}={}
					(1-\tau_k)\prox_{\gamma\sDRS}(\bar s^{k+1})
					{}+{}
					\tau_k\prox_{\gamma\sDRS}(s^k+d^k).
				\]
				Then, for \(\tau_k=\nicefrac12\) instead of computing directly \(u^{k+1}\) one can evaluate \(\prox_{\gamma\sDRS}(\bar s^{k+1})\) and obtain \(u^{k+1}\) by linear combination, and similarly all subsequent trials for smaller values of \(\tau_k\) will not require any additional evaluation of \(\prox_{\gamma\sDRS}\).
				In this case, the number of evaluations of \(\prox_{\gamma\sDRS}\) at each iteration of the algorithm is at most \(2\).
			
				Similarly, at most two evaluations of \(\sDRS\) (needed for computing the value of \(\DRE\)) will be enough.
				In fact, denoting \(\bar u^{k+1}\coloneqq\prox_{\gamma\sDRS}(\bar s^{k+1})\) and \(u_0^{k+1}\coloneqq\prox_{\gamma\sDRS}(s^k+d^k)\), one has that
				\[
					\ell(\tau)
				{}\coloneqq{}
					\sDRS((1-\tau)\bar u^{k+1}+\tau u_0^{k+1})
				\]
				is a one-dimensional quadratic function, hence \(\ell(\tau)=a+b\tau+c\tau^2\) with \(a=\ell(0)\), \(b=\ell'(0)\), and \(c=\ell(1)-a-b\), \ie
				\[
					\begin{cases}
						a
					& {}={}
						\sDRS(\bar u^{k+1})
					\\
						b
					& {}={}
						\frac1\gamma\innprod*{\bar s^{k+1}-\bar u^{k+1}}{u_0^{k+1}-\bar u^{k+1}}
					\\
						c
					& {}={}
						\sDRS(u_0^{k+1})-a-b.
					\end{cases}
				\]
				The expression of \(b\) uses the fact that \(\bar s^{k+1}-\bar u^{k+1}\in\gamma\partial\sDRS(\bar u^{k+1})\), cf. \eqref{eq:prox:subdiff}.
				Consequently, once \(\sDRS(u_0^{k+1})\) and \(\sDRS(\bar u^{k+1})\) are computed, all the needed values of \(\ell(\tau_k)\) can be retrieved at negligible cost.
			\end{rem}

		\subsection{Choice of direction}\label{sec:Directions}

				Although the proposed algorithmic framework is robust to any choice of directions \(d^k\), its efficacy is greatly affected by the specific selection.
				This subsection provides an overview on some convenient choices of update directions \(d^k\).
				Specifically, we propose a generalization of Nesterov extrapolation \cite{nesterov1983method} that is suited for nonconvex problems (although merely heuristical, without optimal convergence guarantees), and then discuss three popular Newton-type schemes: (modified) Broyden \cite{broyden1965class,powell1970numerical}, BFGS \cite{broyden1970convergence,fletcher1970new,goldfarb1970family,shanno1970conditioning}, and Andreson acceleration \cite{anderson1965iterative}.
				Although true higher-order information, when available, can also be considered, we prefer to limit our overview to methods that preserve the simple oracle of the original \ref{DRS} and \ref{ADMM}.
				In the convex case, the interested reader can find an extensive collection of generalized Jacobians of proximal mappings in \cite[\S15.6]{themelis2019acceleration}, useful for deriving directions based on \emph{linear Newton approximation schemes} \cite[\S7.5.1]{facchinei2003finite}.
			
				\begin{rem}[Recommended directions]\label{rem:directions}%
					While all of the directions choices outlined in this section enable the subsequential convergence properties of \cref{alg:DRS} and \cref{alg:ADMM} (cf. \cref{sec:Convergence}), some choice may perform better than the others on specific problems.
					In our experience, and as confirmed by the numerical evidence in \cref{sec:Simulations}, Broyden and L-BFGS directions consistently provided fast convergence: the Broyden method is supported by the theory, in that superlinear convergence can be proved under mild assumptions (cf. \cref{sec:Superlinear}), while L-BFGS clearly scales better with the problem dimension, given its limited-memory nature.
					On the other hand, Nesterov directions proved effective only in case \(\sDRS\) (or \(f\)) is convex, while Anderson acceleration directions only performed well in one example in our experiments.
					These are purely empirical observations, and a thorough understanding of when (and why) these directions perform well should be the subject of future investigation.
				\end{rem}
				
				\subsubsection{Nesterov acceleration}\label{sec:Nesterov}%
					It was shown in \cite{patrinos2014douglas} that when \(\nsDRS\) is convex and \(\sDRS\) is convex quadratic, then \(\DRE\) is convex and continuously differentiable for \(\gamma<\nicefrac{1}{L_{\sDRS}}\).
					This enabled the possibility, for this specific case, to extend the employment of the optimal Nesterov acceleration techniques \cite{nesterov1983method} to \ref{DRS}.
					By using duality arguments, this fact was extended in \cite{pejcic2016accelerated} where an \emph{accelerated} \ref{ADMM} scheme is proposed for problems in the form \eqref{eq:CP} with \(B=\I\), \(\sADMM\) strongly convex quadratic, and \(\nsADMM\) convex.
					
					Although not supported by the theory, extensive numerical evidence suggests that such extrapolations perform quite effectively regardless of what function \(\nsDRS\) or \(\nsADMM\) in problems \eqref{eq:P} and \eqref{eq:CP} are, while convexity of $\sDRS$ and $\sADMM$ seems instead to play an important role, as mentioned in \cref{rem:directions}.
					The main limitation in a direct employment of the acceleration is that convergence itself is not guaranteed.
					However, thanks to the arbitrarity of \(d^k\) in \cref{thm:DRS_subseq,thm:ADMM_subseq}, we can enforce these updates into \cref{alg:DRS,alg:ADMM} and thus obtain a (damped, monotone) extrapolation that is guaranteed to be globally (subsequentially) convergent.
					\begin{proofitemize}
					\item
						Fast \ref{DRS}: in \cref{alg:DRS}, start with \(d^0=-\lambda r^0\), and for \(k\geq1\) select
						\[
							d^k
						{}\coloneqq{}
							\tfrac{k-1}{k+2}
							(\bar s^{k+1}-\bar s^k)
							{}-{}
							\lambda r^k
						\]
					\item
						Fast \ref{ADMM}: in \cref{alg:ADMM}, start with \(d^0=-\lambda r^0\), and for \(k\geq1\) select
						\[
							d^k
						{}\coloneqq{}
							-\lambda r^k
							{}-{}
							\tfrac{k-1}{k+2}
							\bigl(
								Bz^k-Bz^{k-1}
								{}+{}
								(\bar y^{k+\nicefrac12}-\bar y^{k-\nicefrac12})\nicefrac{}\beta
							\bigr).
						\]
					\end{proofitemize}
			
					These extensions differ from the approach proposed in \cite{li2015accelerated} for the proximal gradient method, as we do not discard the candidate fast direction when sufficient decrease is not satisfied but rather dampen it with a backtracking.

				\subsubsection{Quasi-Newton methods}
					The termination criterion for both \cref{alg:DRS,alg:ADMM} is based on (the norm of) the fixed-point residual of the underlying splitting schemes, namely
					\(
						u-v
					{}={}
						\tfrac1\lambda(\bar s-\bar s^+)
					\)
					for \ref{DRS}, which corresponds to
					\(
						Ax+Bz-b
					\)
					in \ref{ADMM}.
			% 		It is well known that when both \(\sDRS\) and \(\nsDRS\) are convex, having \(u^k-v^k=0\) is equivalent to having \(u^k\in\argmin\sDRS+\nsDRS\).
					Under some assumptions such as \emph{prox-regularity} \cite[\S13.F]{rockafellar2011variational}, eventually the updates \(s\mapsto s^+\) are uniquely determined, that is, the inclusion \(v\in\prox_{\gamma\nsDRS}(2u-s)\) in \ref{DRS} becomes an equality.
					It then turns out that one ends up solving a system of nonlinear equations, namely finding \(s\) such that \(s-s^+=0\) in \ref{DRS}, and similarly for the residual \(Ax+Bz-b=0\) in \ref{ADMM}.
					
					Parallel to what done in \cite{themelis2018forward} for (nonconvex) forward-backward splitting, we may then consider directions stemming from fast methods for nonlinear equations, namely
					\(
						d^k=-H_k\Res[][]
					\),
					where \(\Res[][]\) is the fixed-point residual map (\(\Res(s)=u-v\) for \ref{DRS} and \(\Res[\beta][admm](y^{\nicefrac-2},z^{-1})=Ax+Bz-b\) for \ref{ADMM}), and \(H_k\) is some approximation to the inverse of its (generalized) Jacobian.
					To maintain the simplicity of the original \ref{DRS} and \ref{ADMM}, this can be done efficiently by means of quasi-Newton methods, which starting from an invertible matrix \(H_0\) perform low-rank updates based on available quantities.
					Such quantities are pairs of vectors \((p_k,q_k)\), where \(p_k\) is the difference between consecutive iterates and \(q_k\) is the difference of the fixed-point residuals.
					Namely,
					\begin{subequations}\label{subeq:pq}%
						\begin{proofitemize}
						\item
							quasi-Newton \ref{DRS}: in \cref{alg:DRS} use
							\begin{equation}\label{eq:DRS:pq}
								d^k
							{}={}
								-H_kr^k
							\quad\text{and}\quad
								\begin{cases}[l@{{}={}}l]
									p_k &
			% 						s^{k+1}-s^k
			% 						s_0^{k+1}-s^k=
									d^k
								\\
									q_k &
			% 						r^{k+1}-r^k
									r_0^{k+1}-r^k,
								\end{cases}
							\end{equation}
							where \(r_0^{k+1}=u_0^{k+1}-v_0^{k+1}\) with \((u_0^{k+1},v_0^{k+1})\in\DRS(s^k+d^k)\) is the residual computed in the first linesearch trial.%
						\item
							quasi-Newton \ref{ADMM}: start with \(H_0=\mu\I\) for some \(\mu>0\), and in \cref{alg:ADMM} use
							\begin{equation}\label{eq:ADMM:pq}
								d^k
							{}={}
								-H_kr^k
							\quad\text{and}\quad
								\begin{cases}[l@{{}={}}l]
									p_k &
			% 						Bz^k-Bz^{k-1}
			% 						{}+{}
			% 						(y^{k+\nicefrac12}-y^{k-\nicefrac12})\nicefrac{}\beta
									d^k
								\\
									q_k & r_0^{k+1}-r^k,
								\end{cases}
							\end{equation}
							where \(r_0^{k+1}\) is the residual computed in the first linesearch trial, as in the \ref{DRS} case.
			% 				\(r_0^{k+1}=Ax_0^{k+1}+Bz_0^{k+1}-b\) with \((u_0^{k+1},v_0^{k+1})\in\DRS(s^k+d^k)\) is the residual computed in the first linesearch trial.%
						\end{proofitemize}
						In textbook quasi-Newton methods employed in smooth optimization, where the objective is to find a point where the gradient of the cost function is zero, \(q_k\) would amount to the difference of gradient values, see \eg \cite[pp. 24--26]{nocedal2006numerical}.
						Here, where the objective is to find a zero of the residual operator, \(q_k\) is the difference of residual values.
						Nevertheless, a naive adaptation from the smooth case would suggest considering \(p_k=s^{k+1}-s^k\) and \(q_k=r^{k+1}-r^k\), whereas here the values obtained from the first linesearch iteration trial are considered, regardless of whether stepsize \(\tau_k=1\) is accepted or not.
						As will be clear in the superlinear convergence proof of \cref{thm:Broyden}, the one considered here is a more educated adaptation that complies with the proposed innovative linesearch.
					\end{subequations}
					
					We will now list a few update rules for \(H_k\) based on the indicated pairs \((p_k,q_k)\).
					\begin{directions}
					\item[BFGS]
						Start with \(H_0\succ 0\) and update as follows:
						\[
							H_{k+1}
						{}={}
							H_k
							{}+{}
							\frac{\innprod{p_k}{q_k}+\innprod{H_kq_k}{q_k}}{(\innprod{p_k}{q_k})^2}p_k\trans{p_k}
							{}-{}
							\frac{H_kq_k\trans{s_k}+s_k\trans{q_k}H_k}{\innprod{p_k}{q_k}}.
						\]
						Whenever \(\innprod{p_k}{q_k}\leq 0\), one can either set \(H_{k+1}=H_k\) or use a different vector \(p_k\) as proposed in \cite{powell1978fast}.
						The limited-memory variant L-BFGS \cite[Alg. 7.4]{nocedal2006numerical}, which does not require storage of full matrices \(H_k\) or matrix-vector products but only storage of the last few pairs and scalar products, can conveniently be considered.
						
						Although very well performing in practice, to the best of our knowledge fast convergence of BFGS can only be shown when the Jacobian of \(\Res[][]\) at the limit point is symmetric, which hardly ever holds in our framework.
						We suspect, however, that the well performance of BFGS derives from the observation that, when it exists, the Jacobian of \(\Res[][]\) at a local minimum is diagonalizable and has all positive eigenvalues, and is thus similar to a symmetric positive semidefinite matrix (cf. \eqref{eq:FBEmin} and \cite[Thm. 4.11]{themelis2018forward}).
			% 			\todo{Powell BFGS trick?}
					\item[Modified Broyden]\label{sec:Broyden}%
						Fix \(\bar\vartheta\in(0,1)\), \eg \(\bar\vartheta=0.2\), an invertible matrix \(H_0\), and update as follows:
						\begin{subequations}\label{subeq:Broyden}
						\begin{equation}
							H_{k+1}
						{}={}
							H_k
							{}+{}
							\frac{
								p_k-H_kq_k
							}{
								\innprod{p_k}{(\nicefrac{1}{\vartheta_k}-1)p_k+H_kq_k}
							}
							\trans{p_k}H_k
						\end{equation}
						where
						\begin{equation}\label{eq:PowellTheta}
							\vartheta_k
						{}\coloneqq{}
							\begin{cases}[l@{~~\text{if }}l]
									1
								&
									|\delta_k| \geq \bar\vartheta
								\\
									\frac{
										1-\sign(\delta_k)
										\bar\vartheta
									}{
										1-\delta_k
									}
								&
									|\delta_k| < \bar\vartheta
							\end{cases}
						\qquad\text{and}\quad
							\delta_k
						{}\coloneqq{}
							\frac{
								\innprod*{H_kq_k}{p_k}
							}{
								\|p_k\|^2
							},
						\end{equation}
						\end{subequations}
						with the convention that \(\sign 0=1\).
						The original Broyden formula \cite{broyden1965class} corresponds to \(\vartheta_k\equiv 1\), while \(\vartheta_k\) as in \eqref{eq:PowellTheta} ensures nonsingularity of all matrices \(H_k\) \cite{powell1970numerical}.
					\item[Anderson acceleration]
						Fix a buffer size \(m\geq 1\) and start with \(H_0=\I\).
						For \(k\geq1\), let
						\[
							H_k
						{}={}
							\I
							{}+{}
							(\mathcal P_k-\mathcal Q_k)
							(\trans{\mathcal Q_k}\mathcal Q_k)^{-1}
							\trans{\mathcal Q_k},
						\]
						where the columns of matrix \(\mathcal P_k\) are the last vectors \(p_{k-M},\cdots,p_{k-1}\) and those of \(\mathcal Q_k\) are the last vectors \(q_{k-M},\cdots,q_{k-1}\), with \(M=\min\set{k,m}\).
						If \(\mathcal Q_k\) is not full-column rank, for \(x\in\R^M\) the product \((\trans{\mathcal Q_k}\mathcal Q_k)^{-1}x\) is meant in a least-square sense.
						This is a limited-memory scheme, as it requires only the storage of few vectors and the solution of a small \(M\times M\) linear system.
						Anderson acceleration originated in \cite{anderson1965iterative}; here we use the interpretation well explained in \cite{fang2009two} of \emph{(inverse) multi-secant} update: \(H_k\) is the matrix closest to the identity (in the Frobenius norm) among those satisfying \(H_k\mathcal Q_k=\mathcal P_k\).
					\end{directions}

		\subsection{Adaptive variants}\label{sec:Adaptive}

			One drawback of \cref{alg:DRS,alg:ADMM} is that both the stepsize \(\gamma\) in the former and the penalty \(\beta\) in the latter have to be chosen offline based either on a Lipschitz constant or on a strong convexity modulus.
			In practice, the estimation of these quantities is often challenging and prone to yield very conservative approximations, dooming the algorithms to slow convergence in early iterations (that is, in the globalization stage when the effect of the fast local directions is not triggered yet).
			Moreover, even when such constants are known the algorithms may potentially work also with less conservative estimates which better reflect the local geometry.
			
			In order to circumvent these issues and allow for out-of-the-box implementations, we may resort to the adaptive variants of the \ref{DRS} and \ref{ADMM} oracles as described in \cite[\S4.1 and \S5.3]{themelis2020douglas}, where \(\gamma\) and \(\beta\) are tuned online in such a way to ensure the needed sufficient decrease conditions and preserve convergence.

			\subsubsection{Adaptive \texorpdfstring{\cref{alg:DRS}}{Algorithm \ref*{alg:DRS}}}
				For the \ref{DRS}-based \cref{alg:DRS}, it suffices to initialize \(\gamma\) according to an estimate of \(L_{\sDRS}\) or \(\mu_{\sDRS}\), and simply add the following routine after \cref{step:DRS:nominal}:
			
				\begin{algorithmic}[1]\small
					\makeatletter
						\renewcommand\theALG@line{\ref{step:DRS:nominal}a}%
					\makeatother
					\renewcommand\alglinenumber[1]{\ref*{step:DRS:nominal}\fillwidthof[c]ba:}%
					\State
						Evaluate \(\DRE(\bar s^{k+1})\) using \((\bar u^{k+1},\bar v^{k+1})\in\DRS(\bar s^{k+1})\) as in \eqref{eq:DRE}
					\makeatletter
						\renewcommand\theALG@line{\ref{step:DRS:nominal}b}%
					\makeatother
					\renewcommand\alglinenumber[1]{\ref*{step:DRS:nominal}b:}%
					\If{~
						\(
							\pi\DRE(\bar s^{k+1})
						{}\geq{}
							\pi\DRE(s^k)-\tfrac c\gamma\|r^k\|^2
						\)
						~{\bf or}~
						\(
							\pi\DRE(\bar s^{k+1})
						{}<{}
							\pi\varphi_{\text{\sc lb}}
						\)
					~}\label{state:DRS:adaptive}%
						\Statex\hspace*{\algorithmicindent}
							\(\gamma\gets2^{-\pi}\gamma\),~
							\((u^k,v^k)\in\DRS(s^k)\),~
							recompute \(\DRE(s^k)\) and go to \cref{step:DRS:1}.
					\EndIf{}
				\end{algorithmic}
			
				Here, recall that \(\pi\in\set{\pm1}\) is set at algorithm initialization according to whether \cref{ass:DRS} or \cref{ass:DRS*} holds.
				The only new term is \(\varphi_{\text{\sc lb}}\), to be set offline equal to a known quantity that lower bounds \(\inf\varphi\), which in practice is typically easily estimable.
				Its role is however a pure technicality that the not-too-fussy user can neglect; the interested reader can instead find the reasoning for the additional condition it enforces in \cite[\S4.1]{themelis2020douglas}.
				As documented in the reference, this backtracking on \(\gamma\) can happen only a finite number of times, as eventually the condition at \cref{state:DRS:adaptive} is never satisfied.
				In fact, under \cref{ass:DRS}, by halvening \(\gamma\) enough yet finitely many times, the value will eventually fall under the (unknown) threshold prescribed by \eqref{eq:DRE:gamma}.
				Similarly, under \cref{ass:DRS*}, by doubling \(\gamma\) finitely many times it will eventually exceed the (unknown) threshold \(\nicefrac{1}{\mu_{\sDRS}}\) prescribed in \cref{thm:DRE:SD*}.
				Either way, \cref{thm:DRE:SD,thm:DRE:SD*} will guarantee that no more changes to \(\gamma\) will be triggered after that point.

			\subsubsection{Adaptive \texorpdfstring{\cref{alg:ADMM}}{Algorithm \ref*{alg:ADMM}}}
				The same arguments as in the previous paragraph apply to the \ref{ADMM}-based \cref{alg:ADMM}, in which case it suffices to initialize \(\beta\) according to an estimate of \(L_{\epicomp Af}\) or \(\mu_{\epicomp Af}\), setting offline a lower bound
				\(
					\Phi_{\text{\sc lb}}
				{}\leq{}
					\inf\set{f(x)+g(z)}[Ax+Bz=b]
				\)
				(if known), and adding the following check after \cref{step:ADMM:nominal}:
			
				\begin{algorithmic}[1]\small
					\makeatletter
						\renewcommand\theALG@line{\ref{step:ADMM:nominal}a}%
					\makeatother
					\renewcommand\alglinenumber[1]{\ref*{step:ADMM:nominal}\fillwidthof[c]ba:}%
					\State
						Evaluate~
						\(
							\LL_\beta(\bar x^{k+1},\bar z^{k+1},\bar y^{k+1})
						\),
						~where~
						\(
							(\bar x^{k+1},\bar z^{k+1},\bar y^{k+1})
						{}\in{}
							\ADMM(\bar y^{k+\nicefrac12},z^k)
						\)
					\makeatletter
						\renewcommand\theALG@line{\ref{step:ADMM:nominal}b}%
					\makeatother
					\renewcommand\alglinenumber[1]{\ref*{step:ADMM:nominal}b:}%
					\If{~
						\(
							\pi\LL_\beta(\bar x^{k+1},\bar z^{k+1},\bar y^{k+1})
						{}\geq{}
							\pi\LL_\beta(x^k,z^k,y^k)-\beta c\|r^k\|^2
						\)
						~{\bf or}~
						\(
							\pi\LL_\beta(\bar x^{k+1},\bar z^{k+1},\bar y^{k+1})
						{}<{}
							\pi\Phi_{\text{\sc lb}}
						\)
					~}\label{state:ADMM:adaptive}%
						\Statex\hspace*{\algorithmicindent}
							\(\beta\gets2^\pi\beta\),~
							\((x^k,z^k,y^k)\in\ADMM(y^{k-\nicefrac12},z^{k-1})\),
							recompute \(\LL_\beta(x^k,z^k,y^k)\) and go to \cref{step:ADMM:1}.
					\EndIf{}
				\end{algorithmic}

	\section{Convergence results}\label{sec:Convergence}

			This section is dedicated to the convergence properties of the proposed \cref{alg:DRS,alg:ADMM}.
			We begin by addressing their well definedness, namely, that the iterations cannot get stuck in an infinite backtracking loop at \cref{step:DRS:tau/2,step:ADMM:tau/2}.
			We also show that for any strictly positive tolerance \(\varepsilon>0\) the termination criterion is satisfied after finitely many iterations, and provide properties of the output quantities.
			
			\begin{thm}[Well definedness and finite termination of \cref{alg:DRS}]\label{thm:DRS}%
				Suppose that either \cref{ass:DRS} or \cref{ass:DRS*} is satisfied.
				Then, the following hold for the iterates generated by \cref{alg:DRS}:%
				\begin{enumerate}
				\item\label{thm:DRS:LS}
					At every iteration the number of backtrackings at \cref{step:DRS:tau/2} is finite (regardless of whether \(i_{\rm max}\) is finite or not).
				\item\label{thm:DRS:K}
					The algorithm terminates in
					\(
						K
					{}\leq{}
						\frac{1}{c\gamma\varepsilon^2}
						|\DRE(s^0)-\min\varphi|
					\)
					iterations.
				\item\label{thm:DRS:output}
					The last iterate of the algorithm yields
					\begin{itemize}[leftmargin=*]
					\item
						a point \(z\coloneqq v^K\) satisfying
						\(
							\dist(0,\hat\partial\varphi(z))
						{}\leq{}
							2\varepsilon
						\)
						if \cref{ass:DRS} holds,
					\item
						a triplet \((x,y,z)\coloneqq (u^K,\gamma^{-1}(u^K-s^K),v^K)\) satisfying the approximate KKT
						\[
							-y\in\partial\sDRS(x),\quad
							\dist(y,\partial\nsDRS(z))\leq\varepsilon,\quad
							\|x-z\|\leq\gamma\varepsilon,
						\]
						if \cref{ass:DRS*} holds.
					\end{itemize}
				\end{enumerate}
				\begin{proof}
					We first consider the case in which \cref{ass:DRS} holds, so that \(\pi=1\) and consequently \(\pi\DRE=\DRE\).
					\begin{proofitemize}
					\item\ref{thm:DRS:LS}~
						Testing the condition at \cref{step:DRS:LS} assumes \(\|r^k\|>0\), for otherwise the entire algorithm would have stopped at \cref{step:DRS:1}.
						Let \(\C[]=\C\) be as in \eqref{eq:DRE:C}, so that the nominal \ref{DRS}-update \(\bar s^{k+1}\) satisfies
						\[
							\pi\DRE(\bar s^{k+1})
						{}\leq{}
							\pi\DRE(s^k)
							{}-{}
							\tfrac{\C[]}{\gamma}\|r^k\|^2,
						\]
						as shown in \cref{thm:DRE:SD}.
						Since \(c<\C[]\) by algorithm initialization, one has
						\[
							\pi\DRE(s^k)
							{}-{}
							\tfrac c\gamma\|r^k\|^2
						{}>{}
							\pi\DRE(s^k)
							{}-{}
							\tfrac{\C[]}{\gamma}\|r^k\|^2.
						\]
						Continuity of \(\DRE\) (cf. \cref{thm:DREequiv}) at \(\bar s^{k+1}\) thus entails the existence of \(\epsilon>0\) such that
						\(
							\pi\DRE(s)
						{}\leq{}
							\pi\DRE(s^k)
							{}-{}
							\tfrac c\gamma\|r^k\|^2
						\)
						holds for every point \(s\) \(\epsilon\)-close to \(\DRE(s^k)\).
						Since \((1-\tau_k)\bar s^{k+1}+\tau_k(s^k+d^k)\to\bar s^{k+1}\) as \(\tau_k\to 0\), by halvening enough times \(\tau_k\) at \cref{step:DRS:tau/2},
						\begin{itemize}
						\item
							either the candidate point \(s^{k+1}\) is eventually \(\epsilon\)-close to \(\bar s^{k+1}\) so that the needed condition at \cref{step:DRS:LS} holds,
						\item
							or the maximum number of backtrackings \(i_{\rm max}\) is reached (provided that \(i_{\rm max}\) is finite), in which case \(s^{k+1}=\bar s^{k+1}\).
						\end{itemize}
						Either way, only a finite number of backtrackings is performed, and the inequality
						\begin{equation}\label{eq:DRS:SD}
							\pi\DRE(s^{k+1})
						{}\leq{}
							\pi\DRE(s^k)
							{}-{}
							\tfrac c\gamma\|r^k\|^2
						\end{equation}
						holds for every \(k\in\N\).
					\item\ref{thm:DRS:K}~
						By combining the decrease condition \eqref{eq:DRS:SD} with the termination criterion at \cref{step:DRS:1}, we have that \(K\) iterations of the algorithm result in a decrease of the DRE by at least
						\[\textstyle
							\pi\DRE(s^{K+1})
						{}\leq{}
							\pi\DRE(s^0)
							{}-{}
							\tfrac c\gamma
							\sum_{k=0}^K\|r^k\|^2.
						\]
						Since \(\pi\min\varphi=\min\pi\DRE\) (cf. \cref{thm:DREequiv:inf}) and \(\|r^k\|\geq\gamma\varepsilon\) for \(k<K\), we conclude that
						\begin{equation}\label{eq:telescope}
						\textstyle
							\pi\DRE(s^0)
							{}-{}
							\pi\min\varphi
						{}\geq{}
							\tfrac c\gamma
							\sum_{k=0}^{K-1}\|r^k\|^2
						{}\geq{}
							cK\gamma\varepsilon^2,
						\end{equation}
						resulting in the claimed bound on \(K\).
						Finally, the bound on \(\dist\bigl(0,\hat\partial\varphi(z)\bigr)\) follows from \cite[Thm. 4.3(ii)]{themelis2020douglas}.
					\end{proofitemize}
					Suppose now that \cref{ass:DRS*} holds.
					Then, since \(\pi=-1\), it follows from \cref{thm:selfdual} that \(\pi\DRE(s)=\DRE*(s_*)\), where \(\gamma_*=\nicefrac1\gamma\) and \(s_*=-\nicefrac s\gamma\).
					Moreover, one has that \(u-v=\frac{u_*-v_*}{\gamma_*}\), where \(u_*\) and \(v_*\) are as in \cref{thm:selfdual}.
					The previous arguments can thus be replicated in terms of the dual formulation \eqref{eq:D} to prove assertion \ref{thm:DRS:LS} and the bound on the number of iterations \(K\).
			
					To conclude, let the triplet \((x,y,z)=(u^K,\frac{u^K-s^K}{\gamma},v^K)\) be as in the statement.
					From \eqref{eq:prox:subdiff} one has that \(-y\in\partial\sDRS(x)\) and that \(y+\tfrac1\gamma r^K\in\partial\nsDRS(z)\), hence that \(\dist(y,\partial\nsDRS(z))\leq\varepsilon\), owing to the termination criterion \(\|r^K\|\leq\gamma\varepsilon\).
					Finally, \(\|x-z\|=\|r^K\|\leq\varepsilon\), completing the proof.
				\end{proof}
			\end{thm}

			\begin{thm}[Well definedness and finite termination of \cref{alg:ADMM}]\label{thm:ADMM}%
				Suppose that either \cref{ass:ADMM} or \cref{ass:ADMM*} is satisfied.
				Then, the following hold for the iterates generated by \cref{alg:ADMM}:%
				\begin{enumerate}
				\item\label{thm:ADMM:LS}
					At every iteration the number of backtrackings at \cref{step:ADMM:tau/2} is finite (regardless of whether \(i_{\rm max}\) is finite or not).
				\item\label{thm:ADMM:K}
					The algorithm terminates in
					\(
						K
					{}\leq{}
						\frac{\beta}{c\varepsilon^2}|\LL_\beta(x^0,z^0,y^0)-\min\Phi|
					\)
					iterations and yields a triplet \((x,y,z)\coloneqq(x^K,y^K,z^K)\) satisfying the approximate KKT conditions
					\[
						\|Ax+Bz-b\|\leq\nicefrac\varepsilon\beta,
					\quad
						-\trans Ay\in\hat\partial\sADMM(x),
					\quad
						\dist(-\trans By,\hat\partial\nsADMM(z))\leq\|B\|\varepsilon.
					\]
				\end{enumerate}
				\begin{proof}
					Assertion \ref{thm:ADMM:LS} and the bound on the number of iterations \(K\) follow from \cref{thm:DRS}, owing to the equivalence of \cref{alg:DRS,alg:ADMM} stated in \cref{thm:DRS=ADMM}.
					In turn, the conditions
					\(-\trans Ay\in\hat\partial\sADMM(x)\)
					and
					\(\dist(-\trans By,\hat\partial\nsADMM(z))\leq\|B\|\varepsilon\)
					follow from \cref{thm:ADMM:-Ay,thm:ADMM:-By}, since for every \(k\) and independently of the choice of the direction \(d^k\) the triplet \((x^k,y^k,z^k)\) is the result of an \ref{ADMM}-step with penalty \(\beta\).
				\end{proof}
			\end{thm}

		\subsection{Subsequential convergence}

			The remainder of the section will focus on more theoretical aspects of the algorithms, and specifically on asymptotic behaviors.
			To this end, an idealistic tolerance \(\varepsilon=0\) will be considered, so that the algorithms may run infinitely many iterations.
			In this first subsection, without imposing additional requirements other than either one among
			{%
				\renewcommand{\creflastconjunction}{, or }%
				\cref{ass:DRS,ass:DRS*,,ass:ADMM,,ass:ADMM*},%
			}
			we show that every limit point of the generated sequences is stationary, and also give a sufficient condition ensuring boundedness.
			Under additional assumptions, we will later demonstrate and rigorously define the speed-up triggered by suitable choices of directions, which will ultimately be backed up by numerical evidence in \cref{sec:Simulations}.
			Some results are based on auxiliary material presented in \cref{sec:Appendix}.
			
			\begin{thm}[Subsequential convergence of \cref{alg:DRS}]\label{thm:DRS_subseq}%
				Suppose that either \cref{ass:DRS} or \cref{ass:DRS*} is satisfied, and consider the (possibly infinite) iterates generated by \cref{alg:DRS} with tolerance \(\varepsilon=0\).
				Then, the sequence of squared residuals \(\seq{\|r^k\|^2}\) has finite sum.
				Moreover,
				\begin{itemize}[label={}, leftmargin=0pt]
				\item{\sc under \cref{ass:DRS}:}
					\begin{enumerate}
					\item\label{thm:DRS:omega}%
						The sequences \(\seq{u^k}\) and \(\seq{v^k}\) have the same cluster points, all of which are stationary for \(\varphi\) and on which \(\varphi\) and \(\DRE\) have the same (constant) value, this being the limit of the monotonic sequence \(\seq{\DRE(s^k)}\).
					\item\label{thm:DRS:bounded}%
						If \(\varphi\) is level bounded, then \(\seq{s^k,u^k,v^k}\) remains bounded.
					\end{enumerate}
				\item{\sc under \cref{ass:DRS*}:}
					\begin{enumerate}[resume]
					\item\label{thm:DRS*:omega}%
						\(\DRE(s^k)\nearrow\min\varphi\) as \(k\to\infty\), and \(\seq{u^k}\) and \(\seq{v^k}\) converge to the unique (global) minimizer \(x_\star\) of \(\varphi\).
					\item\label{thm:DRS*:bounded}%
						If \(0\in\interior(\dom\sDRS-\dom\nsDRS)\) (equivalently, if the dual cost \(\psi\) as in \eqref{eq:D} is level bounded), then \(\seq{s^k}\)  is bounded.
					\end{enumerate}
				\end{itemize}
				\begin{proof}
					To rule out trivialities we may assume that the stopping criterion \(\|r^k\|=0\) is never reached, hence that the algorithm generates infinitely many iterates.
					That the squared residuals \(\seq{\|r^k\|^2}\) have finite sum follows by letting \(K\to\infty\) in \eqref{eq:telescope}.
					In turn, this implies that \(u^k-v^k\to0\) as \(k\to\infty\).
					We now consider the two separate cases.%
					\begin{proofitemize}
					\item{\sc \cref{ass:DRS}.}~
						Having shown the vanishing of the residual, the proof is identical to that of \cite[Thm. 4.3(ii)-(iii)]{themelis2020douglas}.
					\item{\sc \cref{ass:DRS*}.}~
						In this case, the dual formulation \eqref{eq:D} satisfies \cref{ass:DRS}.
						In particular, adopting the \(*\)-notation of \cref{thm:selfdual} we have that \(\DRE*(s_*^k)=-\DRE(s^k)\) is decreasing and bounded below by \(\inf\psi=-\min\varphi\) (cf. \eqref{eq:strongDuality}), and \(u_*^k-v_*^k=\frac{v^k-u^k}{\gamma}\to0\).
						Note that inequality \eqref{eq:Qlb} guarantees that \(\seq{u^k}\) and \(\seq{v^k}\) are bounded; in what follows, we will make use of \cref{ass:DRS*:relint} to show that they actually converge to the global minimizer of \(\varphi\).
						To this end, it will suffice to show that \(\psi(v_*^k)\to\inf\psi\), since
						\[
							\DRE*(s_*^k)
							{}-{}
							\inf\psi
						{}\leq{}
							\psi(v_*^k)
							{}+{}
							\tfrac{1+\gamma_*L_{\sDRS*}}{2\gamma_*}
							\|u_*^k-v_*^k\|^2
							{}-{}
							\inf\psi
						{}={}
							\psi(v_*^k)
							{}+{}
							\tfrac{1+\mu_{\sDRS}\gamma}{2\mu_{\sDRS}}
							\|u_*^k-v_*^k\|^2
							{}-{}
							\inf\psi
						\]
						and \(\|u_*^k-v_*^k\|\to0\), where the inequality uses \(\bar u=v_*^k\) in \cref{thm:QG} applied to the dual formulation \eqref{eq:D}.
						We start by observing that \eqref{eq:prox:subdiff} implies that \(\frac{2u_*^k-s_*^k-v_*^k}{\gamma_*}\in\partial\nsDRS*(v_*^k)\) and \(\frac{s_*^k-u_*^k}{\gamma_*}=\nabla\sDRS*(u_*^k)\), which together result in
						\(
							\partial\psi(v_*^k)
						{}\ni{}
							\frac{u_*^k-v_*^k}{\gamma_*}
							{}-{}
							(\nabla\sDRS*(u_*^k)-\nabla\sDRS*(v_*^k))
						\),
						owing to the calculus rule of \cite[Thm. 23.8]{rockafellar1970convex}.
						Therefore,
						\(
							\dist(0,\partial\psi(v_*^k))
						{}\leq{}
							\|u_*^k-v_*^k\|
						{}\to{}
							0
						\).
						Adopting the terminology of \cite[Def. 4.2.1]{auslender2002asymptotic}, this means that \(\seq{v_*^k}\) is a \emph{stationary sequence} for \(\psi\).
						We will now show that \(\psi\) is \emph{asymptotically well behaved}, so that stationarity of \(\seq{v_*^k}\) guarantees that \(\psi(v_*^k)\to\inf\psi\), cf. \cite[Def. 4.2.2]{auslender2002asymptotic}.
						To this end, observe that \(\conj\psi=\conj{(\conj{(\refl\sDRS)}+\conj\nsDRS)}=(\refl\sDRS)\infconv g\), see \cite[Prop. 13.24(i) and Thm. 13.37]{bauschke2017convex}, where \(\infconv\) denotes the infimal convolution operator.
						It then follows from \cite[Prop. 12.6(ii)]{bauschke2017convex} that \(\dom\conj\psi=\dom(\refl\sDRS)+\dom\nsDRS=\dom\nsDRS-\dom\sDRS\).
						The validity of \cref{ass:DRS*:relint} then ensures through \cite[Thm. 3.2.1(a)-(d) and Cor. 4.2.1(a)]{auslender2002asymptotic} the sought asymptotic well behavior of \(\psi\), hence that \(\psi(v_*^k)\to\inf\psi=-\min\varphi\).
						This concludes the proof of assertion \ref{thm:DRS*:omega}.
						
						Finally, the equivalence of the condition \(0\in\interior(\dom\sDRS-\dom\nsDRS)\) and level boundedness of the dual cost \(\psi\) is shown in \cite[Thm. 5.2.1(a)]{auslender2002asymptotic}.
						Assertion \ref{thm:DRS:bounded} applied to the dual formulation then ensures boundedness of \(\seq{s_*^k}\), and that of \(\seq{s^k}\) then follows from the identity \(s_k=-\nicefrac{s_*^k}{\gamma_*}\).
					\qedhere
					\end{proofitemize}
				\end{proof}
			\end{thm}
			
			\begin{thm}[Subsequential convergence of \cref{alg:ADMM}]\label{thm:ADMM_subseq}
				Suppose that either \cref{ass:ADMM} or \cref{ass:ADMM*} is satisfied, and consider the (possibly infinite) iterates generated by \cref{alg:ADMM} with tolerance \(\varepsilon=0\).
				Then, the sequence of squared residuals \(\seq{\|r^k\|^2}\) has finite sum.
				Moreover,
				\begin{itemize}[label={}, leftmargin=0pt]
				\item{\sc under \cref{ass:ADMM}:}
					\begin{enumerate}
					\item\label{thm:ADMM:omega}
						All cluster points \((x,y,z)\) of \(\seq{x^k,y^k,z^k}\) satisfy the KKT conditions
						\begin{equation}\label{eq:ADMM_KKT}
							-\trans Ay\in\partial f(x),
						\quad
							-\trans By\in\partial g(z),
						\quad
							Ax+Bz=b,
						\end{equation}
						and attain the same (finite) cost \(f(x)+g(z)\), this being the limit of \(\seq{\LL_\beta(x^k,z^k,y^k)}\).
					\item\label{thm:ADMM:bounded}
						The sequence \(\seq{Ax^k,y^k,Bz^k}\) is bounded provided that the cost function \(\Phi\) is level bounded.
						If, additionally, \(f\in C^{1,1}(\R^m)\), then the sequence \(\seq{x^k,y^k,z^k}\) is bounded.
					\end{enumerate}
				\item{\sc under \cref{ass:ADMM*}:}
					\begin{enumerate}[resume]
					\item\label{thm:ADMM*:omega}%
						\(\seq{\LL_\beta(x^k,z^k,y^k)}\nearrow\min\Phi\) as \(k\to\infty\), \(\seq{Ax^k}\) and \(\seq{Bz^k}\) are convergent, and all the cluster points \((x,y,z)\) of \(\seq{x^k,y^k,z^k}\) are KKT-optimal triplets (satisfying \eqref{eq:ADMM_KKT}).
					\item\label{thm:ADMM*:bounded}%
						If \(b\in\interior(A\dom\sADMM+B\dom\nsADMM)\) (equivalently, if the dual cost \(\Psi(y)\coloneqq\conj f(-\trans Ay)+\conj g(-\trans By)+\innprod by\) is level bounded), then \(\seq{x^k,y^k,z^k}\) remains bounded.
					\end{enumerate}
				\end{itemize}
				\begin{proof}
					We shall again invoke the equivalence of \cref{alg:DRS,alg:ADMM} stated in \cref{thm:DRS=ADMM}, and import the same notation for convenience.
					Note that no clash occurs in using \(r^k\) to express the residual in both algorithms, having \(r^k=u^k-v^k=Ax^k+Bz^k-b\) through the relations of \cref{thm:DRS=ADMM}.
					In particular, that the squared residuals \(\seq{\|r^k\|^2}\) have finite sum is a direct consequence of \cref{thm:DRS_subseq}.
					We now consider the two separate cases.%
					\begin{proofitemize}
					\item{\sc\cref{ass:ADMM}.}~
						Having shown the vanishing of the residual, the proof is identical to that of \cite[Thm. 5.6]{themelis2020douglas}.
					\item{\sc\cref{ass:ADMM*}.}~
						Notice that \(\dom\sDRS-\dom\nsDRS=A\dom\sADMM+B\dom\nsADMM-b\) for \(\sDRS\) and \(\nsDRS\) as in \eqref{eq:CP2P}, owing to the fact that \(\dom\epicomp Af=A\dom f\) and \(\dom\epicomp Bg(b-{}\cdot{})=b-B\dom g\), see \cite[Prop. 12.36(i)]{bauschke2017convex}.
						Notice further that
						\[
							\conj\sDRS=f(-\trans A{}\cdot{})
						\quad\text{and}\quad
							\conj\nsDRS=g(-\trans B{}\cdot{})+\innprod{b}{{}\cdot{}},
						\]
						see \cite[Prop.s 13.23(iii) and 13.24(iv)]{bauschke2017convex}.
						Therefore, it suffices to show that any limit triplet \((x,y,z)\) is KKT-optimal, for all other claims are direct translations of \cref{thm:DRS*:omega}.
						To this end, suppose that a subsequence \(\seq{x^k,y^k,z^k}[k\in K]\) converges to \((x,y,z)\).
						Since \(r^k\to0\), necessarily \(Ax+Bz=b\).
						Moreover,
						\[
							\inf\Phi
						{}={}
							\lim_{K\ni k\to\infty}\LL_\beta(x^k,z^k,y^k)
						{}={}
							\lim_{K\ni k\to\infty}f(x^k)+g(z^k)
						{}\geq{}
							f(x)+g(z)
						{}={}
							\Phi(x,z)
						{}\geq{}
							\inf\Phi,
						\]
						where the second equality follows from the fact that \(\seq{y^k}[k\in K]\) is bounded (since it converges), the first inequality from lsc of \(f\) and \(g\), and the last equality from the fact that \(Ax+Bz-b=0\).
						In fact, notice that having \(f(x^k)+g(z^k)\to f(x)+g(z)\) implies that \(f(x^k)\to f(x)\) and \(g(z^k)\to g(z)\).
						Therefore, the continuity of the convex subdifferential \cite[Thm. 24.4]{rockafellar1970convex} together with \cref{thm:ADMM:-Ay,thm:ADMM:-By} implies the sought inclusions \(-\trans Ay\in\partial f(x)\) and \(-\trans By\in\partial g(z)\), yielding the claimed KKT-optimality.
						\qedhere
					\end{proofitemize}
				\end{proof}
			\end{thm}

		\subsection{Superlinear convergence}\label{sec:Superlinear}

			We now provide sufficient conditions that ensure superlinear convergence of the proposed algorithms.
			For the sake of simplicity we limit the analysis to \cref{alg:DRS} under \cref{ass:DRS}; the other cases can be inferred through the equivalence between \cref{alg:DRS} and \cref{alg:ADMM} and the self-duality of \ref{DRS} proven in \cref{thm:selfdual,thm:DRS=ADMM}.
			
			As a measure of ``quality'' of the oracle producing the update directions, borrowing the terminology of \cite[\S7.5]{facchinei2003finite}, we say that \(\seq{d^k}\) is a sequence of \DEF{superlinear directions} (relative to a sequence \(\seq{s^k}\) converging to \(s_\star\)) if
			\begin{equation}\label{eq:superd}
				\lim_{k\to\infty}{
					\frac{
						\|s^k+d^k-s_\star\|
					}{
						\|s^k-s_\star\|
					}
				}
			{}={}
				0.
			\end{equation}
			The next theorem shows that whenever the algorithm converges to a strong local minimum and the directions comply with the qualitative criterion \eqref{eq:superd}, the linesearch condition at \cref{step:DRS:LS} is eventually always passed at the first trial, and the iterates reduce to \(s^{k+1}=s^k+d^k\) and converge superlinearly.
			To do so, we will use the following lemma showing how strong local minimality on the original cost reflects on the DRE.
			
			\begin{lem}\label{thm:DRS:strmin}%
				Suppose that \cref{ass:DRS} holds, and consider the iterates generated by \cref{alg:DRS} with \(\varepsilon=0\), and let \(\varphi_\star\) be the limit of the sequence \(\seq{\DRE(s^k)}\).
				Suppose that \(\seq{u^k}\) converges to a strong local minimum \(u_\star\) of \(\varphi\).
				Then, there exists \(\delta>0\) such that \(\DRE(s^k)-\varphi_\star\geq\tfrac\delta2\|s^k-s_\star\|^2\) holds for every \(k\) large enough, where
				\(
					s_\star
				{}\coloneqq{}
					u_\star+\gamma\nabla\sDRS(u_\star)
				\).
				\begin{proof}
					We begin by observing that for every \(k\in\N\) and \(\epsilon>0\) it holds that
					\begin{equation}\label{eq:Young}
						\tfrac{1}{(1+\gamma L_{\sDRS})^2}
						\|s^k-s_\star\|^2
					{}\leq{}
						\|u^k-u_\star\|^2
					{}\leq{}
						\tfrac{\epsilon+1}{\epsilon}\|u^k-v^k\|^2
						{}+{}
						(1+\epsilon)\|v^k-u_\star\|^2,
					\end{equation}
					where the first inequality follows from \(\frac{1}{1+\gamma L_{\sDRS}}\)-strong monotonicity of \(\prox_{\gamma\sDRS}\) \cite[Prop. 2.3(ii)]{themelis2020douglas} and the fact that \(u_\star=\prox_{\gamma\sDRS}(s_\star)\) as it follows from \eqref{eq:prox:subdiff}, and the second one uses Young's inequality.
					As ensured by \cref{thm:DRS:omega}, \(u^k-v^k\to0\) and \(\varphi(u_\star)=\varphi_\star\).
					Therefore, \(v^k\to u_\star\) and because of strong local minimality there exists \(\mu>0\) and \(K\in\N\) such that \(\varphi(v^k)-\varphi_\star\geq\tfrac\mu2\|v^k-u_\star\|^2\) for all \(k\geq K\).
					For all \(\epsilon>0\) and \(k\geq K\) we thus have
					\begin{align*}
						\DRE(s^k)
					{}\overrel[\geq]{\ref{thm:sandwich}}{} &
						\varphi(v^k)
						{}+{}
						\tfrac{1-\gamma L_{\sDRS}}{2\gamma}
						\|v^k-u^k\|^2
					\\
					{}\geq{} &
						\varphi_\star
						{}+{}
						\tfrac\mu2
						\|v^k-u_\star\|^2
						{}+{}
						\tfrac{1-\gamma L_{\sDRS}}{2\gamma}
						\|v^k-u^k\|^2
					\\
					{}\overrel[\geq]{\eqref{eq:Young}}{} &
						\varphi_\star
						{}+{}
						\tfrac{\mu}{2(1+\epsilon)(1+\gamma L_{\sDRS})^2}
						\|s^k-s_\star\|^2
						{}+{}
						\left(
							\tfrac{1-\gamma L_{\sDRS}}{2\gamma}
							{}-{}
							\tfrac{\mu}{2\epsilon}
						\right)
						\|v^k-u^k\|^2
					\\
					{}={} &
						\varphi_\star
						{}+{}
						\tfrac{1-\gamma L_{\sDRS}}{(1-\gamma(L_{\sDRS}-\mu))(1+\gamma L_{\sDRS})^2}
						\tfrac\mu2
						\|s^k-s_\star\|^2,
					\end{align*}
					where the last equality uses
					\(
						\epsilon
					{}={}
						\tfrac{\gamma\mu}{1-\gamma L_{\sDRS}}
					\).
				\end{proof}
			\end{lem}
			
			\begin{thm}[Acceptance of the unit stepsize]\label{thm:tau1}%
				Suppose that \cref{ass:DRS} holds, and consider the iterates generated by \cref{alg:DRS}.
				Suppose that \(\seq{u^k}\) converges to a strong local minimum \(u_\star\) of \(\varphi\) and that \(\seq{d^k}\) are superlinear directions as in \eqref{eq:superd}.
				Then, eventually unit stepsize \(\tau_k=1\) is always accepted, hence the iterates reduce to \(s^{k+1}=s^k+d^k\) and converge superlinearly.
				\begin{proof}
					In light of \cref{thm:DRS:strmin}, by possibly discarding the first iterates we may assume that
					\[
						\DRE(s^k)-\DRE(s_\star)
					{}\geq{}
						\tfrac\delta2\|s^k-s_\star\|^2
					\quad
						\text{for all \(k\)'s}
					\]
					for some \(\delta>0\), where \(s_\star\coloneqq u_\star+\gamma\nabla\sDRS(u_\star)\) is the limit point of \(\seq{s^k}\).
					Combined with \cref{thm:QG}, we obtain that
					\[
						\epsilon_k
					{}\coloneqq{}
						\frac{
							\DRE(s^k+d^k)-\DRE(s_\star)
						}{
							\DRE(s^k)-\DRE(s_\star)
						}
					{}\leq{}
						\frac{1+\gamma L_f}{\gamma\delta}
						\frac{
							\|s^k+d^k-s_\star\|^2
						}{
							\|s^k-s_\star\|^2
						}
					{}\to{}
						0.
					\]
					Since \(u^k-v^k\to0\), we have that \(\bar s^k\to s_\star\).
					Therefore, eventually \(\DRE(\bar s^k)\geq\DRE(s_\star)\) and \(\epsilon_k\leq1\).
					Then, denoting \(\C[]=\C\) as in \eqref{eq:DRE:C} we have
					\begin{align*}
						\DRE(s^k+d^k)-\DRE(s^k)
					{}={} &
						-(1-\epsilon_k)
						\bigl(
							\DRE(s^k)-\DRE(s_\star)
						\bigr)
					\\
					{}\leq{} &
						-(1-\epsilon_k)
						\bigl(
							\DRE(s^k)-\DRE(\bar s^k)
						\bigr)
					\\
					{}\overrel[\leq]{\eqref{eq:DRE:C}}{} &
						-(1-\epsilon_k)
						\tfrac{\C[]}{\gamma}
						\|r^k\|^2.
					\end{align*}
					Since \(\epsilon_k\to0\) and \(c<\C[]\) as required in \cref{alg:DRS}, eventually
					\(
						1-\epsilon_k
					{}\geq{}
						\frac{c}{\C[]}
					\),
					resulting in
					\(
						\DRE(s^k+d^k)
					{}\leq{}
						\DRE(s^k)
						{}-{}
						\tfrac c\gamma\|r^k\|^2
					\),
					proving that \(\tau_k=1\) passes the condition at \cref{step:DRS:LS}.
				\end{proof}
			\end{thm}
			
			Although the \ref{DRS}-update \(s\mapsto s^+\) is not uniquely determined owing to the multi-valuedness of \(\prox_{\gamma\nsDRS}\), under some regularity assumptions not only is it single valued, but even differentiable, when close to solutions.
			To see this, observe that
			\[
				\Res
			{}={}
				\Res[\gamma][fb]\circ\prox_{\gamma\sDRS},
			\]
			where
			\(
				\Res
			{}\coloneqq{}
				\prox_{\gamma\sDRS}-\prox_{\gamma\nsDRS}(2\prox_{\gamma\sDRS}-\id)
			\)
			is the Douglas-Rachford residual and
			\[
				\Res[\gamma][fb](u)\coloneqq u-\FB u
			\]
			is the residual of \emph{forward-backward splitting} (FBS), see \eg \cite{attouch2013convergence,bolte2014proximal,themelis2018forward}.
			By combining \cite[Prop. 13.24 and Ex. 13.35]{rockafellar2011variational} and \cite[Thm. 4.4(c)-(f) and Cor. 4.7]{poliquin1996generalized}, it follows that \(\prox_{\gamma\sDRS}(s)\) is continuously differentiable around \(s_\star\) provided that \(\sDRS\) is twice continuously differentiable around \(s_\star\).
			Thus, under this assumption, from the chain rule of differentiation we conclude that \(\Res\) is (strictly) differentiable at \(s_\star\) provided that \(\Res[\gamma][fb]\) is (strictly) differentiable at \(u_\star=\prox_{\gamma\sDRS}(s_\star)\).
			Sufficient conditions for this latter property to hold are documented in \cite[Thm. 4.10]{themelis2018forward}, namely twice (Lipschitz-) continuous differentiability of \(\sDRS\) around \(u_\star\), and prox-regularity and (strict) twice epi-differentiability of \(\nsDRS\) at \(u_\star\) for \(-\nabla\sDRS(u_\star)\) \cite[\S13.B and 13.F]{rockafellar2011variational}.
			
			\begin{thm}[Dennis-Moré criterion for superlinear directions]\label{thm:DM}%
				Suppose that \cref{ass:DRS} holds, and consider the iterates generated by \cref{alg:DRS}.
				Suppose that \(\seq{s^k}\) converges to a point \(s_\star\) at which \(\Res\) is strictly differentiable and with nonsingular Jacobian \(J\Res(s_\star)\).
				If the Dennis-Moré condition
				\begin{equation}\label{eq:DM}
					\lim_{k\to\infty}{
						\frac{
							\|\Res(s^k)+J\Res(s_\star)d^k\|
						}{
							\|d^k\|
						}
					}
				{}={}
					0
				\end{equation}
				holds, then \(\seq{d^k}\) are superlinear directions and the claims of \cref{thm:tau1} hold.
				\begin{proof}
					Since \(J\Res\) exists and is nonsingular at \(s_\star\), for \(k\) large enough \(\Res(s^k)\) is single valued and satisfies
					\(
						\|\Res(s^k)\|\geq\alpha\|s^k-s_\star\|
					\)
					for some \(\alpha>0\).
			% 		It follows from \cref{thm:DRS:strmin} that there exists \(\delta>0\) such that
			% 		\[
			% 			\DRE(s^k)-\varphi_\star
			% 		{}\geq{}
			% 			\tfrac\delta2\|s^k-s_\star\|^2
			% 		\]
			% 		for all \(k\)'s large enough.
			% 		Therefore, from the linesearch condition at \cref{step:DRS:LS} we have
			% 		\begin{equation}\label{eq:Rcoco}
			% 			\color{red}
			% 			c\|\Res(s^k)\|^2
			% 		{}\geq{}
			% 			\DRE(s^k)-\varphi(s^{k+1})
			% 		{}\geq{}
			% 			\DRE(s^k)-\varphi_\star
			% 		{}\geq{}
			% 			\tfrac\delta2\|s^k-s_\star\|^2.
			% 		\end{equation}
					Due to strict differentiability,
					\[
						\lim_{k\to\infty}{
							\frac{
								\|\Res(s^k+d^k)-\Res(s^k)-J\Res(s_\star)d^k\|
							}{
								\|d^k\|
							}
						}
					{}={}
						0,
					\]
					and from the Dennis-Moré condition \eqref{eq:DM} it then follows that
					\(
						\frac{
							\|\Res(s^k+d^k)\|
						}{
							\|d^k\|
						}
					{}\to{}
						0
					\).
					Since \(\|\Res(s^k+d^k)\|\geq\alpha\|s^k+d^k-s_\star\|\), we conclude that
					\(
						\frac{
							\|s^k+d^k-s_\star\|
						}{
							\|d^k\|
						}
					{}\to{}
						0
					\)
					too.
					Therefore,
					\[
						\frac{
							\|s^k+d^k-s_\star\|
						}{
							\|s^k-s_\star\|
						}
					{}\leq{}
						\frac{
							\|s^k+d^k-s_\star\|
						}{
							\|d^k\|
							\bigl|
								1-\frac{\|s^k+d^k-s_\star\|}{\|d^k\|}
							\bigr|
						}
					{}\to{}
						0,
					\]
					proving that \(\seq{d^k}\) are superlinear directions with respect to \(\seq{s^k}\).
				\end{proof}
			\end{thm}
			% 
			% \begin{cor}[Dennis-Moré condition]%
			% 	Consider the iterates generated by \cref{alg:ADMM}.
			% 	Suppose that \(\seq{(x^k,z^k,y^k)}\) converges to a point satisfying the strong local minimality condition \eqref{eq:ADMMstrmin}, and that \(\Res[\gamma][ADMM]\) is strictly differentiable at \(y_\star\).
			% 	If the Dennis-Moré condition
			% 	\begin{equation}
			% 		\lim_{k\to\infty}{
			% 			\frac{
			% 				\|\Res[\gamma][ADMM](y^k)+J\Res[\gamma][ADMM](y_\star)d^k\|
			% 			}{
			% 				\|d^k\|
			% 			}
			% 		}
			% 	{}={}
			% 		0
			% 	\end{equation}
			% 	holds, then \(\seq{d^k}\) are superlinear directions.
			% \end{cor}

			We conclude by showing that the modified Broyden scheme described in \cref{sec:Broyden} enables superlinear rates when some regularity requirements are met at the limit point.
			These include \emph{Lipschitz semidifferentiability} of the residual \(\Res\), a condition that entails classical differentiability at the limit point but not necessarily around it, see \cite{ip1992local}.
			
			\begin{thm}[Superlinear convergence with Broyden directions]\label{thm:Broyden}%
				Suppose that \cref{ass:DRS} holds, and consider the iterates generated by \cref{alg:DRS} with directions being selected with the modified Broyden method (\cref{sec:Broyden}).
				Suppose that the sequence of Broyden matrices \(\seq{H_k}\) is bounded, and that \(\seq{s^k}\) converges to a strong local minimum \(s_\star\) of \(\varphi\) at which \(\Res\) is Lipschitz-semidifferentiable and has a nonsingular Jacobian \(J\Res(s_\star)\).
				Then, the Dennis-Moré condition \eqref{eq:DM} holds, and in particular the unit stepsize \(\tau_k=1\) is eventually always accepted and \(\seq{s^k}\) converges superlinearly.
				\begin{proof}
					Let \(p_k\) and \(q_k\) be as in \eqref{eq:DRS:pq}, \(H_k\) be as in \eqref{subeq:Broyden}, and denote \(G_\star\coloneqq J\Res(s_\star)\).
					Since \(\Res\) is differentiable at \(s_\star\) with nonsingular Jacobian and since \(s^k\to s_\star\), for \(k\) large enough \(\Res(s^k)\) is single valued and satisfies
					\(
						\|\Res(s^k)\|\geq\alpha\|s^k-s_\star\|
					\)
					for some \(\alpha>0\).
					Moreover, it follows from \cite[Lem. 2.2]{ip1992local} that an \(L>0\) exists such that
					\begin{align*}
						\frac{
							\|q_k-G_\star p_k\|
						}{
							\|p_k\|
						}
					{}={} &
						\frac{
							\|\Res(s^k+d^k)-\Res(s^k)-G_\star d^k\|
						}{
							\|d^k\|
						}
					{}\leq{}
						L
						\max{
							\set{
								\|s^k+d^k-s_\star\|,
								\|s^k-s_\star\|
							}
						}
					\\
					\numberthis\label{eq:frac_pq}
					{}\leq{} &
						L
						\bigl(
							\|s^k-s_\star\|
							{}+{}
							\|d^k\|
						\bigr)
					\end{align*}
					for \(k\) large enough.
					Since \(\prox_{\gamma\sDRS}\) is \(\frac{1}{1-\gamma L_{\sDRS}}\)-Lipschitz continuous \cite[Prop. 2.3(ii)]{themelis2020douglas}, denoting \(u_\star\coloneqq\prox_{\gamma\sDRS}(s_\star)\) we have that
					\[
						\|\Res(s^k)\|^2
					{}\geq{}
						\alpha^2\|s^k-s_\star\|^2
					{}\geq{}
						\alpha^2(1-\gamma L_{\sDRS})^2\|u^k-u_\star\|^2
					{}\overrel*[\geq]{\ref{thm:QG}}{}
						\alpha'(\DRE(s^k)-\varphi_\star)
					\]
					for some \(\alpha'>0\).
					Therefore,
					\[
						\DRE(s^{k+1})-\varphi_\star
					{}\leq{}
						\DRE(s^k)-\varphi_\star
						{}-{}
						c\|\Res(s^k)\|^2
					{}\leq{}
						(1-c\alpha')(\DRE(s^k)-\varphi_\star),
					\]
					and invoking \cref{thm:DRS:strmin} we conclude that \(\seq{\|s^k-s_\star\|}\) converges \(R\)-linearly and thus has finite sum.
					Since \(\Res\) is differentiable at \(s_\star\), there exists \(\alpha''>0\) such that \(\|\Res(s^k)\|=\|\Res(s^k)-\Res(s_\star)\|\leq\alpha''\|s^k-s_\star\|\).
					Therefore, also \(\seq{\|\Res(s^k)\|}\) has finite sum, and in turn so does \(\seq{\|d^k\|}\) owing to boundedness of \(\seq{H_k}\) and the fact that \(d^k=-H_k\Res(s^k)\).
					From \eqref{eq:frac_pq} we conclude that
					\(
						\seq{
							\frac{
								\|q_k-G_\star p_k\|
							}{
								\|p_k\|
							}
						}
					\)
					has finite sum as well, and the claimed Dennis-Moré condition follows by verbatim importing the conclusions of the proof of \cite[Thm. VI.8]{themelis2019supermann}, after observing that
					\(
						\frac{
							\|\Res(s^k)+G_\star d^k\|
						}{
							\|d^k\|
						}
					{}={}
						\frac{
							\|(H_k^{-1}-G_\star)p_k\|
						}{
							\|p_k\|
						}
					\).
					Finally, the acceptance of the unit stepsize and superlinear convergence of \(\seq{s^k}\) follow from \cref{thm:tau1,thm:DM}.
				\end{proof}
			\end{thm}%

	\section{Simulations}\label{sec:Simulations}

		In this section we show the effectiveness of the proposed \cref{alg:DRS,alg:ADMM}, for different choices of the linesearch direction, compared to the standard \ref{DRS} and \ref{ADMM}.
		In all problems, the (generalized) quadratic structure of one component of the cost functions is exploited as described in \cref{rem:linearprox}, thus resulting in at most two evaluations of \(\prox_{\gamma\sDRS}\) for every iteration of \cref{alg:DRS}, or \(x\)-minimizations for every iteration of \cref{alg:ADMM}.
		The implementations in Julia of all the algorithms are available online as part of the
		\emph{ProximalAlgorithms.jl} package.\footnote{\url{https://github.com/JuliaFirstOrder/ProximalAlgorithms.jl}}
		All experiments were run using Julia 1.6.3.

		\subsection{Nonconvex sparse least squares}\label{sec:SLS}

			To find a sparse, least-squares solution $x \in \R^n$ to a linear system $Ax = b$, we consider
			the formulation
			\begin{equation}\label{eq:sparse_least_squares}
				\minimize_{x\in\R^n} \tfrac{1}{2}\|Ax - b\|^2 + r\|x\|_{\nicefrac12}^{\nicefrac12},
			\end{equation}
			where $A\in\R^{m\times n}$, $\|x\|_{\nicefrac12}^{\nicefrac12} = \sum_i \sqrt{|x_i|}$ is the square root of the $\ell_{\nicefrac12}$ quasi-norm,
			and $r > 0$ is a regularization parameter.
			The $\ell_{\nicefrac12}$ penalty term has favorable properties compared to the popular $\ell_1$ regularization, as thoroughly documented in \cite[Sec. II]{xu2012regularization}; yet its nonconvexity makes problem \eqref{eq:sparse_least_squares} more challenging to solve.
			As derived in \cite{xu2012regularization}, the proximal operator for the regularization term can be computed in closed form as follows
			\begin{equation}\label{eq:l_half_prox}
				\left[\prox_{\gamma\|\cdot\|_{\nicefrac12}^{\nicefrac12}}(x)\right]_i
			{}={}
				\begin{ifcases}
					\frac12\left(
						1
						{}+{}
						\cos\frac23\left(
							\pi
							{}-{}
							\arccos\frac{\gamma}{4}\left(
								\frac{x_i}{3}
							\right)^{-\nicefrac32}
						\right)
					\right)x_i
				&
					|x_i| > \frac32\gamma^{\nicefrac23}
				\\
					0\otherwise.
				\end{ifcases}
			\end{equation}
			For the least squares term, \(\prox_{\nicefrac{\gamma}{2}\|A\cdot - b\|^2}(x) = (A^T A + \gamma^{-1} \I)^{-1} (A^T b + \gamma^{-1}x)\) involves solving a positive definite linear system, and is therefore the most computationally demanding operation.
			% {\color{blue}Note that since it is an affine mapping, only (at most) two evaluations of \(\prox_{\nicefrac{\gamma}{2}\|A\cdot - b\|^2}\) will be required in \cref{alg:DRS}, cf. \cref{rem:linearprox}.}
			
			We generated random instances of problem \eqref{eq:sparse_least_squares} similarly to the setup of \cite[Sec. 8.2]{daubechies2010iteratively}:
			matrix $A$ has i.i.d. Gaussian entries with variance $\nicefrac{1}{m}$, while $b = A\hat{x}$ for a random vector $\hat{x}$ with $k$ nonzero coefficients. In the \ref{DRS}-\cref{alg:DRS}, we used
			\(\lambda=1\),
			\(\gamma=0.95 \cdot L_{\sDRS}^{-1}\),
			\(c=\tfrac{1}{2}C(\gamma L_{\sDRS}, \lambda)\),
			cf. \eqref{eq:DRE:C} and \eqref{eq:DRE:gamma}.
			
			\begin{figure}[p]
				\caption{%
					{\rm \S\ref{sec:SLS}: nonconvex sparse least squares problem \eqref{eq:sparse_least_squares}.}
					Comparison between \ref{DRS} and the linesearch variant \cref{alg:DRS} using modified Broyden, L-BFGS, and Nesterov acceleration directions.
				}
				\begin{subfigure}{\linewidth}%
				\includetikz[width=0.92\linewidth]{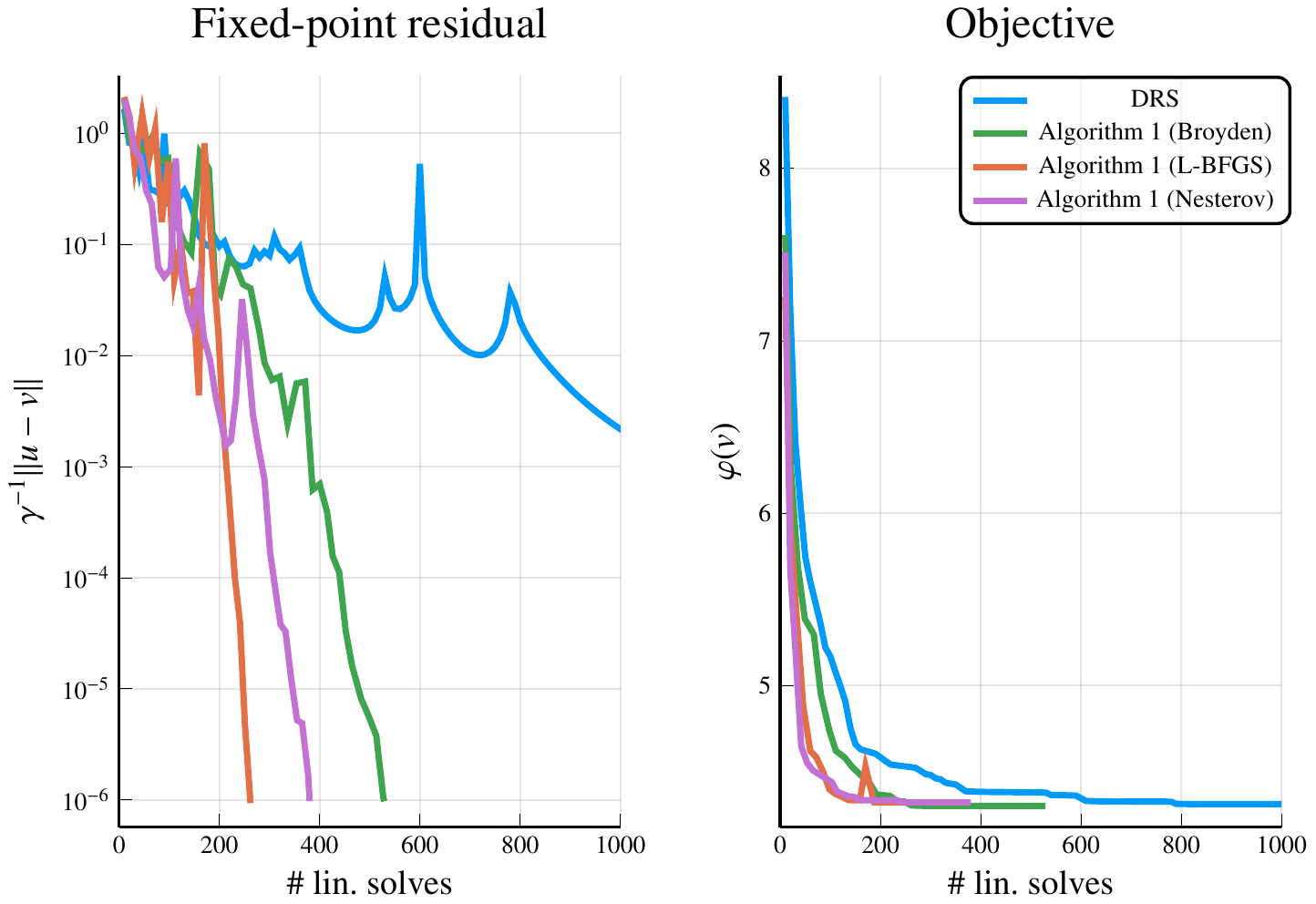}%
				\caption{%
					Comparison of the methods on one randomly generated problem instance: matrix $A\in\R^{100\times 500}$ has random coefficients with variance $0.01$, while $b$ was generated from a known sparse $\hat{x}$ with $50$ nonzero coefficients, and $r=0.1$.
					On the \(x\)-axis, the number of linear systems solved (needed for the \(u\)-update): in the case of \ref{DRS}, this coincides with the number of iterations, while for \cref{alg:DRS} it accounts for all operations performed in the linesearch.
				}%
				\label{fig:sparse_linreg_single}%
				\end{subfigure}
				\begin{subfigure}{\linewidth}%
					\includetikz[width=0.92\linewidth]{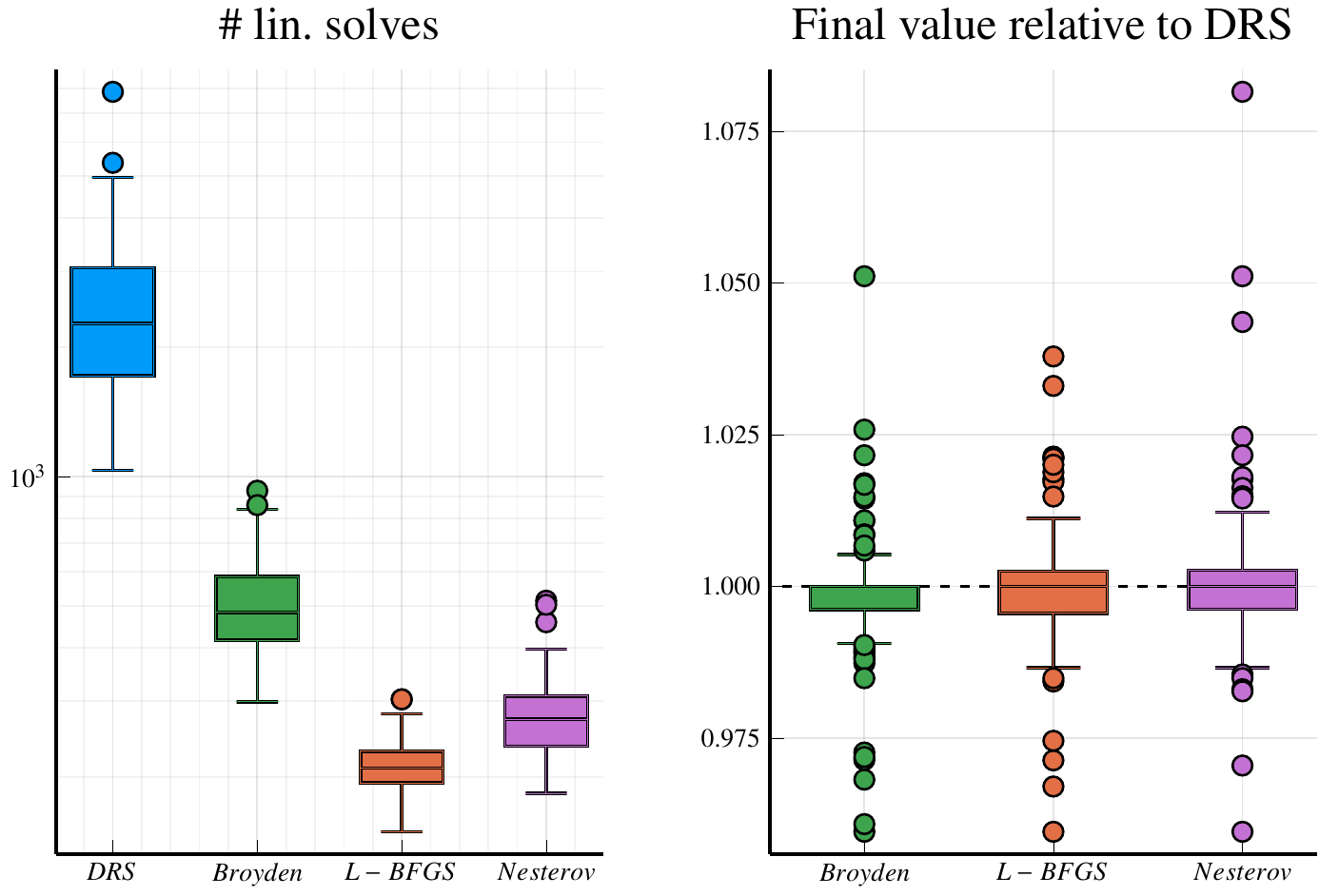}%
					\caption{%
						Comparison of the methods on $100$ randomly generated problems.
						Left: box plot of the number of operations to reach a fixed-point residual of norm $10^{-6}$.
						Right: box plot of the objective value achieved, relative to the one obtained by \ref{DRS}.
						The boxes represent the interquartile range (the range between the \(\text{P25}\) and \(\text{P75}\) percentiles),
						with the median highlighted; the bars extend below to \(\text{P25} - 1.5(\text{P75} - \text{P25})\), and above to \(\text{P75} + 1.5(\text{P75} - \text{P25})\); all observations outside of this range are displayed as stand-alone points.%
					}%
					\label{fig:sparse_linreg_multi}%
				\end{subfigure}
			\end{figure}
			
			The convergence of the proposed \cref{alg:DRS} compared to \ref{DRS} is exemplified in \cref{fig:sparse_linreg_single}, where the algorithms were applied on a randomly generated problem instance with $n = 500$, $m = 100$, and $k = 50$.
			For this experiment, in \cref{alg:DRS} we used Broyden, L-BFGS (with memory $5$) and Nesterov directions, cf. \cref{sec:Directions}: these choices of directions significantly accelerate convergence, compared to \ref{DRS}.
			Using Anderson acceleration directions, in this case, did not perform well.
			In \cref{fig:sparse_linreg_multi}, the result of running the algorithms on $100$ randomly generated problems is displayed: there, the distribution of the number of linear solves (i.e. evaluations of \(\prox_{\nicefrac{\gamma}{2}\|A\cdot - b\|^2}\)) over all problems is shown for each algorithm, as well as the distribution of the best objective values reached relative to the one obtained by \ref{DRS}.
			It is apparent how \cref{alg:DRS} converges to critical points in a fraction of the operations required by \ref{DRS}.

		\subsection{Sparse PCA}\label{sec:SPCA}

			Given a dataset of points in \(\R^n\), the goal of \emph{sparse principal component analysis (SPCA)} is to explain as much variability in the data as possible by using only \(k\ll n\) variables.
			Let the data matrix be \(W\in\R^{m\times n}\) (this can be assumed to be centered, \ie, with zero-mean columns), then the problem can be formulated as follows:
			\begin{equation}\label{eq:sparse_PCA}
				\maximize_{x\in\R^n}\tfrac{1}{2m}\|Wx\|^2
			\quad\stt{}
				\|x\|=1,~
				\|x\|_0\leq k,
			\end{equation}
			where the \(\ell_0\)-quasi-norm \(\|x\|_0\) denotes the number of nonzero elements of vector \(x\). Being \(\tfrac{1}{m}\trans WW \in\R^{n\times n}\) the covariance matrix of \(W\), \eqref{eq:sparse_PCA} amounts to a variance maximization problem.
			The formulation \eqref{eq:sparse_PCA} was first introduced in \cite{daspremont2005direct}, where the authors propose solving an SDP relaxation.
			Here, we consider tackling directly the nonconvex problem \eqref{eq:sparse_PCA} instead.
			
			Denoting the set of feasible points by
			\[
				\mathcal S
			{}\coloneqq{}
				\set{x\in\R^n}[\|x\|=1,~\|x\|_0\leq k],
			\]
			problem \eqref{eq:sparse_PCA} takes the form \eqref{eq:P} once we set \(\sDRS(x) = -\tfrac{1}{m}\|Wx\|^2\) and \(\nsDRS(x) = \indicator_{\mathcal S}(x)\) (the \DEF{indicator function} of \(\mathcal S\), namely \(\indicator_{\mathcal S}(x)=0\) if \(x\in\mathcal S\) and \(\indicator_{\mathcal S}(x)=\infty\) otherwise).
			So formulated, the problem complies with \cref{ass:DRS}, therefore \ref{DRS} can be readily applied.
			Note that, in this case,
			\begin{align*}
				\prox_{\gamma\sDRS}(x)
			{}={}
				\argmin_z\set{
					-\tfrac{1}{2m}\|Wz\|^2+\tfrac{1}{2\gamma}\|z-x\|^2
				}
			{}={} &
				\left(\I - \tfrac\gamma m\trans WW\right)^{-1}\!\!x
			\\
			{}={} &
				x - \trans W\left(\tfrac m\gamma\I - W\trans W\right)^{-1}\!Wx,
			\end{align*}
			where the last equality uses the Woodbury identity.
			Whenever \(\gamma < L_{\sDRS}^{-1} = \nicefrac{\|W\|^2}{m}\), evaluating \(\prox_{\gamma\sDRS}\) amounts to solving a square positive definite linear system of dimension either \(n\) or \(m\), depending on which one is smaller: this can be done by computing the Cholesky factor offline once, and caching the factorization throughout the iterations.
			Recall that, since \(\sDRS\) is quadratic, as illustrated in \cref{rem:linearprox} no more than two evaluations of its proximal mapping will be necessary at every iteration.
			In addition, evaluating \(\prox_{\gamma \nsDRS} = \proj_{\mathcal S}\) (the set-valued projection onto \(\mathcal S\)) amounts to setting to zero the \(n-k\) smallest coefficients of \(x\) (in magnitude), and normalizing the resulting vector to project it on the \(\ell_2\)-sphere:
			this has therefore a negligible cost.
			
\begin{figure}
				\includetikz[width=0.95\linewidth]{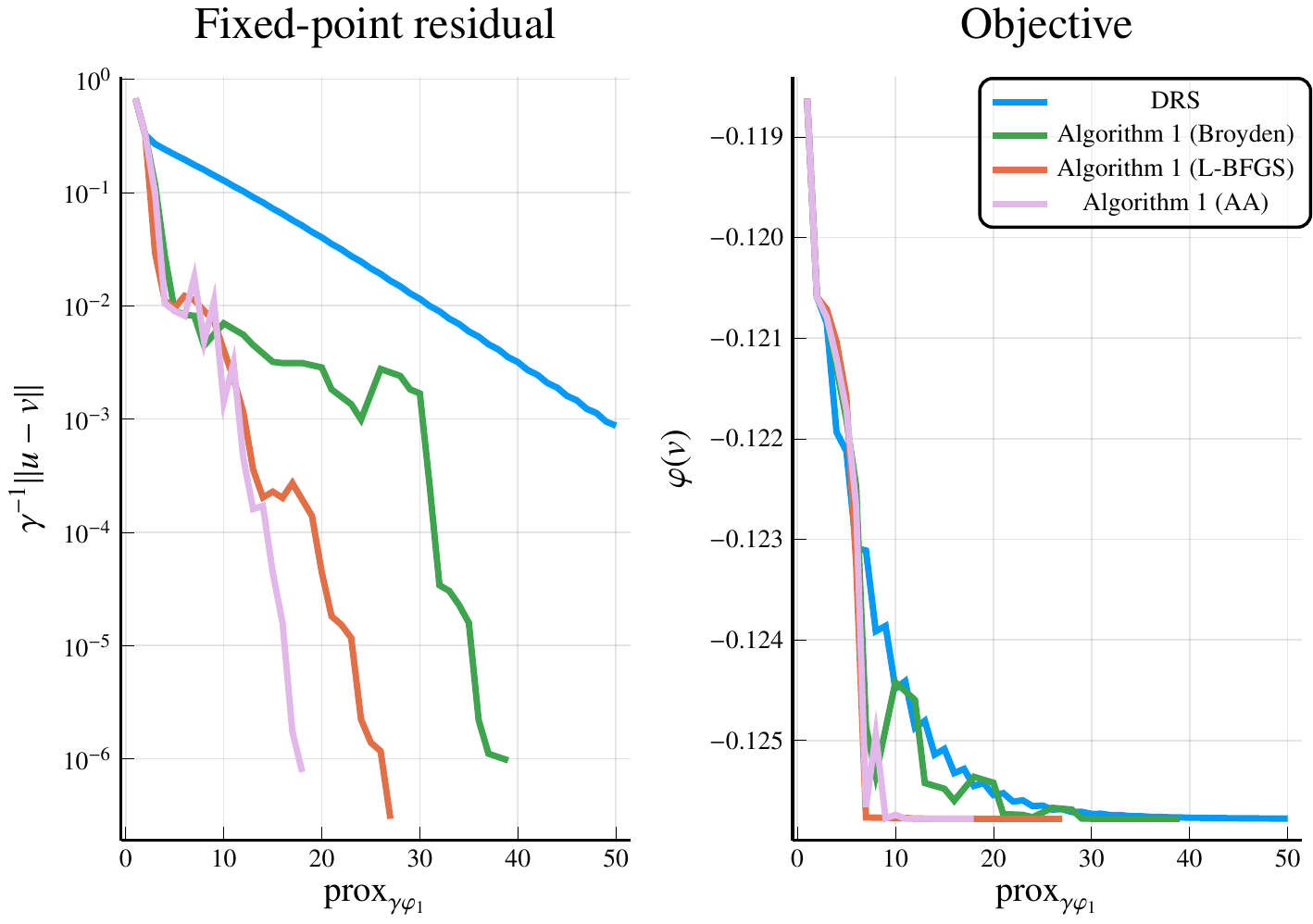}%
				\caption{%
					{\rm\S\ref{sec:SPCA}: sparse PCA problem \eqref{eq:sparse_PCA} on a subset of the \emph{20newsgroup} dataset (100 features only).}
					Comparison between \ref{DRS} and the linesearch variant \cref{alg:DRS}, using modified Broyden, L-BFGS, and Anderson acceleration directions,
					when applied to the sparse PCA problem \eqref{eq:sparse_PCA}.
					On the \(x\)-axis, the number of linear systems solved (needed for the \(u\)-update): in the case of \ref{DRS}, this coincides with the number of iterations, while for \cref{alg:DRS} it accounts for all operations performed in the linesearch.
					We used a memory parameter of 5 for L-BFGS and Anderson acceleration.%
				}%
				\label{fig:SPCA_single}%
			\end{figure}

			\Cref{fig:SPCA_single} shows the results when the algorithm is applied using a small subset of the \emph{20newsgroup} dataset,\footnote{%
				\url{https://cs.nyu.edu/~roweis/data.html}
			}
			which only retains 100 features from the original dataset: \(m = 16242\) and \(n = 100\).
			This is the same dataset that was used in the experiments in \cite{daspremont2005direct}.
			Here, the initial point was chosen as the vector \((1/n, \ldots, 1/n)\in\R^n\), which in our experiments gave consistently better results compared to random initialization, in terms of the objective value reached.
			In the \ref{DRS}-\cref{alg:DRS}, we used
			\(\lambda=1\),
			\(\gamma=\tfrac{0.95}{2} L_{\sDRS}^{-1}\),
			\(c=\tfrac{1}{2}C(\gamma L_{\sDRS}, \lambda)\),
			cf. \eqref{eq:DRE:C}-\eqref{eq:DRE:gamma}, and the directions given by the modified Broyden, L-BFGS and Anderson acceleration (the latter two with memory 5).
			From this experiment, these directions in \cref{alg:DRS} greatly improve the convergence over \ref{DRS}, both in terms of the fixed point residual norm as well as the objective value.
			We did not include results with Nesterov directions, as these did not perform well on this problem, a phenomenon that we believe is related to the nonconvexity of \(\sDRS\) (cf. \cref{rem:directions}).
			
		\subsection{Sparse PCA: consensus formulation}\label{sec:SPCAconsensus}

			As the problem size grows, a big limitation is the need to store and operate with large matrices.
			To account for this issue, we consider the following \emph{consensus formulation}: having fixed a number of \emph{agents} \(N\geq 1\), decompose matrix \(W\) into \(N\) row blocks \(W_1,\ldots,W_N\), so that \(\trans W=[\trans{W_1}~\cdots~\trans{W_N}\,]\) and \(\|Wx\|^2=\sum_{i=1}^N\|W_ix\|^2\), introduce \(N\) copies \(x_1,\ldots,x_N\) of \(x\) (stacked in a vector \(\bm x\in\R^{nN}\)), and solve
			\[
				\minimize_{x\in\R^{Nn},\bm z\in\R^{n}}{
					\sum_{i=1}^N{
						-\tfrac{1}{2m}\|W_ix_i\|^2
					}
				}
			\quad\stt{}
				\|z\|=1,~
				\|z\|_0\leq k,~
				x_i=z,~
				i=1\ldots N.
			\]
			This problem is equivalent to \eqref{eq:sparse_PCA}, and can be expressed in \ref{ADMM} form \eqref{eq:CP} as
			\begin{equation}\label{eq:cSPCA}
				\minimize_{
					\bm x\in\R^{Nn}\!,\,
					z\in\R^n
				}{
					\underbrace{\sum_{i=1}^N-\tfrac{1}{2m}\|W_ix_i\|^2}_{\sADMM(\bm x)}
					{}+{}
					\underbrace{\vphantom{\sum_{i=1}^N}\indicator_{\mathcal S}(z)}_{\nsADMM(z)}
				}
			\quad\stt{}~
				\bm x
				{}-{}
				\setlength\matrixcolsep{2pt}
				\begin{pmatrix}
					& \I &\\
					&\vdots&\\
					&\I&
				\end{pmatrix}
				z
			{}={}
				0.
			\end{equation}
			Apparently, the \ref{ADMM} matrix \(A\) is the \(nN\times nN\) identity, \(B\in\R^{nN\times n}\) is the vertical stacking of \(N\) many \(n\times n\) negative identity matrices, and \(b\) is the zero \(\R^{nN}\) vector.
			Notice that \cref{ass:ADMM} is satisfied, as \(\epicomp A\sADMM=\sADMM\) has Lipschitz-continuous gradient with modulus
			\(
				L_{\epicomp A\sADMM}
			{}={}
				\tfrac{1}{m}\max_{i=1\ldots N}\|W_i\|^2
			{}\leq{}
				\tfrac{1}{m}\|W\|^2
			\)
			and \(A\) has clearly full row rank.
			
			The \(z\)-update as prescribed by \ref{ADMM} comes at negligible cost, since
			\[
				\argmin_{z\in\R^n}\set{\indicator_{\mathcal S}(z)+\tfrac\beta2\|\bm x+Bz\|^2}
			{}={}
				\proj_{\mathcal S}
				\bigl(\textstyle
					\tfrac1N\sum_{i=1}^Nx_i
				\bigr)
			\quad
				\forall\bm x\in\R^{nN}.
			\]
			The \(\bm x\)-update amounts to solving (in parallel) \(N\) (small) linear systems:
			\[
				\argmin_{x_i\in\R^n}\set{-\tfrac{1}{2m}\|W_ix_i\|^2 + \tfrac\beta2\|x_i-z\|^2}
			{}={}
				\left(\I-\tfrac{1}{m\beta}\trans{W_i}W_i\right)^{-1} z
			{}={}
				z + \trans{W_i}(m\beta \I - W_i\trans{W_i})^{-1}W_i z,
			\]
			for \(i=1\ldots N\), where the second equality uses the Woodbury identity.
			The Cholesky factors of the \(m_i\times m_i\) matrices \(m\beta \I-W_i\trans{W_i}\), \(i=1,\ldots,N\), can be computed once offline to efficiently solve the linear systems at each \(\bm x\)-update, resulting in \(O\bigl(\sum_{i=1}^Nm_i^2\bigr)\) memory requirement, as opposed to \(O(N^2)=O\bigl(\sum_{i=1}^Nm_i\bigr)^2\) (let alone the operational cost) needed for the original single-agent problem expression.
			
\begin{figure}[p]
				\centering
				\caption{%
					{\rm \S\ref{sec:SPCAconsensus}: consensus sparse PCA problem \eqref{eq:cSPCA} on full datasets.}
					Comparison between \ref{ADMM} (blue) and the L-BFGS enhancement (red), for different number of agents \(N=10,20,50\).
					Left: \ref{ADMM} residual; Right: cost.
					On the \(x\)-axis, the number of linear systems solved (needed for the \(\bm x\)-update): in \ref{ADMM} this coincides with the number of iterations, while in \cref{alg:ADMM} it accounts for each linesearch step.
					This is the only expensive operation, as the \(z\)-update is negligible.
					Apparently, \ref{ADMM} is severly affected by \(N\), whereas using L-BFGS directions in \cref{alg:ADMM} consistently results in faster convergence.%
				}%
				\label{fig:SPCA_consensus}%
				\begin{subfigure}{\linewidth}
					\includetikz[width=0.99\linewidth]{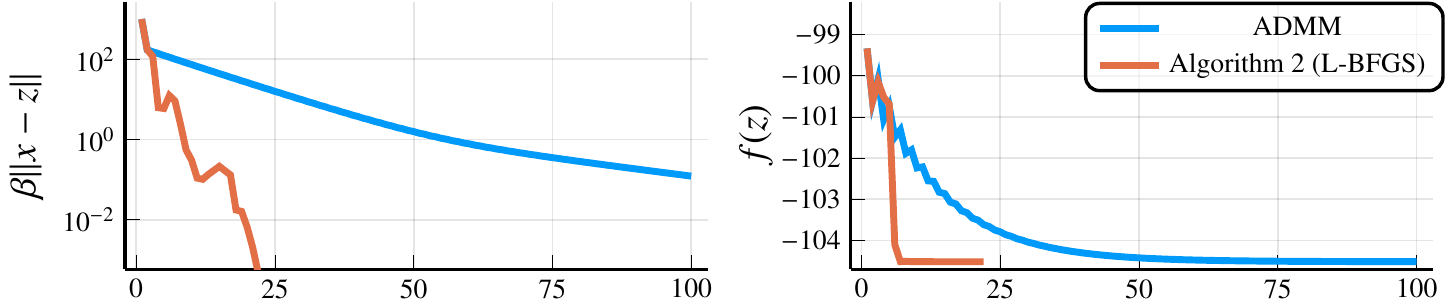}
					
					\includetikz[width=0.99\linewidth]{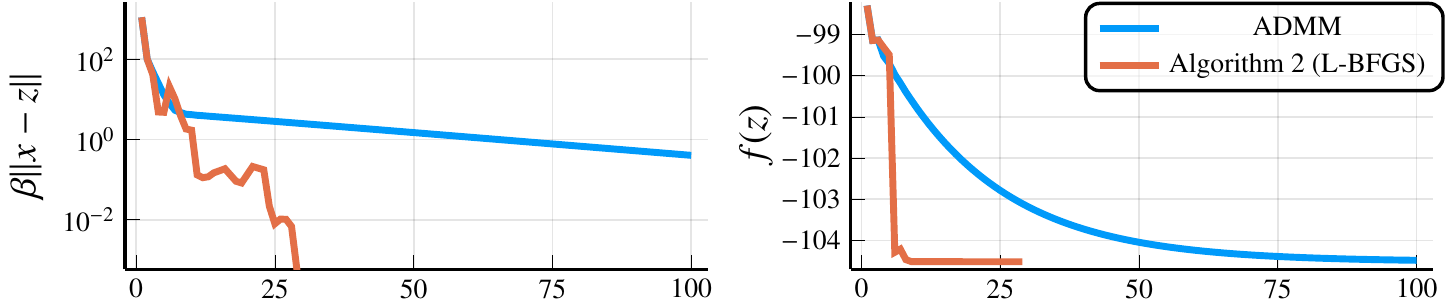}
					
					\includetikz[width=0.99\linewidth]{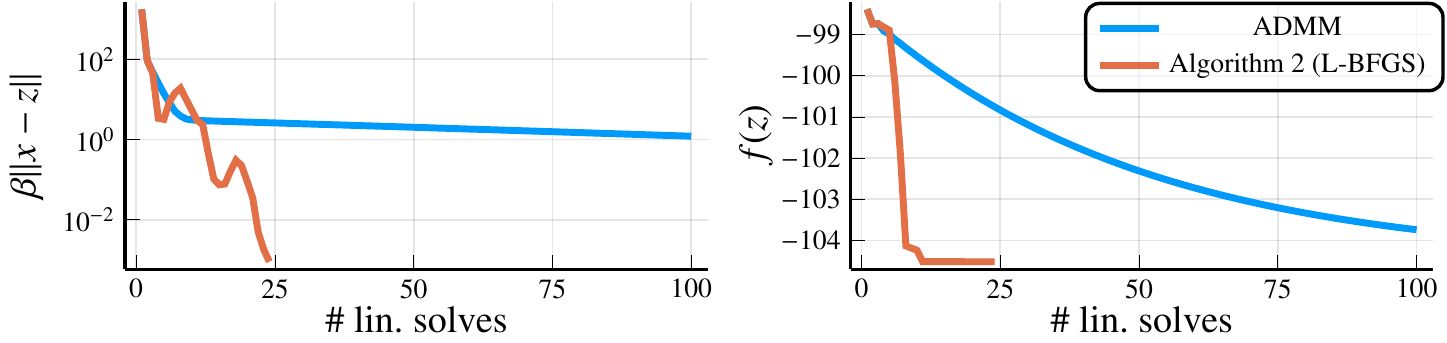}
					\caption{Full \emph{20newsgroup} dataset.}
					\label{fig:SPCA_consensus_20newsgroup}%
				\end{subfigure}
				\begin{subfigure}{\linewidth}
					\includetikz[width=0.99\linewidth]{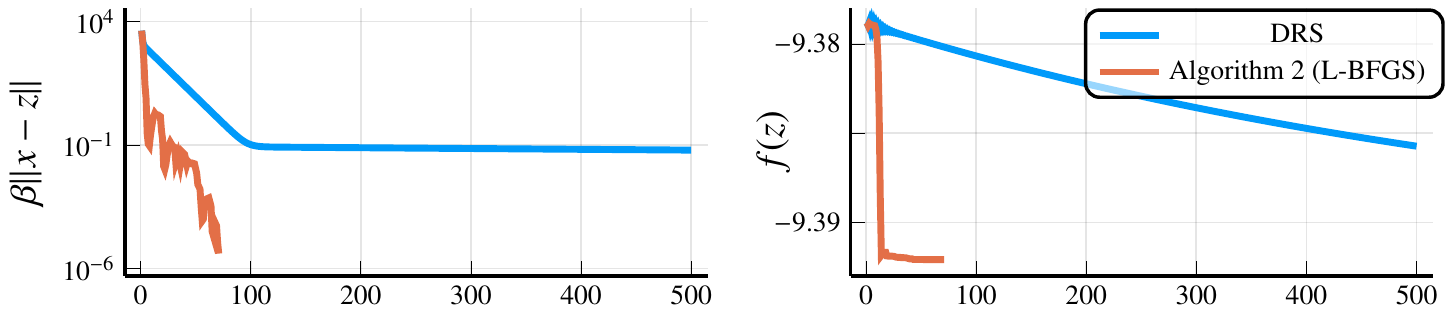}
					
					\includetikz[width=0.99\linewidth]{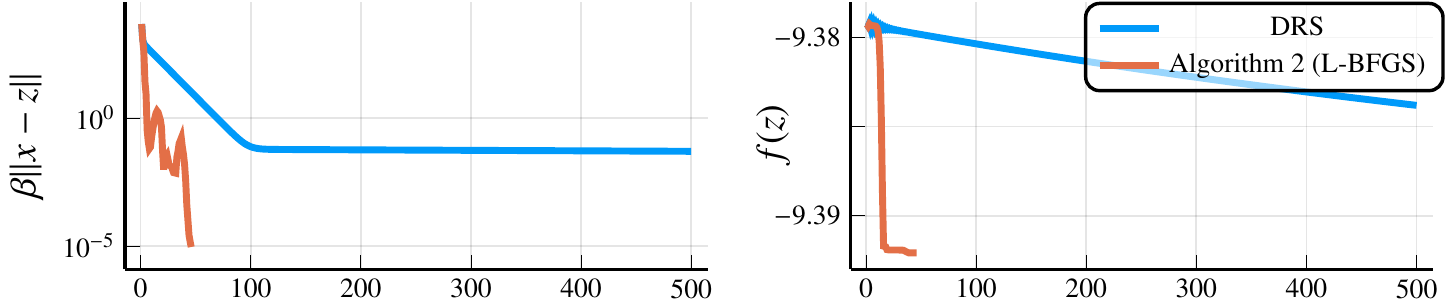}
					
					\includetikz[width=0.99\linewidth]{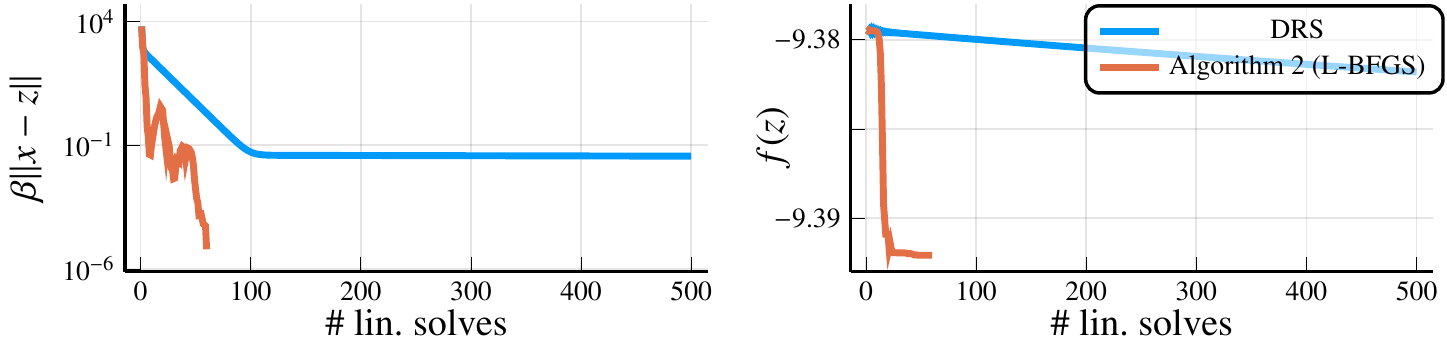}
					\caption{Full \emph{nips\_conference\_papers} dataset.}
					\label{fig:SPCA_consensus_nips}%
				\end{subfigure}
			\end{figure}

			This consensus formulation, however, increases the problem size and thus the ill conditioning, and for moderate values of \(m\), \(n\) and \(N\) the convergence of plain \ref{ADMM} is already prohibitively slow, cf. \cref{fig:SPCA_consensus}.
			On the contrary, the adoption of L-BFGS directions in the \ref{ADMM}-\cref{alg:ADMM} robustifies the performance at the negligible cost of few scalar products per iteration.
			
			\Cref{fig:SPCA_consensus} shows the result of running the proposed algorithm, with L-BFGS directions, to the full \emph{20newsgroup}\footnote{%
				\url{http://www.cad.zju.edu.cn/home/dengcai/Data/TextData.html}
			}
			and \emph{nips\_conference\_papers}\footnote{%
				\url{https://archive.ics.uci.edu/ml/datasets/NIPS+Conference+Papers+1987-2015}
			}
			datasets.
			The former consists of the frequencies of \(n=26214\) words in \(m=18846\) documents, and the latter contains \(n=11463\) word counts from \(m=5811\) NIPS conference papers.
			We split the data in \(N\) subsets (by row) of approximately equal size, for \(N\in\set{10,20,50}\), to put the problem in the form \eqref{eq:cSPCA}.
			In each experiment we used \(\lambda = 1\), the penalty parameter in both \cref{alg:ADMM} and the nominal \ref{ADMM} was set to \(\beta=\tfrac{2}{0.95} L_{\epicomp A\sADMM}\), and \(c=\frac12C(\beta^{-1}L_{\epicomp A\sADMM}, \lambda)\), cf. \eqref{eq:DRE:C} and \eqref{eq:ADMM:beta}.
			In this case, we only considered L-BFGS directions with memory \(5\) since their computation scales better with the problem dimension (cf. \cref{rem:directions}).
			Both algorithms were started at the same initial iterates \(y^0=0\) and \(z^0=(1/n,\ldots,1/n) \in \R^n\).
			Apparently, \cref{alg:ADMM} using L-BFGS converges faster than the nominal \ref{ADMM}, and its convergence speed is significantly less susceptible to the number \(N\) of agents.

		\subsection{Linear model predictive control (strongly convex)}\label{sec:MPC}

			To showcase the performance of the proposed method in the strongly convex case,
			we apply \cref{alg:DRS} to linear model predictive control (MPC) problems \cite{garcia1989model},
			i.e. finite-horizon, discrete-time, linear optimal control problems of the form
			\begin{align}\label{eq:linearMPC}\begin{split}
				\minimize_{\mathclap{\substack{u_0, \ldots, u_{N-1} \in \R^{n_u} \\ x_1, \ldots, x_N \in \R^{n_x}}}}{}
			~&
				\sum_{i=0}^{N-1} \left(\|x_{i+1} - \bar x\|_Q^2 + \|u_i\|_R^2 + g(u_i, x_{i+1})\right)
			\\
				\stt{}
			~&
				\fillwidthof[r]{x_{i+1}}{x_0} \in \R^{n_x} \text{ given}
			\\
			&
				x_{i+1} = A x_i + B u_i,
				~
				i = 0,\dots,N-1.
			\end{split}\end{align}
			The decision variables are the system inputs $u_i$ and states $x_{i+1}$, $i=0,\ldots, N-1$.
			In the quadratic cost terms, $Q \succ 0$ and $R \succ 0$.
			The objective is find the optimal sequence of inputs $u_i$ that drive the system towards the reference state $\bar x$.
			The equality constraints enforce the linear dynamics, while $g$ is a convex functions that models constraints on inputs and states.
			
			Grouping all variables into $z = (u_0, \ldots, u_{N-1}, x_1, \ldots, x_N) \in \R^{N(n_u + n_x)}$,
			and denoting the set $\mathcal{D}_{x_0} = \set{z \in \R^{N(n_u + n_x)}}[x_{i+1} = Ax_i + Bu_i, ~i = 0,\dots,N-1]$,
			problem \eqref{eq:linearMPC} can be solved with \cref{alg:DRS} by setting
			\begin{align*}
				\sDRS(z) {}&{}= \sum_{i=0}^{N-1} \left(\|x_{i+1} - \bar x\|_{Q}^2 + \|u_i\|_{R}^2 + \indicator_{\mathcal D_{x_0}}(s)\right),\\
				\nsDRS(z) {}&{}= \sum_{i=0}^{N-1} g(u_i, x_{i+1}).
			\end{align*}
			Both terms are nonsmooth in general; however, $\sDRS$ is strongly convex and computing $\prox_{\gamma \sDRS}$ amounts to solving a strongly convex, equality-constrained quadratic program: due to the problem structure, this can be done efficiently via the Riccati equation, cf. \cite[Sec. 1.9]{bertsekas1999nonlinear}.
			Therefore, \cref{alg:DRS} can be applied under \cref{ass:DRS*}.
			Notice further that \(\mathcal D_{x_0}\) is an affine subspace, resulting in function \(\sDRS\) being \emph{generalized} quadratic; the linearity of \(\prox_{\gamma\sDRS}\) can be thus conveniently exploited in the linesearch as described in \cref{rem:linearprox}.
			
			\begin{figure}[p]
				\caption{%
					{\rm\S\ref{sec:MPC}: linear MPC problem \eqref{eq:linearMPC} for the AFTI-16 system.}
					Comparison between \ref{DRS} and the linesearch variant \cref{alg:DRS} using modified Broyden, L-BFGS, and Nesterov acceleration directions, to reach a tolerance of $10^{-5}$ at each time step.
				}%
				\begin{subfigure}[t]{\linewidth}%
					\vspace*{0pt}%
					\includetikz[width=\linewidth]{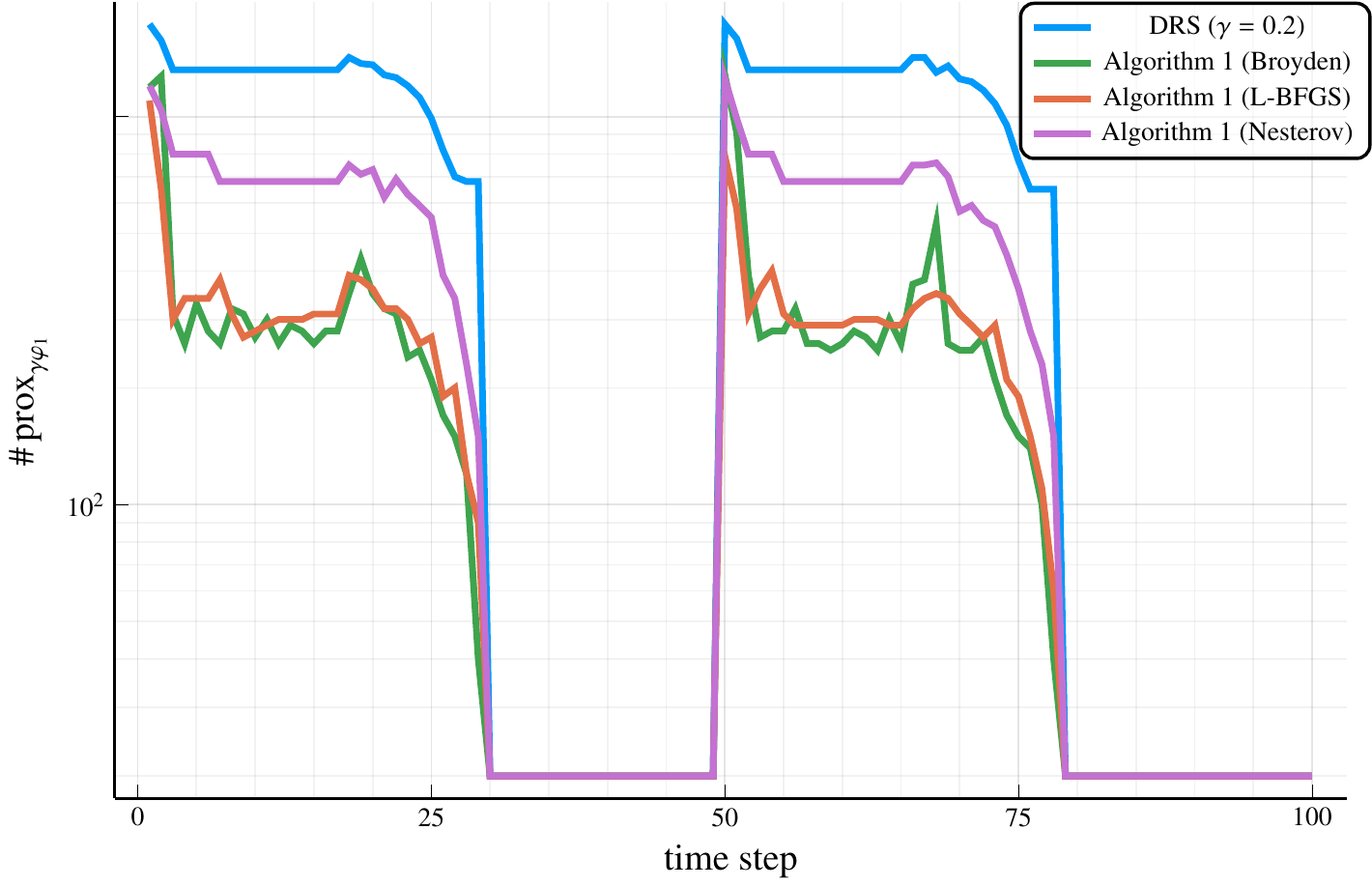}%
					\caption{%
						Number of calls to $\prox_{\gamma\sDRS}$ required by each method.
			% 			throughout the simulation of the AFTI-16 system, to reach a tolerance of $10^{-5}$ at each time step.
						After initially harder problems, the task becomes easier as the system reaches the imposed reference state;
						when the reference state changes at step $50$, the problem becomes again harder to solve.
					}%
					\label{fig:MPC_prox_count}%
					\vspace*{0pt}%
				\end{subfigure}
				\begin{subfigure}[t]{\linewidth}%
					\vspace*{0pt}%
					\includetikz[width=\linewidth]{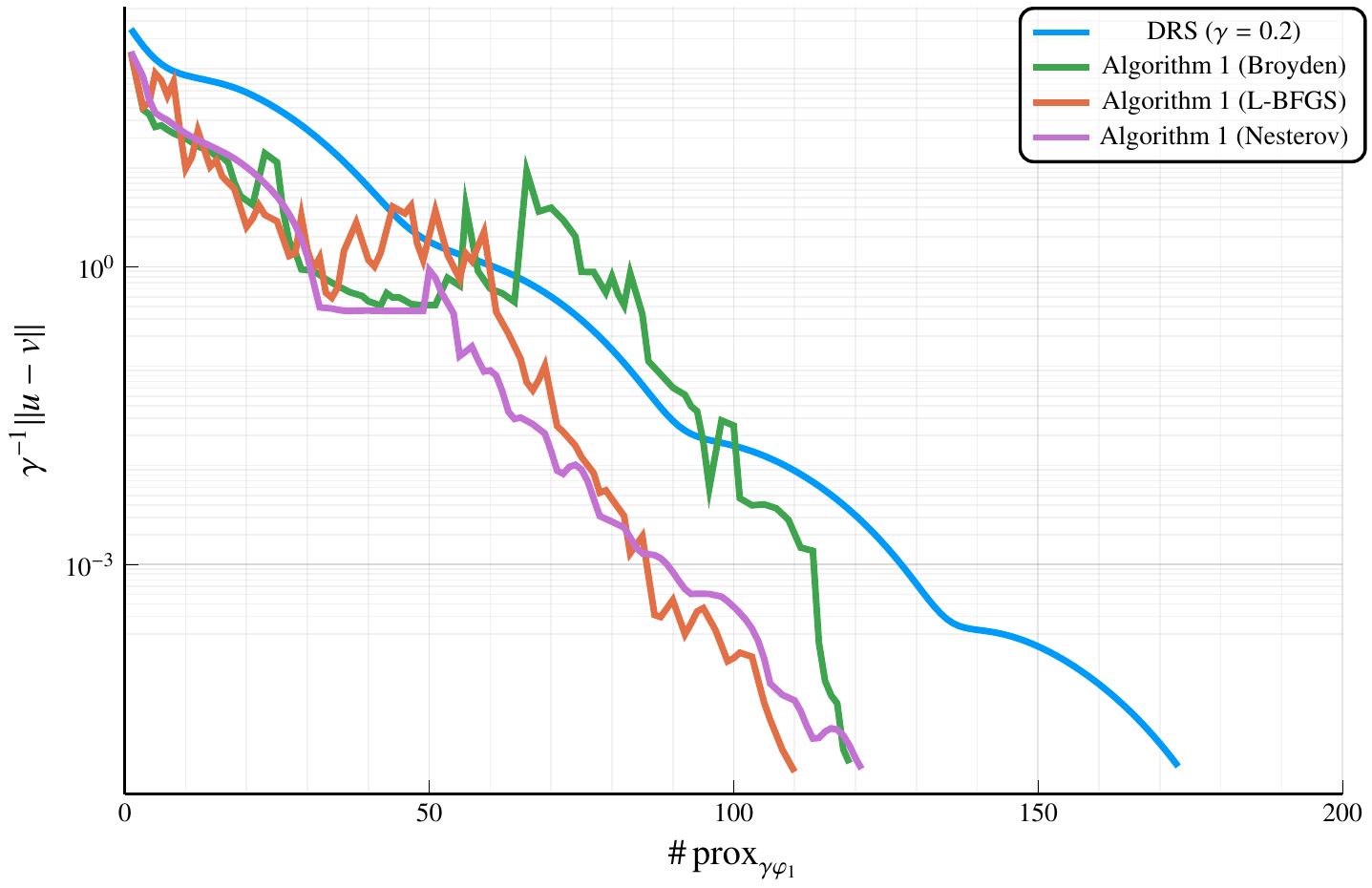}%
					\caption{%
						Convergence of the fixed point residual for the different methods, on the initial problem of the simulation.%
					}%
					\label{fig:MPC_first}%
					\vspace*{0pt}%
				\end{subfigure}%
			\end{figure}
			
			As specific instance of the problem, we considered the AFTI-16 system \cite{bemporad1997nonlinear}:
			this has $n_u = 2$ inputs, $n_x = 4$ states, and unstable dynamics.
			We imposed hard constraint on the input variables and soft constraints on the system states, with
			\[
				g(u_i, x_{i+1}) = \indicator_{[-25, 25]}(u_i) + 10^6 \left[\max\set{0, |x_{i+1}^{(2)}| - 0.5} + \max\set{0, |x_{i+1}^{(4)}| - 100}\right],
			\]
			and set $Q = \diag(10^{-4}, 10^2, 10^{-3}, 10^2)$, $R = \diag(10^{-2}, 10^{-2})$.
			The prediction horizon was set to $N = 10$, and the problem was scaled so as to have identity Hessian in the quadratic cost: this is known to improve significantly the convegence speed of \ref{DRS} \cite{rey2016admm} and proved beneficial for \cref{alg:DRS} as well in our experiments. Therefore the scaled problem has \(\mu_{\sDRS} = 1\).
			In applying \cref{alg:DRS}, we used
			\(\lambda=1\),
			\(\gamma=(0.95\cdot\mu_{\sDRS})^{-1}\),
			\(c=\tfrac{1}{2}C((\gamma \mu_{\sDRS})^{-1}, \lambda)\), cf. \cref{thm:DRE:SD*}.
			For \ref{DRS}, since the problem is convex, it is known that any $\gamma > 0$ is feasible.
			The choice of $\gamma$ that performs best has been made in hindsight; it is well known that the performance of \ref{DRS} is sensitive to this choice, and we are not aware of any generically applicable rule.
			Therefore, after inspecting a grid of values for $\gamma$, we empirically found $\gamma = 0.2$ to give the best performance with \ref{DRS} on this specific problem, and used this value as baseline.
			
			We simulated the system for $100$ time steps, which correspond to $5$ seconds when the original continuous-time dynamics is discretized with a step of $0.05$ seconds.
			At each time step, problem \eqref{eq:linearMPC} is solved with tolerance $10^{-5}$; then, the first optimal input $u_0$ is applied and the system evolves to the next time step, and the next problem is solved.
			On each problem, both \ref{DRS} and \cref{alg:DRS} were warm-started by providing the final $s$-iterate to the previous problem as initial $s$-iterate: this proved beneficial for all algorithms.
			The initial system state is set to $x_0 = (0, 0, 0, 0)$ at the beginning of the simulation; the reference state was set to $\bar x = (0, 0, 0, 10)$ for the first $50$ time steps ($2.5$ seconds), and to $\bar x = (0, 0, 0, 0)$ for the remaining steps.
			
			The performance of \cref{alg:DRS} is illustrated in \cref{fig:MPC_prox_count}, where the total number of calls to $\prox_{\gamma\sDRS}$ at each time step is reported, for all considered algorithms. \cref{fig:MPC_first} shows the convergence of the methods for the very first problem in the simulation.
			While using Anderson acceleration directions did not perform well in this example, it is clear that Broyden, L-BFGS, and Nesterov directions perform significantly better than vanilla \ref{DRS}.

	\section{Conclusions}\label{sec:Conclusions}

		We proposed two linesearch algorithms that allow the employment of Newton-like update directions to enhance \ref{DRS} and \ref{ADMM}.
		The choice of quasi-Newton directions maintains the same low complexity as the original \ref{DRS} and \ref{ADMM} algorithms, as it prescribes only additional direct linear algebra.
		Simulations confirm that L-BFGS considerably robustifies the convergence, rendering these first-order algorithms extremely fast and unaffected by problem size and ill conditioning.
		The proposed algorithms are tuning-free and out-of-the-box, as the needed stepsizes and parameters can adaptively be retrieved without prior knowledge.
		Last but not least, they are suited for fully nonconvex problems, and maintain the same worst-case convergence properties of \ref{DRS} and \ref{ADMM}.

	% ~~~~~~~~~~~~~~~~~~~~~~~~~~~~~~~~~~~~~~~~~~~~~~~~~~~~~~~~~~~~~~~~~~~~~~~~~~~ %

	\begin{appendix}
		\section{Auxiliary results}\label{sec:Appendix}%

			This appendix contains some auxiliary results needed for the convergence analysis of \cref{sec:Convergence}.
			As shown in \cite[Eq. (3.4)]{themelis2020douglas}, the DRE can be expressed in terms of the forward-backward envelope \(\FBE\) \cite{patrinos2013proximal,stella2017forward,themelis2018forward} as
			\begin{subequations}
				\begin{align}
				\label{eq:FBEmin}
					\DRE(s)
				{}={} &
					\FBE(u)
				{}\coloneqq{}
					\min_{w\in\R^p}\set{
						\sDRS(u)
						{}+{}
						\nsDRS(w)
						{}+{}
						\innprod{\nabla\sDRS(u)}{w-u}
						{}+{}
						\tfrac{1}{2\gamma}\|w-u\|^2,
					}
				\shortintertext{%
					where \(u=\prox_{\gamma\sDRS}(s)\) and the minimum is attained at any \(v\in\prox_{\gamma\nsDRS}(2u-s)\).
					Equivalently,
				}
				\label{eq:FBEMoreau}
					\DRE(s)
				{}={} &
					\sDRS(u)
					{}-{}
					\tfrac\gamma2\|\nabla\sDRS(u)\|^2
					{}+{}
					\nsDRS^\gamma(\Fw u).
				\end{align}
			\end{subequations}
			
			\begin{fact}[{\cite[Prop. 3.3]{themelis2020douglas}}]\label{thm:sandwich}%
				Suppose that \cref{ass:DRS} holds and let \(\gamma<\nicefrac{1}{L_{\sDRS}}\) be fixed.
				Then, for all \(s\in\R^p\) and \((u,v)\in\DRS(s)\) it holds that
				\begin{qedequation*}
					\varphi(v)
					{}+{}
					\tfrac{1-\gamma L_{\sDRS}}{2\gamma}\|v-u\|^2
				{}\leq{}
					\DRE(s)
				{}\leq{}
					\varphi(u).
				\end{qedequation*}
			\end{fact}
			
			\grayout{%
				\begin{lem}[Equivalence of (strong) local minimality]\label{thm:DREequiv:strmin}%
				% 	Suppose that \cref{ass:DRS} holds, and let \(\bar s\) be a fixed point for \ref{DRS} with stepsize \(\bar\gamma\leq\nicefrac{1}{L_{\sDRS}}\).
				% 	Then, for every \(\gamma<\bar\gamma\) \(\bar s\) is a fixed point for \ref{DRS} with stepsize \(\gamma\), and is a strong local minimum for \(\DRE\) whenver \(\bar u\coloneqq\prox_{\gamma\sDRS}(\bar s)\) is a strong local minimum for \(\varphi\).
					Suppose that \cref{ass:DRS} holds, and let \(\bar s\) be a fixed point for \ref{DRS} with stepsize \(\gamma<\nicefrac{1}{L_{\sDRS}}\).
					If \(\bar u\coloneqq\prox_{\gamma\sDRS}(\bar s)\) is a strong local minimum for \(\varphi\) and the inclusion
					\(
						\bar u
					{}\in{}
						\prox_{\gamma\sDRS}(2\bar u-\bar s)
					\)
					holds as equality (namely,
					\(
						\prox_{\gamma\sDRS}(2\bar u-\bar s)
					{}={}
						\set{\bar u}
					\)),
					then \(\bar s\) is a strong local minimum for \(\DRE\).
					\begin{proof}
						By assumption, there exists \(\mu>0\) such that
						\(
							\varphi(u)
						{}\geq{}
							\varphi(\bar u)
							{}+{}
							\tfrac\mu2\|u-\bar u\|^2
						\)
						for all \(u\) sufficiently close to \(\bar u\).
						Let
						\(
							\delta
						{}\coloneqq{}
							\min\set{
								\tfrac\mu2,
								\tfrac{1-\gamma L_{\sDRS}}{2\gamma}
							}
						{}>{}
							0
						\).
						To arrive at a contradiction, suppose that for all \(k\in\N\) there exists \(s^k\in\ball{\bar s}{\nicefrac1k}\) such that
						\(
							\DRE(s^k)
						{}<{}
							\DRE(\bar s)
							{}+{}
							\tfrac\delta2\|s^k-\bar s\|^2
						\).
						Let \(u^k\coloneqq\prox_{\gamma\sDRS}(s^k)\) and \(v^k\in\prox_{\gamma\nsDRS}(2u^k-s^k)\).
						Then, since \(\prox_{\gamma\sDRS}\) is continuous and \(\prox_{\gamma\nsDRS}\) is outer semicontinuous \(u^k\to\bar u\) and \(v^k\to\bar u\), where the second limit owes to the fact that \(\prox_{\gamma\nsDRS}(2\bar u-\bar s)=\set{\bar u}\).
						\def\lb{\nicefrac1\gamma-L_{\sDRS}}
						We have
						\begin{align*}
							\varphi(v^k)
						{}\overrel[\leq]{\ref{thm:sandwich}}{} &
							\DRE(s^k)
							{}-{}
							\tfrac{1-\gamma L_{\sDRS}}{2\gamma}\|v^k-u^k\|^2
						\\
						{}<{} &
							\DRE(\bar s)
							{}+{}
							\tfrac\delta2\|s^k-\bar s\|^2
							{}-{}
							\tfrac{1-\gamma L_{\sDRS}}{2\gamma}\|v^k-u^k\|^2
						\\
						{}={} &
							\varphi(\bar x)
							{}+{}
							\tfrac\delta2\|x^k-\bar x\|^2
							{}-{}
							\tfrac{\lb}{2}\|x^k-\bar x^k\|^2.
						\shortintertext{%
							By using the inequality
							\(
								\tfrac12\|a-c\|^2
							{}\leq{}
								\|a-b\|^2
								{}+{}
								\|b-c\|^2
							\)
							holding for all vectors \(a,b,c\in\R^n\), we have
						}
							\varphi(\bar x^k)
						{}<{} &
							\varphi(\bar x)
							{}+{}
							\delta\|\bar x^k-\bar x\|^2
							{}+{}
							\bigl(\delta-\tfrac{\lb}{2}\bigr)\|x^k-\bar x^k\|^2
						\\
						{}\leq{} &
							\varphi(\bar x)
							{}+{}
							\tfrac\mu2\|\bar x^k-\bar x\|^2,
						\end{align*}
						where the last inequality follows from the definition of \(\delta\).
				% 		Thus, we obtain
				% 		\(
				% 			\varphi(\bar x^k)
				% 		{}<{}
				% 			\varphi(\bar x)
				% 			{}+{}
				% 			\tfrac\mu2\|\bar x^k-\bar x\|^2
				% 		\)
				% 		for all \(k\in\N\), hence the contradiction since \(\bar x^k\) is arbitrarily close to \(\bar x^k\).
					\end{proof}
				\end{lem}
			
				The previous result appears in \cite[Thm. 6.5]{themelis2018proximal} as a byproduct of a more general analysis.
				The interested reader is referred to Section 3.2.2 therein for a discussion regarding the necessity of single-valuedness of \(\prox_{\gamma\nsDRS}(2u_\star-s_\star)\).
				A more general version of the following result is also discussed in \cite[Prop. 3.9]{themelis2018proximal}.
				We conclude with another useful property that relates the DRE and the original function \(\varphi\).
			
				We now show that for small enough stepsizes \(\gamma\) the DRE decreases after one iteration of \ref{DRS}.
				Although the result can be improved (see \cite{themelis2020douglas} for tight bounds on \(\gamma\)), we sacrifice generality for the sake of simplicity.
			}%
			
			\begin{lem}\label{thm:QG}%
				Suppose that \cref{ass:DRS} holds.
				Then, for all \(s,\bar u\in\R^p\) and \(\gamma<\nicefrac{1}{L_{\sDRS}}\)
				\[
					\DRE(s)-\varphi(\bar u)
				{}\leq{}
					\tfrac{1+\gamma L_{\sDRS}}{2\gamma}
					\|\prox_{\gamma\sDRS}(s)-\bar u\|^2.
				\]
				\grayout{%
					In particular,
					\[
						\DRE(s)-\min\varphi
					{}\leq{}
						\tfrac{1+\gamma L_{\sDRS}}{2\gamma}
						\dist\bigl(\prox_{\gamma\sDRS}(s),\,\argmin\varphi\bigr)^2.
					\]
				}%
				\begin{proof}
					Let \(u\coloneqq\prox_{\gamma\sDRS}(s)\) for brevity.
					By plugging \(w=\bar u\) into \eqref{eq:FBEmin} we obtain
					\begin{align*}
						\DRE(s)
					{}\leq{} &
						\nsDRS(\bar u)
						{}+{}
						\underbracket[0.5pt]{
							\sDRS(u)
							{}+{}
							\innprod{\nabla\sDRS(u)}{\bar u-u}
						}{}
						{}+{}
						\tfrac{1}{2\gamma}\|\bar u-u\|^2
					\\
					{}\leq{} &
						\nsDRS(\bar u)
						{}+{}
						\overbracket[0.5pt]{
							\fillwidthof[l]{
								\sDRS(u)
								{}+{}
								\innprod{\nabla\sDRS(u)}{\bar u-u}
							}{
								\sDRS(\bar u)
								{}+{}
								\tfrac{L_{\sDRS}}{2}\|\bar u-u\|^2
							}
						}{}
						{}+{}
						\tfrac{1}{2\gamma}\|\bar u-u\|^2,
					\end{align*}
					where the second inequality uses the known quadratic upper bound \cite[Prop. A.24]{bertsekas1999nonlinear} for functions with Lipschitz-continuous gradient.
				\end{proof}
			\end{lem}
			
			\begin{lem}
				Suppose that \cref{ass:DRS*} holds and let \(\gamma>1\nicefrac{}{\mu_{\sDRS}}\) be fixed.
				Then,
				\begin{equation}\label{eq:strongDuality}
					\inf\varphi
				{}={}
					-\inf\psi
				{}={}
					-\inf\DRE*
				{}={}
					\sup\DRE
				\end{equation}
				where \(\gamma_*=\nicefrac1\gamma\).
				Moreover, for any \(s\in\R^p\) it holds that
				\begin{align}\label{eq:Qlb}
					\tfrac{1}{2\gamma}\|x_\star-v\|^2
					{}+{}
					\tfrac{\gamma\mu_{\sDRS}-1}{2\gamma}\|x_\star-u\|^2
				{}\leq{}
					\inf\varphi
					{}-{}
					\DRE(s)
				{}={} &
					\DRE*(s_*)
					{}-{}
					\inf\psi
			% 	\\
			% 	{}\leq{} &
			% 		\psi(v_*)
			% 		{}-{}
			% 		\inf\psi
			% 		{}+{}
			% 		\tfrac{1+\gamma\mu_{\sDRS}}{2\mu_{\sDRS}}\|u_*-v_*\|^2,
				\end{align}
				where \(x_\star\) is the unique minimizer of \(\varphi\), \(s_*=-\nicefrac s\gamma\), and \((u,v)=\DRS(s)\).
				\begin{proof}
					Due to strong convexity, the set of primal solutions \(\argmin\varphi\) is a singleton, ensuring strong duality \(\inf\varphi=-\inf\psi\) through \cite[Thm. 5.2.1(b)-(c)]{auslender2002asymptotic}.
					Since \(\sDRS*\) is \(\nicefrac{1}{\mu_{\sDRS}}\)-smooth, it follows from \cref{thm:DREequiv:inf} that \(\inf\DRE*=\inf\psi\) for every \(\gamma_*<\mu_{\sDRS}\); combined with the identity \(\DRE(s)=-\DRE*(-s\nicefrac{}\gamma)\) holding for \(\gamma_*=1\nicefrac{}\gamma\) (cf. \cref{thm:selfdual}), \eqref{eq:strongDuality} is obtained.
					Let now \(s\in\R^p\) be fixed and consider \((u,v)=\DRS(s)\).
					From the inclusion \(\tfrac{s-u}{\gamma}\in\partial\sDRS(u)\) (cf. \eqref{eq:prox:subdiff}) and strong convexity of \(\sDRS\) one has
					\begin{align*}
						\sDRS(x_\star)
					{}\geq{} &
						\sDRS(u)
						{}+{}
						\tfrac1\gamma\innprod{s-u}{x_\star-u}
						{}+{}
						\tfrac{\mu_{\sDRS}}{2}\|x_\star-u\|^2.
					\shortintertext{%
						Similarly, since \(\tfrac{2u-s-v}{\gamma}\in\partial\nsDRS(v)\) and \(\nsDRS\) is convex, one has
					}
						\nsDRS(x_\star)
					{}\geq{} &
						\nsDRS(v)
						{}+{}
						\tfrac1\gamma\innprod{2u-s-v}{x_\star-v}.
					\end{align*}
					Summing the two inequalities yields
					\begin{align*}
						\inf\varphi
					{}\geq{} &
						\sDRS(u)
						{}+{}
						\nsDRS(v)
						{}+{}
						\tfrac1\gamma\innprod{s-u}{x_\star-u}
						{}+{}
						\tfrac1\gamma\innprod{2u-s-v}{x_\star-v}
						{}+{}
						\tfrac{\mu_{\sDRS}}{2}\|x_\star-u\|^2
					\\
					{}={} &
						\sDRS(u)
						{}+{}
						\nsDRS(v)
						{}+{}
						\tfrac1\gamma\innprod{s-u}{v-u}
						{}+{}
						\tfrac1\gamma\innprod{u-v}{x_\star-v}
						{}+{}
						\tfrac{\mu_{\sDRS}}{2}\|x_\star-u\|^2
					\\
					{}={} &
						\sDRS(u)
						{}+{}
						\nsDRS(v)
						{}+{}
						\tfrac1\gamma\innprod{s-u}{v-u}
						{}+{}
						\tfrac{1}{2\gamma}\|u-v\|^2
						{}+{}
						\tfrac{1}{2\gamma}\|x_\star-v\|^2
						{}+{}
						\tfrac{\gamma\mu_{\sDRS}-1}{2\gamma}\|x_\star-u\|^2
					\\
					{}={} &
						\DRE(s)
						{}+{}
						\tfrac{1}{2\gamma}\|x_\star-v\|^2
						{}+{}
						\tfrac{\gamma\mu_{\sDRS}-1}{2\gamma}\|x_\star-u\|^2.
					\end{align*}
					The claim now follows from the identity \(\DRE*(s_*)=-\DRE(s)\) shown in \cref{thm:selfdual}.
				\end{proof}
			\end{lem}
			
	\end{appendix}

	% ~~~~~~~~~~~~~~~~~~~~~~~~~~~~~~~~~~~~~~~~~~~~~~~~~~~~~~~~~~~~~~~~~~~~~~~~~~~ %

	\bibliographystyle{plain}
	\bibliography{Bibliography.bib}

	% ~~~~~~~~~~~~~~~~~~~~~~~~~~~~~~~~~~~~~~~~~~~~~~~~~~~~~~~~~~~~~~~~~~~~~~~~~~~ %

\end{document}